\documentclass[10pt,openbib]{article}
\usepackage{amsmath}
\usepackage{amsfonts}
\usepackage{amssymb}
\usepackage{graphicx}
\newcommand{\R}{\mathbb{R}} 
\newcommand{\C}{\mathbb{C}}
 
\newcommand{\CP}{\mathbb{C}\mathbb{P}}

\newcommand{\Z}{\mathbb{Z}}

\newcommand{\dee}{\mathop{\! \,\mathrm{d} \!}\nolimits}
\newcommand{\comp}{\raisebox{0pt}{$\scriptstyle\circ \, $}}
\newcommand{\setrule}{\, \rule[-4pt]{.5pt}{13pt}\, }
\newcommand{\smallspace}{\smallskip\par\noindent}

\newcommand{\onehalf}{\mbox{$\frac{\scriptstyle 1}{\scriptstyle 2}\,$}} 
 
\newcommand{\lefthook}{\mbox{$\, \rule{8pt}{.5pt}\rule{.5pt}{6pt}\, \, $}}
\newcommand{\ttfrac}[2]{\mbox{$\frac{{\scriptstyle #1}}{{\scriptstyle #2}}$}}
\newcommand{\circdot}{\raisebox{2pt}{\tiny$\, \bigodot$}}

\allowdisplaybreaks

\begin{document}
\thispagestyle{empty}
\begin{center}
{\Large \textbf{An affine model of a Riemann surface associated \\ 
\mbox{}\vspace{-.2in} \\
to a Schwarz-Christoffel mapping}}
\mbox{}\vspace{.05in} \\ 
\mbox{}\\
Richard Cushman\footnotemark 
\end{center}
\footnotetext{Department of Mathematics and Statistics, University of Calgary \\ 
\rule{.22in}{0in}e-mail: rcushman@ucalgary.ca  \newline
\rule{.22in}{0in}printed: \today
}

\begin{abstract}
In this paper we construct an affine model of a Riemann surface with a flat Riemannian metric associated to a Schwarz-Christoffel mapping of the upper half plane onto a rational triangle. We explain the relation between the geodesics on this Riemann surface and billiard motions in a regular stellated $n$-gon in 
the complex plane.
\end{abstract}

\section{Introduction}

Here we give a detailed description of the contents of this paper. \medskip

Consider the conformal Schwarz-Christoffel mapping 
\begin{align}
& F_T: {\C }^{+} \rightarrow T \subseteq \C : \xi \mapsto \int^{\xi }_{0} \frac{\dee \xi }{\eta }. 
\tag*{$(\mathrm{I}\, 1 )$}
\end{align}
where 
\begin{align}
&{\eta }^n = {\xi }^{n-n_0}(1-\xi )^{n-n_1}. 
\tag*{$(\mathrm{I}\, 2 )$}
\end{align}
The map $F_T$ sends the closed upper half plane ${\C }^{+}$ onto the rational triangle $T = T_{n_0n_1n_{\infty}}$, 
where $n_0+n_1+n_{\infty} = n$ and $1 \le n_0\le n_1\le n_{\infty}$. Because ${F_T}{|[0,1]}$ has real 
values, using the Schwarz reflection principal we extend $F_T$ to the conformal map 
\begin{displaymath}
F_Q: \C \setminus \{ 0,1\} \rightarrow Q = T \cup \overline{T}
\end{displaymath}
of $\C \setminus \{ 0,1\}$ onto the quadrilateral $Q$. \medskip 

Following Aurell and Itzykson \cite{aurell-itzykson} we associate to the map $F_Q$ the affine Riemann surface 
$\mathcal{S} \subseteq {\C }^2$ defined by ($\mathrm{I}\, 2$). Then ${\mathcal{S}}_{\mathrm{reg}} = 
\mathcal{S} \setminus \{ (0,0), (1,0) \}$ is a smooth submanifold of ${\C}^2 \setminus \{ \eta = 0 \}$. 
To determine the geometry of ${\mathcal{S}}_{\mathrm{reg}}$, we think of $\mathcal{S}$ as the 
$n$-fold branched covering $\pi : \mathcal{S} \subseteq {\C }^2 \rightarrow \C : (\xi ,\eta ) \mapsto \xi $.  
The map $\pi $ has branch points at $0$, $1$, and $\infty$ of degree $\frac{n}{d_0}$, $\frac{n}{d_1}$, and 
$\frac{n}{d_{\infty}}$, respectively, where $d_j = \gcd (n, n_j )$ for $j = 0,1, \infty$. Using the 
Riemann-Hurwitz formula, see McKean and Moll \cite{mckean-moll}, it follows that the genus of the compact 
Riemann surface $\mathrm{cl}(\mathcal{S}) \subseteq {\mathbb{CP}}^2$ is 
$\onehalf (n+2 -(d_0+d_1+d_{\infty}))$.  Here $\mathrm{cl}$ denotes closure. Thus 
${\mathcal{S}}_{\mathrm{reg}}$, which is $\mathrm{cl}(\mathcal{S})$ less 
three points, has the same genus as $\mathrm{cl}(\mathcal{S})$. \medskip 

We now give a more geometric description of ${\mathcal{S}}_{\mathrm{reg}}$. The abelian group 
$\widehat{\mathcal{G}}$ generated by 
\begin{displaymath}
\mathcal{R}: {\mathcal{S}}_{\mathrm{reg}} \subseteq {\C }^2 \rightarrow 
{\mathcal{S}}_{\mathrm{reg}} \subseteq {\C }^2: (\xi , \eta ) \mapsto (\xi , {\mathrm{e}}^{2\pi i/n}\eta )
\end{displaymath}
is the group of covering transformations of the holomorphic covering map 
\begin{align}
& \widehat{\pi }: {\mathcal{S}}_{\mathrm{reg}} \subseteq {\C }^2 \rightarrow \C \setminus \{ 0,1 \}: (\xi ,\eta ) 
\mapsto \xi .
\tag*{$(\mathrm{I}\, 3 )$}
\end{align}
Let $\mathcal{D}$ be a fundamental domain for the $\widehat{\mathcal{G}}$ action on 
${\mathcal{S}}_{\mathrm{reg}}$, which is a ``sheet'' of the covering map $\widehat{\pi }$ ($\mathrm{I}\, 3$). 
Its image under the map 
\begin{align}
{\delta }_Q: \mathcal{D} \subseteq {\mathcal{S}}_{\mathrm{reg}} \rightarrow Q \subseteq \C: 
(\xi , \eta ) \mapsto (F_Q \comp \widehat{\pi })(\xi ,\eta ),
\tag*{$(\mathrm{I}\, 4 )$}
\end{align}
which is a holomorphic diffeomorphism of $\mathrm{int}\, \mathcal{D}$ onto $\mathrm{int}\, Q$ and 
a homeomorphism of $\partial \mathcal{D}$ onto $\partial Q$, is the quadrilateral $Q$. \medskip 

Let $K^{\ast} = \coprod_{0 \le j \le n-1} R^j \delta (Q)$, where $R: \C \rightarrow \C: z \mapsto 
{\mathrm{e}}^{2\pi i/n}z$. Then $K^{\ast }$ is a regular stellated $n$-gon, which is invariant under 
the action of the dihedral group $G$ generated by the rotation $R$ and the reflection $U: \C \rightarrow \C: 
z \mapsto \overline{z}$, that are subject to the relation $RU = UR^{-1}$. 
Using $\mathrm{cl}(K^{\ast})$ we build a model ${\widetilde{\mathcal{S}}}_{\mathrm{reg}}$ of 
the affine Riemann surface ${\mathcal{S}}_{\mathrm{reg}}$ following Richens and Berry 
\cite{richens-berry}. We say that two closed edges $E$ and $E'$ of $\mathrm{cl}(K^{\ast })$ are 
equivalent $\sim $ if they are not adjacent and $E'$ is the reflection in the diagonal 
$R^m{\ell }^j$, where ${\ell }^j = R^{n_j}U\ell $ and $\ell$ is the edge of $Q$ contained in the ray 
$R^{\pi n_0/n}({\R }_{>0})$. The $G$ orbit space formed by first identifying equivalent points of 
$\mathrm{cl}(K^{\ast })$, which are on equivalent edges in $\partial K^{\ast }$ or are points 
in $\mathrm{int}\, \mathrm{cl}(K^{\ast })$, and then acting on the identification space 
$(\mathrm{cl}(K^{\ast })\setminus \{ O \})^{\sim}$ by the induced action of the group $G$ gives 
${\widetilde{\mathcal{S}}}_{\mathrm{reg}}$. Since the action of $G$ on the identification space is 
free and proper, ${\widetilde{\mathcal{S}}}_{\mathrm{reg}}$ is a smooth 
$1$-dimensional complex manifold. Its genus is $\onehalf (n+2-(d_0+d_1+d_{\infty}))$. 
So ${\widetilde{\mathcal{S}}}_{\mathrm{reg}}$ is a model of the affine Riemann surface 
${\mathcal{S}}_{\mathrm{reg}}$. \medskip %

We construct an affine model of ${\widetilde{\mathcal{S}}}_{\mathrm{reg}}$ as follows. Reflecting 
in the edges of $K^{\ast } \setminus \{ O \}$, which is $\mathrm{cl}(K^{\ast })$ less the vertices and center $O$, 
and then in the edges of the reflected $K^{\ast } \setminus \{ O \}$ et cetera , gives 
$\C \setminus {\mathbb{V}}^{+}$, which is certain translations of $K^{\ast } \setminus \{ O \}$ that generate 
the abelian group $\mathcal{T}$. Here ${\mathbb{V}}^{+}$ is the union of translations of the vertices 
of $\mathrm{cl}(K^{\ast })$ and its center $O$ by elements of $\mathcal{T}$. The group 
$\mathfrak{G} = G \ltimes \mathcal{T}$ acts freely, properly, and transitively on the identification 
space $(\C \setminus {\mathbb{V}}^{+})^{\sim}$ of equivalent points, which are either on equivalent edges of $\C \setminus {\mathbb{V}}^{+}$ or lie in the interior of some $\mathcal{T}$ translate of $K^{\ast } \setminus \{ O \}$. The $\mathfrak{G}$ orbit space $(\C \setminus {\mathbb{V}}^{+})^{\sim}/\mathfrak{G}$ 
of the induced action of $\mathfrak{G}$ is holomorphically differomorphic to 
${\widetilde{\mathcal{S}}}_{\mathrm{reg}}$. It is an affine model of ${\mathcal{S}}_{\mathrm{reg}}$ being 
the space of $\mathfrak{G}$ orbits on $\C \setminus {\mathbb{V}}^{+}$, where 
$\mathfrak{G}$ is a discrete subgroup of the $2$-dimensional Euclidean group. \medskip 

We now look at dynamics on the affine Riemann surface ${\mathcal{S}}_{\mathrm{reg}}$. The vector 
\begin{align}
& X(\xi ,\eta ) = \eta \frac{\partial }{\partial \xi } + \ttfrac{n-n_0}{n}\,  {\xi }^{n-n_0-1}(1 - \xi )^{n-n_1-1} 
(1 - \ttfrac{2n-n_0-n_1}{n-n_0} \xi ) \frac{\partial }{\partial \eta }
\tag*{$(\mathrm{I}\, 5 )$}
\end{align}
is tangent to ${\mathcal{S}}_{\mathrm{reg}}$ at every $(\xi , \eta ) \in \mathcal{D}$ and defines a 
nowhere vanishing holomorphic vector field on the fundamental domain $\mathcal{D}$. Since 
$\frac{\partial }{\partial z} = T_{\xi }F_Q\big( \eta \frac{\partial }{\partial \xi } \big)$ for every 
$(\xi ,\eta ) \in \mathcal{D}$ we get $T_{(\xi , \eta )} \delta \, X(\xi , \eta ) = 
\frac{\partial }{\partial z}\rule[-6pt]{.5pt}{14pt}\raisebox{-5pt}{$\, {\scriptscriptstyle z = \delta (\xi , \eta )}$ } $,  
where 
\begin{align}
& \delta : \mathcal{D} \subseteq {\mathcal{S}}_{\mathrm{reg}} \rightarrow Q \subseteq \C : 
(\xi , \eta ) \mapsto (F_Q \comp \widehat{\pi })(\xi , \eta ),  
\tag*{$(\mathrm{I}\, 6 )$}
\end{align}
the map $\delta $ ($\mathrm{I}\, 6 $) straightens the holomorphic vector field $X$ on $\mathcal{D}$. 
Since $\mathcal{D}$ is a connected open subset of ${\mathcal{S}}_{\mathrm{reg}}$, 
the map ${\delta }_Q$ ($\mathrm{I}\, 3$) straightens the holomorphic vector field $X$ on 
${\mathcal{S}}_{\mathrm{reg}}$ determined by $X$ on $\mathcal{D}$. \medskip 

Let $u = \mathrm{Re}\, z$ and $v= \mathrm{Im}\, z$. Then $\gamma = \dee u \circdot \dee u + 
\dee v \circdot \dee v = \dee z \circdot \dee \overline{z}$ is the Euclidean metric on $\C $. Pulling 
${\gamma }_Q = {\gamma }_{|Q}$ back by the map ${\delta }_Q$ ($\mathrm{I}\, 3$) gives a 
Riemannian metric $\Gamma  = \frac{1}{\eta } \dee z \circdot \frac{1}{\overline{\eta }} \dee \overline{z} $ on 
${\mathcal{S}}_{\mathrm{reg}}$. Since the metric ${\gamma }_Q$ is flat on $Q$, the metric 
$\Gamma $ on ${\mathcal{S}}_{\mathrm{reg}}$ is flat. In other words, the map 
${\delta }_Q: ({\mathcal{S}}_{\mathrm{reg}}, \Gamma ) \rightarrow (Q, {\gamma }_Q)$ is an isometry. Thus 
${\delta }_Q$ is a developing map in the sense of differential geometry, see Spivak \cite[note 12, vol. 2]{spivak} 
and Gauss \cite{gauss}. Since the vector field $X$ on ${\mathcal{S}}_{\mathrm{reg}}$ preserves the 
metric $\Gamma $, the vector field $X$ ($\mathrm{I}\, 5$) on ${\mathcal{S}}_{\mathrm{reg}}$ is the 
geodesic vector field for the metric $\Gamma $. However, $X$ is incomplete, since the image 
of a geodesic on ${\mathcal{S}}_{\mathrm{reg}}$ under the map ${\delta }_Q$ is a straight line on $Q$, which is parallel to the $u$ axis on $\C $, that runs off $Q$ in finite time. The group $\mathcal{G}$ generated by the mappings  $\mathcal{R}: {\mathcal{S}}_{\mathrm{reg}} \rightarrow {\mathcal{S}}_{\mathrm{reg}}: 
(\xi ,\eta ) \mapsto (\xi , {\mathrm{e}}^{2\pi i/n}\eta ) $ and 
$\mathcal{U}: {\mathcal{S}}_{\mathrm{reg}} \rightarrow {\mathcal{S}}_{\mathrm{reg}}: 
(\xi ,\eta ) \mapsto (\overline{\xi }, \overline{\eta })$ 
preserves the metric $\Gamma $. The map ${\delta }_Q $ ($\mathrm{I}\, 3$) extends to the developing map 
\begin{align}
& {\delta }_{K^{\ast }}: ({\mathcal{S}}_{\mathrm{reg}}, \Gamma ) \rightarrow (K^{\ast }, {\gamma }_{K^{\ast }}), 
\tag*{$( \mathrm{I}\, 7 )$}
\end{align}
which is an isometry that intertwines the action of $\mathcal{G}$ on ${\mathcal{S}}_{\mathrm{reg}}$ 
with the action of $G$ on $K^{\ast }$. Since the geodesic vector field $X$ on ${\mathcal{S}}_{\mathrm{reg}}$ 
is invariant under the action of $\mathcal{G}$ and the vector field $\frac{\partial }{\partial z}$ on 
$K^{\ast }$ is invariant under the action of $G$, the map ${\delta }_{K^{\ast }}$ sends geodesics 
on ${\mathcal{S}}_{\mathrm{reg}}$ to geodesics on $K^{\ast }$. However, incompleteness of the vector field $X$ 
remains. \medskip 

Following Richens and Berry \cite{richens-berry} we impose the condition that when a geodesic, starting 
at a point in $\mathrm{int}( \mathrm{cl}(K^{\ast }) \setminus \{ O \})$, meets $\partial K^{\ast }$ it undergoes a reflection in the edge of $K^{\ast }$ that it meets. Such geodesics never meet a vertex of $\mathrm{cl}(K^{\ast })$. 
Thus this type of geodesic becomes a billiard motion in $\mathrm{cl}(K^{\ast })\setminus \{ O \}$, which is defined for all time. Billiard motions in polygons have been extensively studied. For a nice overview see Berger 
\cite[chpt. XI ]{berger} and references therein. An argument shows that $\widehat{\mathcal{G}}$ 
invariant geodesics on $({\mathcal{S}}_{\mathrm{reg}}, \Gamma )$ correspond, under the map 
${\delta }_{K^{\ast }\setminus \{ O\} }$ ($\mathrm{I}\, 7$), to billiard motions on 
$(\mathrm{cl}(K^{\ast }) \setminus \{ O \}, {\gamma }_{\mathrm{cl}(K^{\ast }) \setminus \{ O \}})$. \medskip 

Repeatedly reflecting a billiard motion in an edge of $\mathrm{cl}(K^{\ast })$ and suitable edges of 
suitable $\mathcal{T}$ translations of $\mathrm{cl}(K^{\ast })$ gives a straight line motion $\lambda $ on 
$\C \setminus {\mathbb{V}}^{+}$, which is invariant under the action of $\widehat{G} \ltimes \mathcal{T}$. 
Use the union of $\lambda $ and $U\lambda $, whose intersection with $\mathrm{cl}(K^{\ast })$ is a 
segment of an extended billiard motion. The image of this extended billiard motion in the orbit space 
$(\C \setminus {\mathbb{V}}^{+})^{\sim }/\mathfrak{G} = {\widetilde{\mathcal{S}}}_{\mathrm{reg}}$ is a geodesic. 
Here we use the Riemannian metric $\widehat{\gamma }$, 
which is induced by the $\mathfrak{G}$ invariant Euclidean metric $\gamma $ on $\C \setminus 
{\mathbb{V}}^{+}$ restricted to $\mathrm{cl}(K^{\ast })\setminus \{ O \}$. Consequently, 
$({\widetilde{\mathcal{S}}}_{\mathrm{reg}}, 
\widehat{\gamma })$ is an affine analogue of the affine Riemann surface ${\mathcal{S}}_{\mathrm{reg}}$ 
thought of as the orbit space of a discrete subgroup of $\mathrm{PGl}(2, \C)$ acting on $\C $ with 
the Poincar\'{e} metric, see Weyl \cite{weyl}. %

\section{A Schwarz-Christoffel mapping}%

Consider the conformal Schwarz-Christoffel mapping 
\begin{align}
F_T: {\C}^{+} & =  \{ \xi \in \C \setrule \, \mathrm{Im}\, \xi \ge 0 \} \rightarrow T = 
T_{n_0n_1n_{\infty}} \subseteq \C: \notag \\
&\hspace{.25in} \xi \mapsto \int^{\xi }_0 \frac{\dee w}{w^{1-\frac{n_0}{n}}{(1-w)^{1-\frac{n_1}{n}}}}= z  
\label{eq-s1one}
\end{align}
of the upper half plane ${\C}^{+}$ to the rational triangle $T = T_{n_0n_1n_{\infty}}$ with interior 
angles $\frac{n_0}{n} \pi $, $\frac{n_1}{n} \pi $, and $\frac{n_{\infty}}{n} \pi $, see figure $1$.  Here $n_0+n_1+n_{\infty} =n$ and $n_i \in {\Z }_{\ge 1}$ for $i=0,1$ and 
$\infty$ with $1\le n_0 \le n_1 \le n_{\infty}$. Because $n_{\infty}$ is greater than or equal to either $n_0$ or $n_1$, it follows that $OC$ is the longest side of the triangle \linebreak  
\par \noindent \hspace{1.2in}\begin{tabular}{l}
\includegraphics[width=180pt]{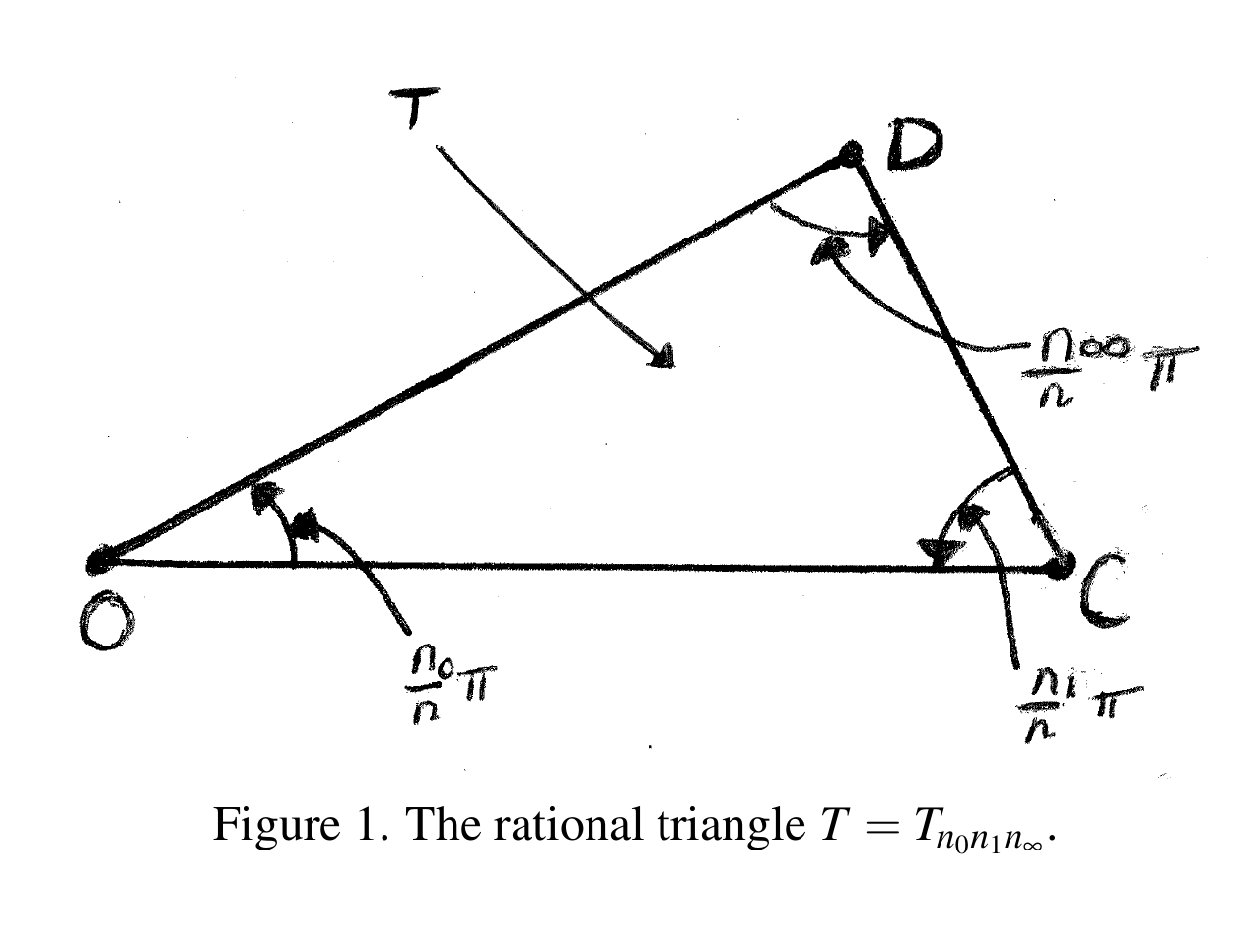}
\end{tabular} 
\par \noindent $T = \bigtriangleup OCD$. In the integrand of 
(\ref{eq-s1one}) we use the following choice of complex $n^{\mathrm{th}}$ root. Suppose that $w \in \C \setminus \{ 0,1\} $. Let $w = r_0{\mathrm{e}}^{i{\theta }_0}$ and 
$1 -w = r_1 {\mathrm{e}}^{i{\theta }_1}$ where $r_0,$ $r_1 \in {\R }_{>0}$ and 
${\theta }_0$, ${\theta }_1  \in [0 , 2\pi )$. For $w \in (0,1)$ on the real axis we have 
${\theta }_0 = {\theta }_1 =0$, $w = r_0 >0$, and $1-w =r_1 >0$. So 
$\big( w^{n-n_0}(1-w)^{n-n_1} \big)^{\raisebox{-2pt}{$\scriptstyle {1/n}$} } = 
(r^{n-n_0}_0r^{n-n_1}_1)^{1/n}$. In general for $w \in \C \setminus \{ 0,1 \}$, we have 
\begin{displaymath}
\big( w^{n-n_0}(1-w)^{n-n_1} \big)^{\raisebox{-2pt}{$\scriptstyle {1/n}$}} = 
(r^{n-n_0}_0r^{n-n_1}_1)^{1/n}{\mathrm{e}}^{i ((n-n_0){\theta }_0+(n-n_1){\theta }_1)/n}.    
\end{displaymath} 

From (\ref{eq-s1one}) we get   
\begin{displaymath}
F_{T}(0) = 0, \, \, F_T(1) = C, \, \, \mathrm{and} \, \,  F_T(\infty) = D, 
\end{displaymath} 
where $C = \int^1_0 \frac{\dee w}{w^{1-\frac{n_0}{n}}(1-w)^{1-\frac{n_1}{n}}}$ and 
$D ={\mathrm{e}}^{\frac{n_0}{n} \pi i}\big( \frac{\raisebox{4pt}{${\scriptscriptstyle \sin \frac{n_1}{n}\pi }$}}{\rule{0pt}{7pt}\raisebox{0pt}{$\scriptscriptstyle \sin \frac{n_{\infty}}{n} \pi $}} \big) C$.  Consequently, the bijective holomorphic mapping $F_T$ sends $\mathrm{int}( {\C }^{+}\setminus \{ 0, 1 \} )$, the interior of the upper half plane 
less $0$ and $1$, onto 
$\mathrm{int}\, T$, the interior of the rational triangle $T = T_{n_0n_1n_{\infty}}$, 
and sends the boundary of ${\C }^{+} \setminus \{ 0, 1 \}$ to the edges of 
$\partial T$ less their end points $O$, $C$ and $D$, see figure 1. Thus the image 
of ${\C }^{+} \setminus \{ 0, 1\} $ under $F_T$ is $\mathrm{cl}(T) \setminus 
\{ O, C, D \} $. Here $\mathrm{cl}(T)$ is the closure of $T$ in $\C $. \medskip 

Because $F_T|_{[0,1]}$ is real valued, we may use the Schwarz reflection principle to extend $F_T$ to the holomorphic diffeomorphism
\begin{align}
& \hspace{-5pt} F_Q : \C \setminus \{ 0, 1 \} \rightarrow Q = T \cup \overline{T} \subseteq \C  :  
\xi \mapsto z = \left\{ \begin{array}{rl} 
F_T(\xi ), & \mbox{if $\xi \in {\C}^{+} \setminus \{ 0,1 \} $} \\
\rule{0pt}{14pt} \overline{F_T(\overline{\xi })}, & \mbox{if $\xi  \in \overline{{\C}^{+} \setminus \{ 0,1 \}}$.} \end{array} \right. 
\label{eq-s1two} 
\end{align} 
Here $Q= Q_{n_0n_1n_{\infty}}$ is a quadrilateral with internal angles 
$2\pi \frac{n_0}{n}$, $\pi \frac{n_{\infty}}{n}$, 
$2\pi \frac{n_1}{n}$, and $\pi \frac{n_{\infty}}{n}$ and vertices at $O$, $D$, $C$, and $\overline{D}$, 
see figure 2. The conformal mapping $F_Q$ \linebreak  
\par\noindent \hspace{1in}\begin{tabular}{l}
\includegraphics[width=180pt]{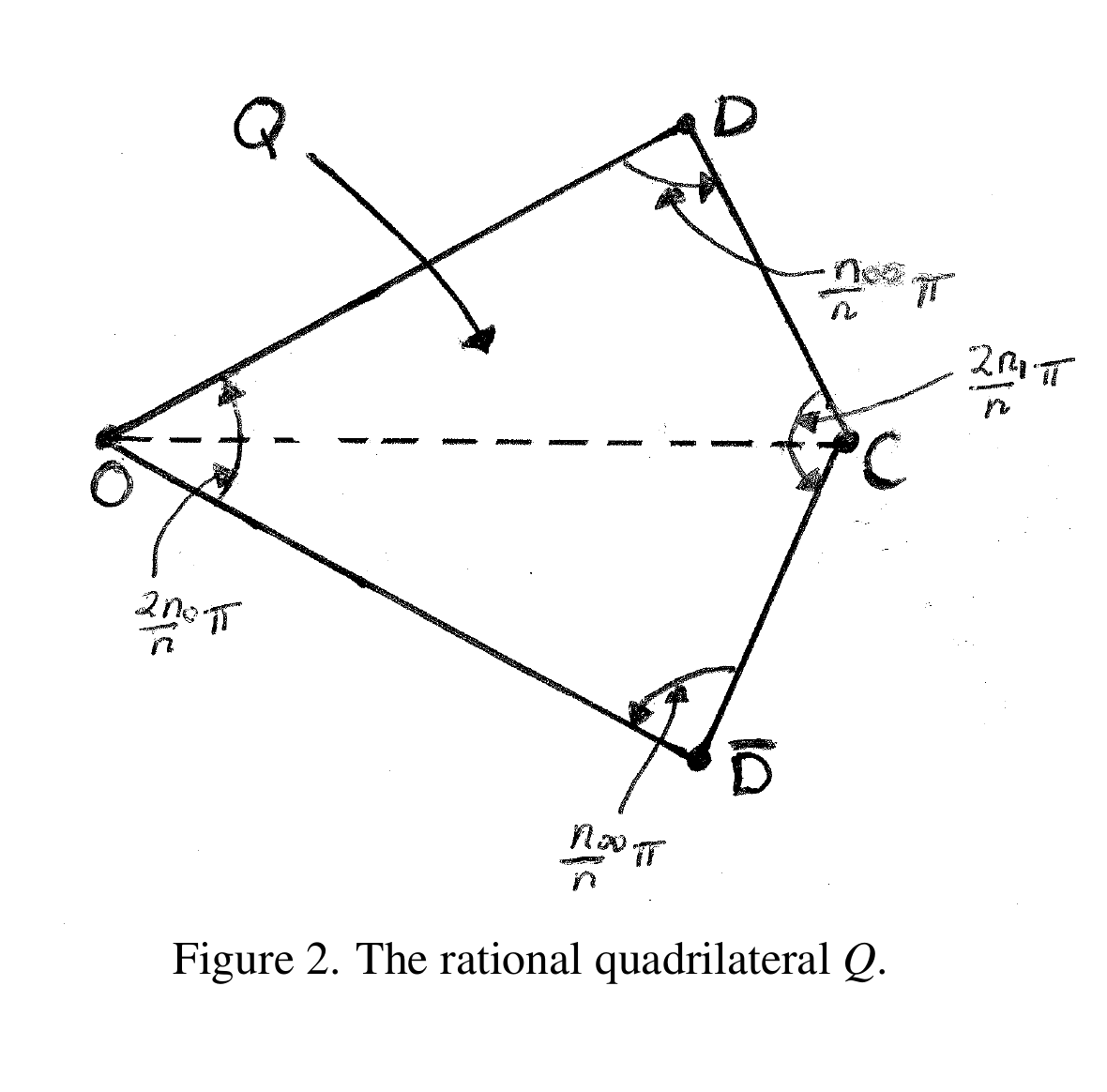}
\end{tabular}
\par \noindent sends $\C \setminus \{ 0,1 \}$ onto $\mathrm{cl}(Q) \setminus \{ O, D,C, \overline{D} \} $. 

\section{The geometry of an affine Riemann surface}

Let $\xi $ and $\eta $ be coordinate functions on ${\C }^2$. Consider the affine Riemann surface 
$\mathcal{S} = {\mathcal{S}}_{n_0,n_1,n_{\infty}}$ in ${\C }^2$, associated to the holomorphic mapping $F_Q$, defined by 
\begin{equation}
g(\xi ,\eta ) = {\eta }^n - {\xi }^{n-n_0}(1-\xi )^{n-n_1} =0 , 
\label{eq-s2one}
\end{equation}
see \cite{aurell-itzykson}. 

We determine the singular points of $\mathcal{S}$ by solving  
\begin{align}
0 & = \dee g (\xi , \eta ) \notag \\
& = -(n-n_0) {\xi }^{n-n_0-1}(1-\xi )^{n-n_1-1}(1-\ttfrac{2n-n_0-n_1}{n-n_0} \xi) \dee \xi 
+ n{\eta }^{n-1} \dee \eta 
\label{eq-s2onestar} 
\end{align}
For $(\xi ,\eta ) \in \mathcal{S}$, we have $\dee g (\xi ,\eta ) =0$  if and only if $(\xi ,\eta ) = (0,0)$ or $(1,0)$. Thus the set 
${\mathcal{S}}_{\mathrm{sing}}$ of singular points of $\mathcal{S}$ is  
$\{ (0,0), (1,0) \}$. So the affine Riemann surface ${\mathcal{S}}_{\mathrm{reg}} = 
\mathcal{S} \setminus {\mathcal{S}}_{\mathrm{sing}}$ is 
a complex submanifold of ${\C}^2$. Actually, ${\mathcal{S}}_{\mathrm{reg}} \subseteq 
{\C }^2 \setminus \{ \eta =0 \} $, for if $(\xi ,\eta ) \in \mathcal{S}$ and 
$\eta =0$, then either $\xi =0$ or $\xi =1$.  \medskip

\noindent \textbf{Lemma 2.1} Topologically ${\mathcal{S}}_{\mathrm{reg}}$ is a compact Riemann 
surface $\overline{\mathcal{S}} \subseteq {\CP}^2$ of genus $2g = n+2-(d_0+d_1+d_{\infty} )$ 
less three points: $[0:0:1]$, $[1:0:1]$, and $[0:1:0]$. Here $d_i= \gcd (n_i, n)$ for $i=0,1, \infty$,\medskip 

\noindent \textbf{Proof.} Consider the (projective) Riemann surface $\overline{\mathcal{S}} 
\subseteq {\CP }^2$ specified by the condition $[\xi : \eta : \zeta ] \in \overline{\mathcal{S}}$ if and only if \begin{equation}
G(\xi ,\eta , \zeta ) = {\zeta }^{n-n_0-n_1}{\eta }^n - {\xi }^{n-n_0}(\zeta -\xi )^{n -n_1} =0.
\label{eq-s2four}
\end{equation}
Thinking of $G$ as a polynomial in $\eta $ with coefficients which are polynomials in 
$\xi $ and $\zeta $, we may view $\overline{\mathcal{S}}$ as the branched covering 
\begin{equation}
\overline{\pi } : \overline{\mathcal{S}} \subseteq {\CP}^2 \rightarrow \CP : 
[\xi : \eta : \zeta ] \mapsto [\xi : \zeta ].
\label{eq-s2five}
\end{equation}
When $\zeta =1$ we get the affine branched covering 
\begin{equation}
\pi = \overline{\pi}|\mathcal{S}: \mathcal{S} = \overline{\mathcal{S}} \cap \{ \zeta = 1 \} 
\subseteq {\C }^2 \rightarrow \C = \CP \cap \{ \zeta =1 \} : (\xi ,\eta ) \mapsto \xi.
\label{eq-s2six}
\end{equation}
From (\ref{eq-s2one}) it follows that $\eta = {\omega }_k({\xi }^{n-n_0}(1-\xi )^{n-n_1})^{1/n}$, where 
${\omega }_k$ for $k=0,1, \ldots , n-1$ is an $n^{\mathrm{th}}$ root of unity with 
and $(\, \, \, \, )^{1/n}$ is the complex $n^{\mathrm{th}}$ root used in the definition of the mapping $F_T$ (\ref{eq-s1one}). Thus the branched covering mapping $\overline{\pi}$ (\ref{eq-s2five}) has $n$ ``sheets'' except at its branch points. Since
\begin{subequations}
\begin{align}
\eta & = {\xi }^{1-\frac{n_0}{n}}(1-\xi )^{1- \frac{n_1}{n}} = 
{\xi }^{1-\frac{n_0}{n}}\big( 1 - (1-\ttfrac{n_1}{n})\xi + \cdots \big) 
\label{eq-s2sevena} 
\end{align}
and
\begin{align}
\eta & =  ( 1 - \xi )^{1-\frac{n_1}{n}}\big( 1 - (1-\xi ) \big)^{1-\frac{n_0}{n}} \notag \\
& = (1-\xi )^{1- \frac{n_1}{n}}\big( 1 - (1-\ttfrac{n_0}{n})(1-\xi ) + \cdots \big) , 
\label{eq-s2sevenb} 
\end{align}
it follows that $\xi =0$ and $\xi =1$ are branch points of the mapping $\overline{\pi }$ of 
degree $\frac{n}{d_0}$ and $\frac{n}{d_1}$, since $d_j = \mathrm{gcd}(n, n_j) = 
\mathrm{gcd}(n-n_j,n_j)$ for $j =0,1$, see McKean and Moll\cite[p.39]{mckean-moll}. Because
\begin{align}
\eta & =  \big( \ttfrac{1}{\xi } \big)^{-(1-\frac{n_0}{n})}
\big(1- \frac{1}{\frac{1}{\xi }} \big)^{1-\frac{n_1}{n}} =
(-1)^{1-\frac{n_1}{n}} {\xi }^{2-\frac{n_0+n_1}{n}} (1-\ttfrac{1}{\xi })^{1-\frac{n_1}{n}} 
\notag \\
& = (-1)^{1-\frac{n_1}{n}} {\xi }^{1+\frac{n_{\infty}}{n}}\big( 1-(1-\ttfrac{n_1}{n}) 
\ttfrac{1}{\xi } + \cdots \big) , 
\label{eq-s2sevenc} 
\end{align}
\end{subequations}
$\infty $ is a branch point of the mapping $\overline{\pi }$ of degree $\frac{n}{d_{\infty}}$, 
where $d_{\infty} = \mathrm{gcd}(n,n_{\infty})$. Hence the ramification 
index of $0$, $1$, $\infty$ is $d_0(\frac{n}{d_0}-1) = n - d_0$, $n - d_1$, and 
$n-d_{\infty}$, respectively. Thus the map $\overline{\pi}$ has 
$d_0$ fewer sheets at $0$, $d_1$ fewer at $1$, and $d_{\infty}$ fewer at $\infty$ than an $n$-fold covering of $\CP$. Thus the total ramification index $r$ of the mapping $\overline{\pi }$ is $r = (n-d_0)+(n-d_1)+(n-d_{\infty})$. By the Riemann-Hurwitz formula, the genus $g$ of $\overline{\mathcal{S}}$ is $r = 2n + 2g - 2$. In other words, 
\begin{equation}
2g = n+2 -(d_0+d_1+d_{\infty}). 
\label{eq-s2eightstar}
\end{equation}
Consequently, the affine Riemann surface $\mathcal{S}$ is the compact connected surface 
$\overline{\mathcal{S}}$ less the point at $\infty $, namely, $\mathcal{S} = 
\overline{\mathcal{S}}\setminus \{ [0:1:0] \}$. So 
${\mathcal{S}}_{\mathrm{reg}}$ is the compact connected surface $\overline{\mathcal{S}}$ 
less three points: $[0:0:1]$, $[1:0:1]$, and $[0:1:0]$. \hfill $\square $ \medskip

\noindent \textbf{Examples of $\overline{\mathcal{S}} = {\overline{\mathcal{S}}}_{n_0,n_1,n_{\infty}}$}
\par \noindent \qquad 1. \parbox[t]{4in}{$n_0 =1$, $n_1 =1$, $n_{\infty}= 1$; $n =3$. So $d_0=d_1=d_{\infty} =1$. Hence $2g =5 -3 =2$. So $g =1$.}
\smallspace
\par \noindent \qquad 2. \parbox[t]{4in}{$n_0 =1$, $n_0 =1$, $n_{\infty}= 4$; $n =6$. So $d_0=1$, $d_1 =1$, $d_{\infty} =2$. Hence $2g =8 - 4 =4$. So $g =2$.}
\smallspace
\par \noindent \qquad 3. \parbox[t]{4in}{$n_0 =1$, $n_1 =2$, $n_{\infty}= 3$; $n =6$. So $d_0=1$, $d_1=2$, $d_{\infty} =3$. Hence $2g =8 - 6 =2$. So $g =1$.}
\smallspace
\par \noindent \qquad 4. \parbox[t]{4in}{$n_0 =2$, $n_1 =2$, $n_{\infty}= 3$; $n =7$. So $d_0=d_1=d_{\infty} =1$. Hence $2g =9 - 3 =6$. So $g =3$.\hfill $\square $} \medskip 

Below is a table listing all the partitions $\{ n_1, n_0, n_{\infty} \} $ of $n$, which give a low genus Riemann surface $\overline{\mathcal{S}}= 
{\overline{\mathcal{S}}}_{n_0,n_1,n_{\infty}}$ \bigskip

\hspace{.75in}\begin{tabular}{ccll}
\multicolumn{1}{c}{$g$} & \multicolumn{1}{c}{\hspace{-8pt}$n_0,n_1,n_{\infty};n$} &
\multicolumn{1}{l}{$g$} &\multicolumn{1}{l}{\hspace{9pt}$n_0,n_1,n_{\infty};n$} \\ \cline{1-4} 
$1$ & $\hspace{12pt} 1,\hspace{5pt} 1,\hspace{5pt} 1; \hspace{10pt}  3$ \hspace{.2in} \vline & $3$ & $\hspace{12pt} 2,\hspace{5pt}2,\hspace{5pt}3;\hspace{8pt} 7 $ \\
$1$ & $\hspace{12pt} 1,\hspace{5pt}1,\hspace{5pt}2; \hspace{10pt}4$ \hspace{.2in} \vline &
$3$ & $\hspace{12pt}1,\hspace{5pt}3,\hspace{5pt}3; \hspace{8pt} 7$ \\
$1$ & $\hspace{12pt}1,\hspace{5pt}2,\hspace{5pt}3;\hspace{10pt} 6$ \hspace{.2in} \vline &
$3$ & $\hspace{12pt}1,\hspace{5pt}1,\hspace{5pt}5;\hspace{8pt} 7$ \\
$2$ & $\hspace{12pt}1,\hspace{5pt}2,\hspace{5pt}2;\hspace{10pt} 5$ \hspace{.2in} \vline &
$3$ & $\hspace{12pt}2,\hspace{5pt}3,\hspace{5pt}3;\hspace{8pt} 8$ \\
$2$ & $\hspace{12pt}1,\hspace{5pt}1,\hspace{5pt}3;\hspace{10pt} 5$ \hspace{.2in} \vline &
$3$ & $\hspace{12pt}1,\hspace{5pt}2,\hspace{5pt}5;\hspace{8pt} 8$ \\
$2$ & $\hspace{12pt}1,\hspace{5pt}1,\hspace{5pt}4;\hspace{10pt} 6$ \hspace{.2in} \vline & 
$3$ & $\hspace{12pt}1,\hspace{5pt}1,\hspace{5pt}6;\hspace{8pt} 8$ \\
$2$ & $\hspace{12pt}1,\hspace{5pt}3,\hspace{5pt}4;\hspace{10pt} 8$ \hspace{.2in} \vline & 
$3$ & $\hspace{12pt}2,\hspace{5pt}3,\hspace{5pt}4;\hspace{8pt} 9$ \\
$2$ & $\hspace{12pt}2,\hspace{5pt}3,\hspace{5pt}5;\hspace{8pt}10$ \hspace{.16in} \vline & 
$3$ & $\hspace{12pt}1,\hspace{5pt}3,\hspace{5pt}5;\hspace{8pt} 9$ \\
$2$ & $\hspace{12pt}1,\hspace{5pt}4,\hspace{5pt}5;\hspace{8pt}10$ \hspace{.16in} \vline & 
$3$ & $\hspace{12pt}1,\hspace{5pt}2,\hspace{5pt}6;\hspace{8pt} 9$ \\
& \hspace{1.15in} \vline & 
$3$ & $\hspace{12pt}3,\hspace{5pt}4,\hspace{5pt}5;\hspace{8pt} 12$ \\
& \hspace{1.15in} \vline &                          
$3$ & $\hspace{12pt}1,\hspace{5pt}5,\hspace{5pt}6;\hspace{8pt} 12 $ \\
& \hspace{1.15in} \vline &                          
$3$ & $\hspace{12pt}1,\hspace{5pt}3,\hspace{5pt}8;\hspace{8pt} 12$ \\
& \hspace{1.15in} \vline &                           
$3$ & $\hspace{12pt}2,\hspace{5pt}5,\hspace{5pt}7;\hspace{8pt} 14$ \\
& \hspace{1.15in} \vline &                            
$3$ & $\hspace{12pt}1,\hspace{5pt}6,\hspace{5pt}7;\hspace{8pt} 14$ 
\end{tabular}
\vspace{.1in}
\par\noindent \hspace{.65in}\parbox[t]{3.5in}{Table 1.Genus $g$ of 
$\overline{\mathcal{S}} = {\overline{\mathcal{S}}}_{n_0,n_1,n_{\infty}}$. This table is based on the table in 
Aurell and Itzykson \cite[p.193]{aurell-itzykson}.} \medskip

\noindent \textbf{Corollary 2.1a} If $n$ is an odd prime number and $\{ n_1, n_0, n_{\infty} \} $ is a partition of $n$ into three parts, then the genus of $\overline{\mathcal{S}}$ is $\onehalf (n-1)$. \medskip 

\noindent \textbf{Proof.} Because $n$ is prime, we get $d_0=d_1=d_{\infty} =1$. Using 
(\ref{eq-s2eightstar}) we obtain $g = \onehalf (n-1)$. \hfill $\square $ \medskip 

\par\noindent \textbf{Corollary 2.1b} The singular points of the Riemann surface $\overline{\mathcal{S}}$ are 
$[0:0:1]$, $[1:0:1]$, and if $n_{\infty} >1$ then also $[0:1:0]$. \medskip

\noindent \textbf{Proof.} A point $[\xi :\eta :\zeta ] \in 
{\overline{\mathcal{S}}}_{\mathrm{sing}}$ if and only if $[\xi :\eta : \zeta ] \in \overline{\mathcal{S}}$, that is, 
\begin{subequations}
\begin{equation}
0 = G(\xi ,\eta ,\zeta ) = {\zeta }^{n-(n_0+n_1)}{\eta }^n - {\xi }^{n-n_0}(\zeta -\xi )^{n-n_1}
\label{eq-s2eighta}
\end{equation} 
and  
\begin{align}
(0,0,0) & = DG (\xi ,\eta , \zeta  ) \notag \\
& = \big( -{\xi }^{n-n_0 -1}(\zeta - \xi )^{n-n_1 -1}
\big( (n-n_0)(\zeta - \xi ) - (n-n_1)\xi \big) ,   \notag \\
&\hspace{.5in} n{\eta }^{n-1}{\zeta }^{n-(n_0+n_1)}, (n-(n_0+n_1)) {\eta }^n{\zeta }^{n-n_0-n_1-1}   
\notag \\
& \hspace{.75in}  -(n-n_1){\xi }^{n-n_0}(\zeta - \xi )^{n-n_1-1} \big) 
\label{eq-s2eightb} 
\end{align}
\end{subequations}
We need only check the points $[0:0:1]$, $[1:0:1]$ and $[0:1:0]$. Since the first two points are singular points of 
$\mathcal{S} = \overline{\mathcal{S}}\setminus \{ [0:1:0] \}$, they 
are singular points of $\overline{\mathcal{S}}$. Thus we need to see if $[0:1:0]$ is a singular point of $\overline{\mathcal{S}}$. Substituting $(0,1,0)$ into the right hand side of (\ref{eq-s2eightb}) we get {\tiny $\left\{ \! \! \! \! \begin{array}{l} (0,0,1),\mbox{\, if $n_{\infty} = n-(n_0+n_1) = 1$} \\ (0,0,0),\mbox{\, if $n_{\infty} >1$.} \end{array} \right.$} Thus $[0:1:0]$ is a singular point of $\overline{\mathcal{S}}$ only if $n_{\infty} >1$. \hfill $\square $ \medskip 

\noindent \textbf{Lemma 2.2} The mapping  
\begin{equation}
\widehat{\pi } = \pi |{\mathcal{S}}_{\mathrm{reg}}: 
{\mathcal{S}}_{\mathrm{reg}} \subseteq {\C}^2 \rightarrow \C \setminus \{ 0,1 \} : 
(\xi ,\eta ) \mapsto \xi 
\label{eq-s2seven}
\end{equation}
is a surjective holomorphic local diffeomorphism. \medskip 

\noindent \textbf{Proof.} Let $(\xi ,\eta ) \in {\mathcal{S}}_{\mathrm{reg}}$ and let 
\begin{equation}
X(\xi ,\eta ) = \eta \frac{\partial }{\partial \xi } +
\ttfrac{n-n_0}{n} \frac{{\xi }^{n-n_0-1}(1-\xi)^{n-n_1-1}(1-\ttfrac{2n-n_0-n_1}{n-n_0} \xi )}{{\eta }^{n-2} } \frac{\partial }{\partial \eta }.
\label{eq-s2eight}
\end{equation}
By hypothesis $(\xi , \eta ) \in {\mathcal{S}}_{\mathrm{reg}}$ implies that $\eta \ne 0$. 
The vector $X(\xi ,\eta )$ is defined and is nonzero. From $(X \lefthook \dee g) (\xi ,\eta ) =0$ and $T_{(\xi ,\eta )}{\mathcal{S}}_{\mathrm{reg}} = \ker \dee g (\xi ,\eta )$, it follows that $X(\xi ,\eta )  \in T_{(\xi ,\eta )}{\mathcal{S}}_{\mathrm{reg}}$. Using the definition of $X(\xi ,\eta )$ (\ref{eq-s2eight}) and the definition of the mapping $\pi $ (\ref{eq-s2six}), we see that the tangent of the mapping 
$\widehat{\pi }$ (\ref{eq-s2seven}) at $(\xi ,\eta ) \in {\mathcal{S}}_{\mathrm{reg}}$ is given by 
\begin{equation}
T_{(\xi ,\eta )}\widehat{\pi }:T_{(\xi ,\eta )}{\mathcal{S}}_{\mathrm{reg}} \rightarrow T_{\xi }(\C \setminus \{ 0, 1 \} ) = \C:  X(\xi ,\eta ) \mapsto \eta \frac{\partial }{\partial \xi }. 
\label{eq-s2eightdot}
\end{equation}
Since $X(\xi ,\eta )$ and $\eta \frac{\partial }{\partial \xi }$ 
are nonzero vectors, they form a complex basis for $T_{(\xi ,\eta )}{\mathcal{S}}_{\mathrm{reg}}$ and $T_{\xi }(\C \setminus \{ 0,1\})$, respectively. Thus the complex linear mapping 
$T_{(\xi ,\eta )}\widehat{\pi }$ is an isomorphism. Hence 
$\widehat{\pi }$ is a local holomorphic diffeomorphism. \hfill $\square $ \medskip 

\noindent \textbf{Corollary 2.2a} $\widehat{\pi }$ (\ref{eq-s2seven}) is a surjective 
holomorphic $n$ to $1$ covering map. \medskip 

\noindent \textbf{Proof.} We only need to show that $\widehat{\pi }$ is a covering map. First 
we note that every fiber of $\widehat{\pi }$ is a finite set with $n$ elements, 
since for each fixed $\xi \in \C \setminus \{ 0,1 \}$ we have 
${\widehat{\pi }}^{-1}(\xi ) = \{ (\xi ,\eta ) \in {\mathcal{S}}_{\mathrm{reg}} 
\setrule \eta = {\omega }_k  ({\xi }^{n-n_0}(1-\xi )^{n-n_1})^{1/n} \} $. Here ${\omega }_k$ for 
$k=0,1, \ldots , n-1$, is an $n^{\mathrm{th}}$ root of $1$ and $(\, \, \, \, )^{1/n}$ is the complex $n^{\mathrm{th}}$ root used in the definition of the Schwarz-Christoffel map $F_Q$ (\ref{eq-s1two}). Hence the map $\widehat{\pi }$ is a proper surjective holomorphic submersion, because each fiber is compact. Thus the mapping $\widehat{\pi }$ is a presentation of a locally trivial fiber bundle with fiber consisting of $n$ distinct points. In other words, the map 
$\widehat{\pi }$ is a $n$ to $1$ covering mapping. \hfill $\square $ \medskip 

Consider the group $\widehat{\mathcal{G}}$ of linear transformations of ${\C }^2$ generated by  
\begin{displaymath}
\mathcal{R}: {\C }^2 \rightarrow {\C }^2: (\xi ,\eta ) \mapsto (\xi , {\mathrm{e}}^{2\pi i /n} \eta ).
\end{displaymath}
Clearly ${\mathcal{R}}^n = {\mathrm{id}}_{{\C}^2} = e$, the identity element of 
$\widehat{\mathcal{G}}$ and $\widehat{\mathcal{G}} = 
\{ e, \mathcal{R}, \ldots , {\mathcal{R}}^{n-1} \} $.  For each $(\xi , \eta ) \in \mathcal{S}$ we have 
\begin{align}
({\mathrm{e}}^{2\pi i /n}\eta )^n - {\xi }^{n-n_0}(1-\xi )^{n-n_1} & = 
{\eta }^n - {\xi }^{n-n_0}(1- \xi )^{n-n_1}  = 0. \notag
\end{align}
So $\mathcal{R}(\xi ,\eta ) \in \mathcal{S}$. Thus we have an action of $\widehat{\mathcal{G}}$ on the affine Riemann surface $\mathcal{S}$ given by 
\begin{equation}
\Phi : \widehat{\mathcal{G}} \times \mathcal{S} \rightarrow \mathcal{S}:\big( g, (\xi , \eta ) \big) \mapsto 
g( \xi ,\eta ) . 
\label{eq-s2ten}
\end{equation}
Since $\widehat{\mathcal{G}}$ is finite, and hence is compact, the action $\Phi $ is proper. 
For every $g \in \widehat{\mathcal{G}}$ we have ${\Phi }_g(0,0) = (0,0)$ and ${\Phi }_g(1,0) = (1,0)$. So ${\Phi }_g$ maps ${\mathcal{S}}_{\mathrm{reg}}$ into itself. At 
$(\xi ,\eta ) \in {\mathcal{S}}_{\mathrm{reg}}$ the isotropy group 
${\widehat{\mathcal{G}}}_{(\xi ,\eta )}$ is $\{ e \}$, that is, the $\widehat{\mathcal{G}}$-action $\Phi $ on ${\mathcal{S}}_{\mathrm{reg}}$ is free. Thus the orbit space 
${\mathcal{S}}_{\mathrm{reg}}/\widehat{\mathcal{G}}$ is a complex manifold. \medskip 

\noindent \textbf{Corollary 2.2b} The holomorphic $\widehat{\mathcal{G}}$-principal bundle 
\begin{displaymath}
\rho : {\mathcal{S}}_{\mathrm{reg}} \subseteq {\C }^2 \rightarrow {\mathcal{S}}_{\mathrm{reg}}/\widehat{\mathcal{G}} \subseteq {\C }^2:   (\xi ,\eta ) \mapsto 
[(\xi ,\eta )] .
\end{displaymath}
Here $[(\xi ,\eta )]$ is the $\widehat{\mathcal{G}}$-orbit 
$\{ {\Phi }_g(\xi ,\eta ) \in {\mathcal{S}}_{\mathrm{reg}} 
\setrule g \in \widehat{\mathcal{G}} \} $ of $(\xi ,\eta )$ in 
${\mathcal{S}}_{\mathrm{reg}}$. The 
bundle presented by the mapping $\rho $ is isomorphic to the bundle presented by the mapping $\widehat{\pi }$ (\ref{eq-s2seven}). \medskip 

\noindent \textbf{Proof.} We use invariant theory to determine the orbit space 
$\mathcal{S}/\widehat{\mathcal{G}}$. The algebra of polynomials on ${\C }^2$, which are invariant under the $\widehat{\mathcal{G}}$-action $\Phi $, is generated by 
${\pi }_1 = \xi \, \, \, \mathrm{and} \, \, \, {\pi }_2 = {\eta }^n$. 
Since $(\xi , \eta ) \in \mathcal{S}$, these polynomials are subject to the relation 
\begin{equation}
{\pi }_2 - {\pi }^{n-n_0}_1(1-{\pi }_1)^{n-n_1} = {\eta }^n- {\xi }^{n-n_0}(1-\xi )^{n-n_1} =0. 
\label{eq-s2twelve}
\end{equation}
Equation (\ref{eq-s2twelve}) defines the orbit space $\mathcal{S}/\widehat{\mathcal{G}}$ as a complex subvariety of ${\C }^2$. This 
subvariety is homeomorphic to $\C $, because it is the graph of the function 
${\pi }_1 \mapsto {\pi }^{n-n_0}_1(1-{\pi }_1)^{n-n_1}$. Consequently, the orbit space 
${\mathcal{S}}_{\mathrm{reg}}/\widehat{\mathcal{G}}$ is holomorphically diffeomorphic to $\C \setminus \{ 0, 1 \} $. \medskip 

It remains to show that $\widehat{\mathcal{G}}$ is the group of covering transformations of the bundle presented by the mapping $\widehat{\pi }$ (\ref{eq-s2seven}). For each 
$\xi \in \C \setminus \{ 0,1 \} $ look at the fiber ${\widehat{\pi }}^{-1}(\xi )$. If 
$(\xi ,\eta ) \in {\widehat{\pi }}^{-1}(\xi )$, then ${\mathcal{R}}^{\pm 1}(\xi ,\eta ) = 
(\xi , {\mathrm{e}}^{\pm 2\pi i /n}\eta ) \in {\mathcal{S}}_{\mathrm{reg}}$, since 
$(\xi , {\mathrm{e}}^{\pm 2\pi i /n}\eta ) \ne (0,0)$ or $(1,0)$ and 
$(\xi ,{\mathrm{e}}^{\pm 2\pi i /n}\eta ) 
\in \mathcal{S}$. Thus ${\Phi }_{{\mathcal{R}}^{\pm 1}} \big({\widehat{\pi }}^{-1}(\xi ) \big) \subseteq {\widehat{\pi }}^{-1}(\xi )$. So ${\widehat{\pi }}^{-1}(\xi ) \subseteq 
{\Phi }_{\mathcal{R}}\big( {\widehat{\pi }}^{-1}(\xi )\big) \subseteq 
{\widehat{\pi }}^{-1}(\xi )$. 
Hence ${\Phi }_{\mathcal{R}}\big( {\widehat{\pi }}^{-1}(\xi )\big) =
{\widehat{\pi }}^{-1}(\xi )$. Thus ${\Phi }_{\mathcal{R}}$ is a covering transformation for the bundle presented by the mapping $\widehat{\pi }$. So $\widehat{\mathcal{G}}$ is a subgroup of the group of covering transformations. These groups are equal because  
$\widehat{\mathcal{G}}$ acts transitively on each fiber of the mapping $\widehat{\pi }$. \hfill $\square $  

\section{Another model for ${\mathcal{S}}_{\mathrm{reg}}$}

In this section we construct another model ${\widetilde{S}}_{\mathrm{reg}}$ for the smooth part 
${\mathcal{S}}_{\mathrm{reg}}$ of the affine Riemann surface $\mathcal{S}$ (\ref{eq-s2one}). \medskip 

Let $\mathcal{D} \subseteq {\mathcal{S}}_{\mathrm{reg}}$ be a fundamental domain for the $\widehat{\mathcal{G}}$ action $\Phi $ (\ref{eq-s2ten}) on 
${\mathcal{S}}_{\mathrm{reg}}$. So $(\xi ,\eta ) \in \mathcal{D}$ if and only if for 
$\xi \in \C \setminus \{ 0,1 \}$ we have 
$\eta = \big({\xi}^{n-n_0}(1-\xi)^{n-n_1}\big)^{\raisebox{-2pt}{$\scriptstyle 1/n$}}$. 
Here $(\, \, \, \, )^{1/n}$ is the $n^{\mathrm{th}}$ root used in the definition of the mapping $F_Q$ (\ref{eq-s1two}). The domain $\mathcal{D}$ is a connected subset of 
${\mathcal{S}}_{\mathrm{reg}}$ with nonempty interior. Its image under the map $\widehat{\pi }$ (\ref{eq-s2seven}) is $\C \setminus \{ 0, 1 \} $. Thus $\mathcal{D}$ is 
one ``sheet'' of the covering map $\widehat{\pi }$. So $\widehat{\pi }|_{\mathcal{D}}$ is 
one to one. 
\par\noindent \hspace{.5in}\begin{tabular}{l}
\includegraphics[width=250pt]{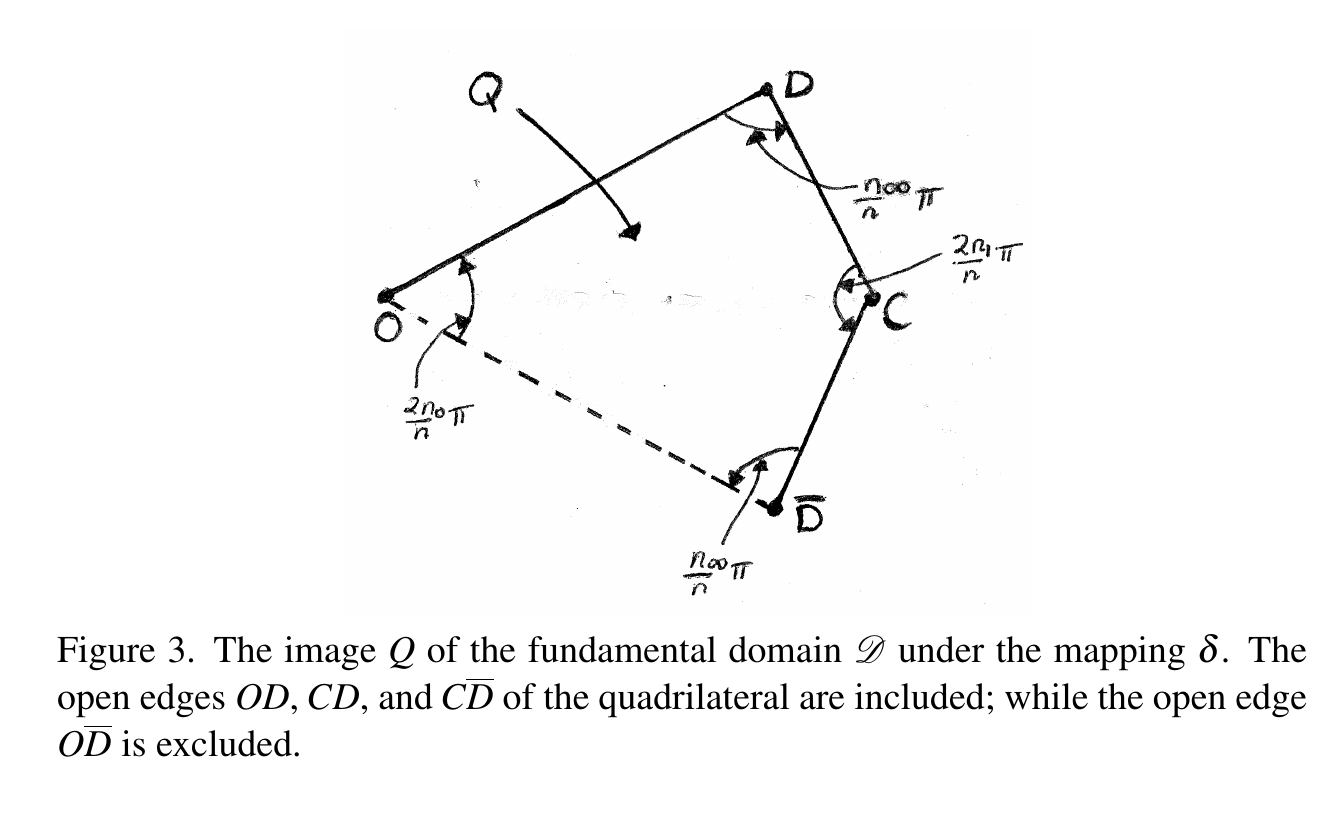}
\end{tabular}

Using the extended Schwarz-Christoffel mapping $F_Q$ (\ref{eq-s1two}), we give a more geometric description of the fundamental domain $\mathcal{D}$. Consider the mapping 
\begin{equation}
\delta : \mathcal{D} \subseteq {\mathcal{S}}_{\mathrm{reg}} \rightarrow Q \subseteq \C : 
(\xi , \eta ) \mapsto F_Q \big( \widehat{\pi }(\xi ,\eta ) \big) ,  
\label{eq-s3onestar}
\end{equation}
where the map $\widehat{\pi}$ is given by equation (\ref{eq-s2seven}). The map $\delta $ is a holomorphic diffeomorphism of $\mathrm{int}\, \mathcal{D}$ onto $\mathrm{int}\, Q$, which sends $\partial \mathcal{D}$ homeomorphically onto $\partial Q$. Look \linebreak 
\vspace{-.1in}\par\noindent \hspace{1.25in}\begin{tabular}{l}
\includegraphics[width=150pt]{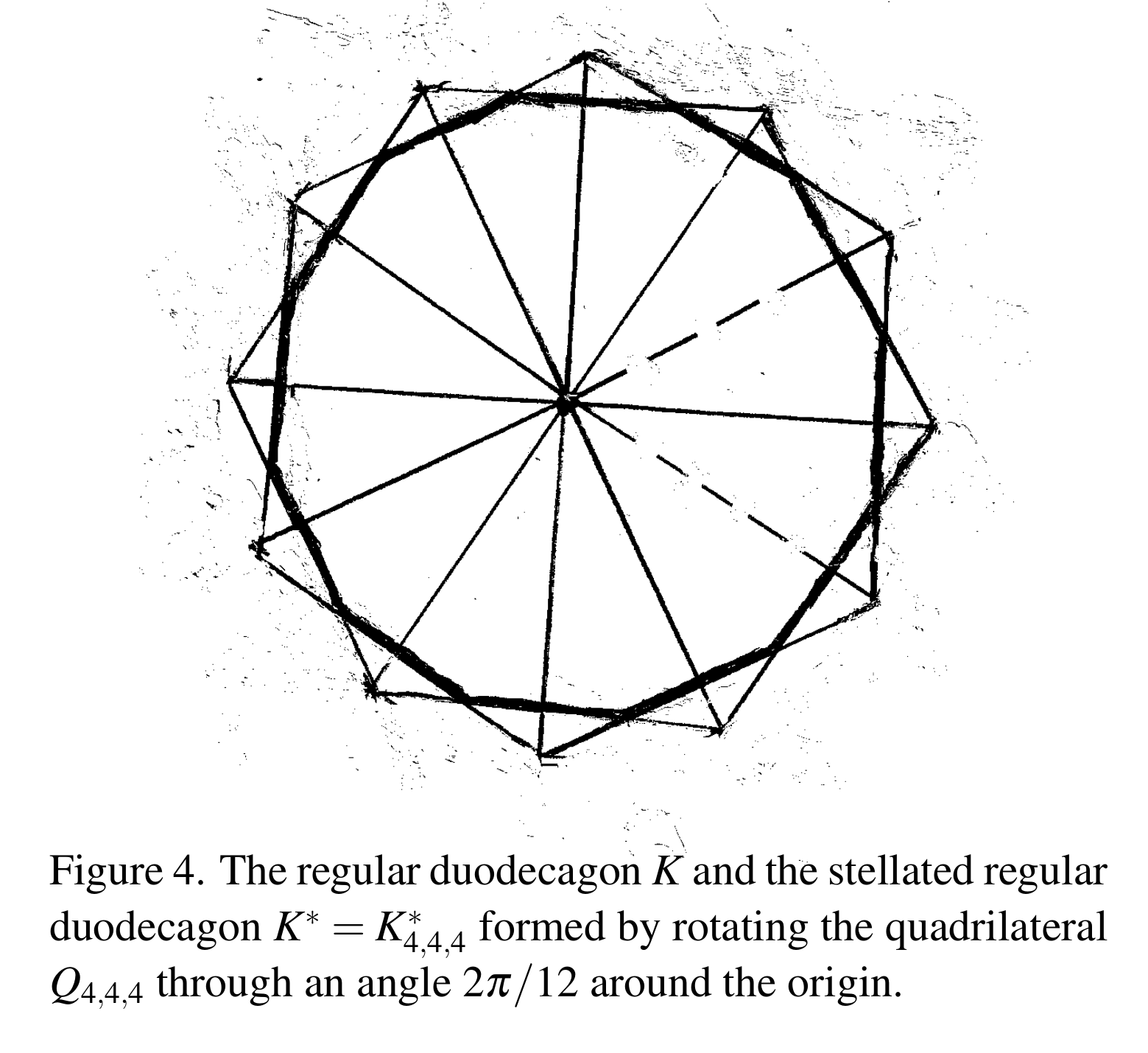}
\end{tabular}
\par \noindent at $\mathrm{cl}(Q)$, which is   a closed quadrilateral with vertices $O$, $D$, $C$, and $\overline{D}$. The set $\delta (\mathcal{D})$ 
contains the open edges $OD$, $DC$, and $C\overline{D}$ but \emph{not} the open edge $O\overline{D}$ of $\mathrm{cl}(Q)$, see figure 3 above. Let $K^{\ast } = 
K^{\ast }_{n_0,n_1,n_{\infty}} = 
{\amalg}_{0 \le j \le n-1} R^j\big( \delta (\mathcal{D}) \big) $ be 
the region in $\C $ formed by repeatedly rotating $Q = \delta (\mathcal{D})$ through an angle $2\pi /n$. Here $R$ is the rotation $\C \rightarrow \C : z\mapsto {\mathrm{e}}^{2\pi i /n} z$. We say that the quadrilateral $Q = Q_{2n_0,n_{\infty}, 2n_1,n_{\infty}}$ \emph{forms} $K^{\ast }$ less its vertices, see figure 4 above. \medskip 
 
\noindent \textbf{Claim 3.1} The connected set $K^{\ast }$ is a regular stellated $n$-gon with its $2n$ vertices omitted, which is formed from the quadrilateral $Q' = OD'C\overline{D'}$, see figure 5. \medskip

\par \noindent \textbf{Proof}. By construction the quadrilateral $Q' = OD'C\overline{D'}$
is contained in the quadrilateral $Q = ODC\overline{D}$. Note that 
$Q \subseteq \bigcup^{[\frac{n_1+1}{2}]}_{j= [-\frac{n_1+1}{2}]} R^j(Q')$. Thus 
\begin{displaymath}
K^{\ast } = \bigcup^n_{j=0} R^j(Q) \subseteq \bigcup^n_{j=0} R^j(Q')  
\subseteq \bigcup^n_{j=0} R^j(Q) = K^{\ast}. 
\end{displaymath}
So $K^{\ast } = \bigcup^n_{j=0} R^j(Q')$. Thus $K^{\ast }$ is the regular stellated $n$-gon, one of whose sides is the diagonal $D' \overline{D'}$ of $Q'$. \hfill $\square $ \medskip 

\par\noindent \hspace{.85in}\begin{tabular}{l}
\includegraphics[width=200pt]{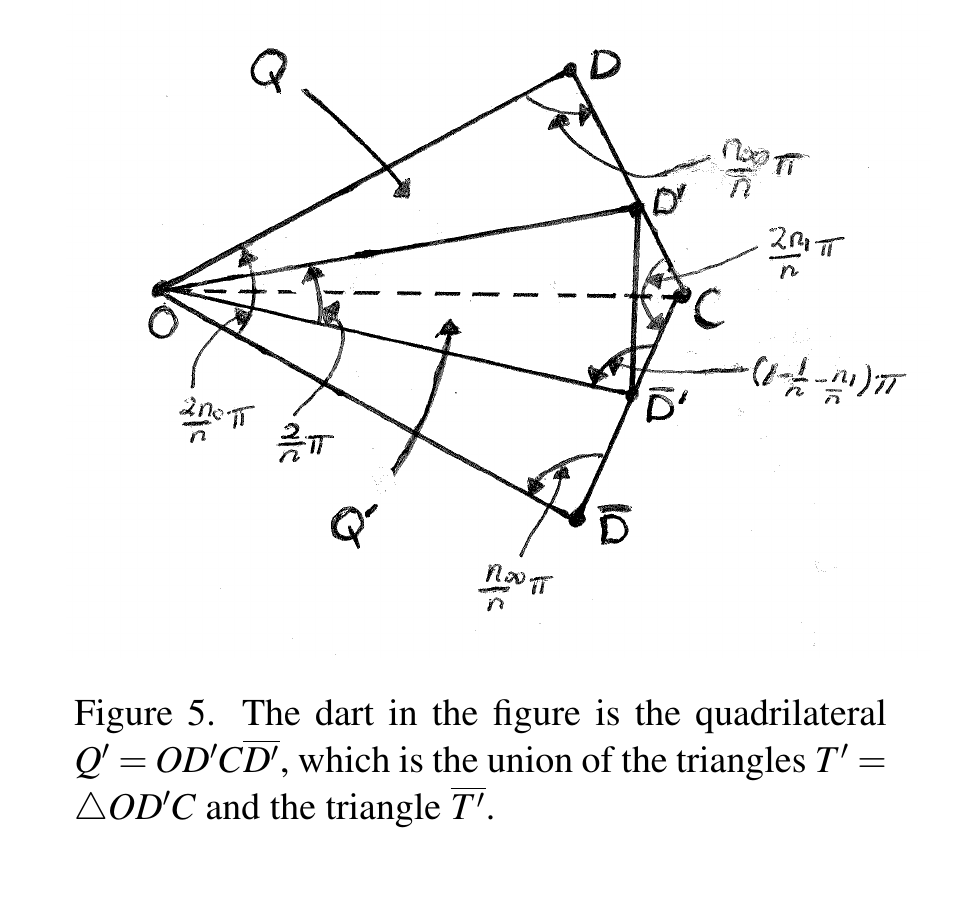} \\ 
\end{tabular} 

\vspace{-.1in}We would like to extend the  mapping $\delta $ (\ref{eq-s3onestar}) to a mapping of 
${\mathcal{S}}_{\mathrm{reg}}$ onto $K^{\ast }$. Let 
\begin{displaymath}
{\delta }_{ {\Phi }_{{\mathcal{R}}^j}(\mathcal{D}) }: {\Phi }_{ {\mathcal{R}}^j}(\mathcal{D})
\subseteq {\mathcal{S}}_{\mathrm{reg} } \rightarrow R^j \big( \delta (\mathcal{D}) \big) 
\subseteq K^{\ast }: (\xi ,\eta ) \mapsto R^j \delta 
\big( {\Phi }_{ {\mathcal{R}}^{-j} }(\xi ,\eta ) \big) ,  
\end{displaymath}  
where $\Phi $ is the $\widehat{\mathcal{G}}$ action defined in equation (\ref{eq-s2ten}). So we have a mapping 
\begin{equation}
{\delta }_{K^{\ast }}: {\mathcal{S}}_{\mathrm{reg}} \subseteq {\C }^2 \rightarrow K^{\ast } \subseteq \C 
\label{eq-s3three}
\end{equation}
defined by $({\delta }_{K^{\ast }})|_{ {\Phi }_{ {\mathcal{R}}^j }(\mathcal{D}) } = 
{\delta }|_{ {\Phi }_{ {\mathcal{R}}^j}(\mathcal{D}) }$. 
The mapping ${\delta }_{K^{\ast }}$ is defined on ${\mathcal{S}}_{\mathrm{reg}}$, 
because ${\mathcal{S}}_{\mathrm{reg}} = 
{\amalg}_{0 \le j \le n-1} {\Phi }_{{\mathcal{R}}^j}(\mathcal{D})$, since $\mathcal{D}$ is a fundamental domain for the $\widehat{\mathcal{G}}$-action 
$\Phi $ (\ref{eq-s2ten}) on ${\mathcal{S}}_{\mathrm{reg}}$. Because 
$K^{\ast } = {\amalg}_{0\le j \le n-1 } R^j \big( \delta (\mathcal{D}) \big) $, 
the mapping ${\delta }_{K^{\ast }}$ is surjective. Hence ${\delta }_{K^{\ast }}$ is holomorphic, since it is continuous on ${\mathcal{S}}_{\mathrm{reg}}$ and is holomorphic on the dense open subset ${\amalg}_{0\le j \le n-1 } {\mathcal{R}}^j(\mathrm{int}\, \mathcal{D})$ of 
${\mathcal{S}}_{\mathrm{reg}}$. \medskip %

Let $U: \C \rightarrow \C : z \mapsto \overline{z}$ and let $G$ be the group generated 
by the rotation $R$ and the reflection $U$ subject to the relations 
$R^n = U^2 =e$ and $RU = UR^{-1}$. Shorthand $G = \langle U, R \setrule 
\, U^2 =e = R^n \, \, \& \, \, RU = UR^{-1} \rangle $. Then $G = \{ e; R^pU^{\ell}, \, \ell = 0,1 \, \, \& \, \, p=0,1,\ldots , n-1 \}$. The group $G$ is the dihedral group $D_{2n}$. \medskip

The closure $\mathrm{cl}(K^{\ast})$ of $K^{\ast } = {\amalg}_{0\le j \le n-1}R^j(Q) $ is invariant under $\widehat{G}$, 
the subgroup of $G$ generated by the rotation $R$. Because the quadrilateral $Q$ is invariant under the reflection $U :z \mapsto \overline{z}$, and $UR^j = R^{-j}U$, it follows that $\mathrm{cl}(K^{\ast })$ is invariant under the reflection $U$. So $\mathrm{cl}(K^{\ast })$ is invariant under the group $G$.  \medskip 

We now look at some group theoretic properties of $K^{\ast }$. \medskip 

\noindent \textbf{Lemma 3.2} If $F$ is a closed edge of the polygon 
$\mathrm{cl}(K^{\ast })$ and $g|_F = {\mathrm{id}}|_F$ for some $g \in G$, then $g = e$. \medskip 

\noindent \textbf{Proof.} Suppose that $g \ne e$. Then $g = R^pU^{\ell }$ for some 
$\ell \in \{ 0,1 \} $ and some $p \in \{ 0,1, \ldots , n-1 \} $. Let $g = R^pU$ and suppose that 
$F$ is an edge of $\mathrm{cl}(K^{\ast}) $ such 
that $\mathrm{int}(F) \cap \R \ne \varnothing$, where $\R = \{ \mathrm{Re}\, z \setrule \, 
z \in \C \}$. Then $U(F) = F$, but $U|_F \ne 
{\mathrm{id}}_F$. So $g|_F = R^p U|_F \ne {\mathrm{id}}_F$. Now 
suppose that $\mathrm{int}(F) \cap \R = \varnothing$. Then $U(F) \ne F$. So 
$U|_F \ne {\mathrm{id}}_F$. Hence $g|_F \ne {\mathrm{id}}_F$. 
Finally, suppose that $g = R^p$ with $p \ne 0$. Then $g(F) \ne F$. So $g|_F \ne 
{\mathrm{id}}|_F$.  \hfill $\square $ \medskip

\noindent \textbf{Lemma 3.3} For $j=0,1, \infty $ put 
$S^{(j)} = R^{n_j} U $. Then $S^{(j)}$ is a reflection in the closed ray 
${\ell }^j = \{ t{\mathrm{e}}^{i \, \pi  n_j/n} \in \C \setrule \, t \in OD \}$. The closed 
ray ${\ell }^0$  is the closure of the side $OD$ of the quadrilateral 
$Q = ODC\overline{D}$ in figure $5$.  \medskip

\noindent \textbf{Proof.} $S^{(j)}$ fixes every point on the closed ray ${\ell }^j$, because  
\begin{align}
S^{(j)}(\{ t{\mathrm{e}}^{i\, \pi n_j/n} \setrule \, t \in OD \} ) & = 
R^{n_j}( \{ t{\mathrm{e}}^{-i\, \pi n_j/n} \setrule \, t \in OD \} ) 
 =\{  t{\mathrm{e}}^{i\, \pi n_j/n} \setrule \, t \in OD \} . \notag 
\end{align}
Since $(S^{(j)})^2 = (R^{n_j}U)(R^{n_j}U) = R^{n_j}(UU)R^{-n_j}  = e$, it follows that $S^{(j)}$ is a reflection in the closed ray ${\ell }^j$.  
\hfill $\square $ \medskip

\noindent \textbf{Corollary 3.3a} For every $j=0,1, \infty$ and every 
$k \in \{ 0,1, \ldots , n-1\}$ let $S^{(j)}_k= R^kS^{(j)} R^{-k}$. Here 
$S^{(j)}_n = S^{(j)}_0 = S^{(j)}$, because 
$R^n =e$. Then $S^{(j)}_k$ is a reflection in the closed ray $R^k{\ell }^j$. \medskip 

\noindent \textbf{Proof.} This follows because $(S^{(j)}_k)^2 =R^k (S^{(j)})^2 R^{-k} = e$ and $S^{(j)}_k$ fixes every point on the closed ray $R^k {\ell }^j$, for   
\begin{align}
S^{(j)}_k\big( R^k (\{ t{\mathrm{e}}^{i \, \pi n_j /n}\setrule \, t \in OD \} ) \big) & = 
R^k S^{(j)}(\{ t{\mathrm{e}}^{ i\,  \pi n_j/n} \setrule \, t \in OD \} )\big) \notag \\ 
& \hspace{-1in} = R^k (\{ t{\mathrm{e}}^{i \, \pi n_j/n} \setrule \, t \in OD \} ) .  
\tag*{$\square $} 
\end{align}%

\noindent \textbf{Corollary 3.3b} For every $j =0,1, \infty$, every $S^{(j)}_k$ with 
$k=0,1, \ldots , n-1$, and every $g \in G$, we have $gS^{(j)}_kg^{-1} = S^{(j)}_r$ for a unique $r \in \{ 0, 1, \ldots , n-1 \}$. \medskip 

\noindent \textbf{Proof.} We compute. For every $k=0,1, \ldots , n-1$ we have 
\begin{equation}
RS^{(j)}_kR^{-1} = R(R^kS^{(j)}R^{-k})R^{-1} = R^{(k+1)}S^{(j)}R^{-(k+1)} = 
S^{(j)}_{k+1}
\label{eq-s3threea}
\end{equation}
and 
\begin{align}
US^{(j)}_kU^{-1} & = U(R^{(k+n_j)}UR^{-(k+n_j)})U = 
R^{-(k+n_j)} U R^{(k+n_j)} \notag \\
& = S^{(j)}_{-(k+2n_j)} = S^{(j)}_t, 
\label{eq-s3threeb}
\end{align}
where $t = -(k+2n_j) \bmod n$. Since $R$ and $U$ generate the group $G$, the 
corollary follows. \hfill $\square $ \medskip 

\noindent \textbf{Corollary 3.3c} For $j=0,1, \infty$ let $G^j$ be the group generated by the reflections $S^{(j)}_k$ for $k = 0,1, \ldots , n-1$. Then 
$G^j$ is a normal subgroup of $G$. \medskip 

\noindent \textbf{Proof.} Clearly $G^j$ is a subgroup of $G$. From equations 
(\ref{eq-s3threea}) and (\ref{eq-s3threeb}) it follows that 
$g S^{(j)}_k g^{-1} \in G^j$ 
for every $g \in G$ and every $k = 0 ,1 \ldots , n-1$, since $G$ is generated 
by $R$ and $U$. But $G^j$ is generated by the reflections $S^{(j)}_k$ for $k=0,1, \ldots , n-1$, that is, every $g' \in G^j$ may be written as 
$S^{(j)}_{i_1} \cdots S^{(j)}_{i_p}$, where for 
$\ell \in \{ 1, \ldots p \} $ we have $i_{\ell } \in \{ 0, 1, \ldots , n-1 \} $. So 
$gg'g^{-1} = g(S^{(j)}_{i_1} \cdots S^{(j)}_{i_p})g^{-1} = 
(gS^{(j)}_{i_1}g^{-1}) \cdots (gS^{(j)}_{i_p}g^{-1}) \in G^j$ for every 
$g \in G$, that is, $G^j$ is a normal subgroup of $G$. \hfill $\square $ \medskip  

As a first step in constructing ${\widetilde{S}}_{\mathrm{reg}}$ from the regular 
stellated $n$-gon $K^{\ast }$ we look at certain pairs of edges of $\mathrm{cl}(K^{\ast})$. We say two distinct closed edges $E$ and $E'$ of $\mathrm{cl}(K^{\ast })$ are \emph{adjacent} if and only if they intersect at a vertex of $\mathrm{cl}(K^{\ast })$. 
For $j=0,1, \infty$ let ${\mathcal{E}}^j$ be the set of unordered pairs of closed edges $E$ and $E'$ of 
$\mathrm{cl}(K^{\ast })$, that is, the edges $E$ and $E'$ are not adjacent and 
$E' = S^{(j)}_m(E)$ for some generator $S^{(j)}_m$ of $G^j$. Recall that for $x$ and $y$ in some set, the unordered pair $[x,y]$ is precisely one of the ordered pairs $(x,y)$ or 
$(y,x)$. Geometrically, two nonadjacent closed edges $E'$ and $E$ of $\mathrm{cl}(K^{\ast})$ are equivalent if and only if $E'$ is obtained from $E$ by reflection in the line $R^m {\ell }^j$ for some 
$m \in \{ 0, 1, \ldots , n-1 \} $. \medskip 

In figure 7, where $K^{\ast } = K^{\ast }_{1,1,4}$, parallel edges of $K^{\ast }$, 
which are labeled with the same letter, are $G^0$-equivalent. This is no coincidence. \medskip 

\par\noindent \hspace{1in}\begin{tabular}{l}
\includegraphics[width=175pt]{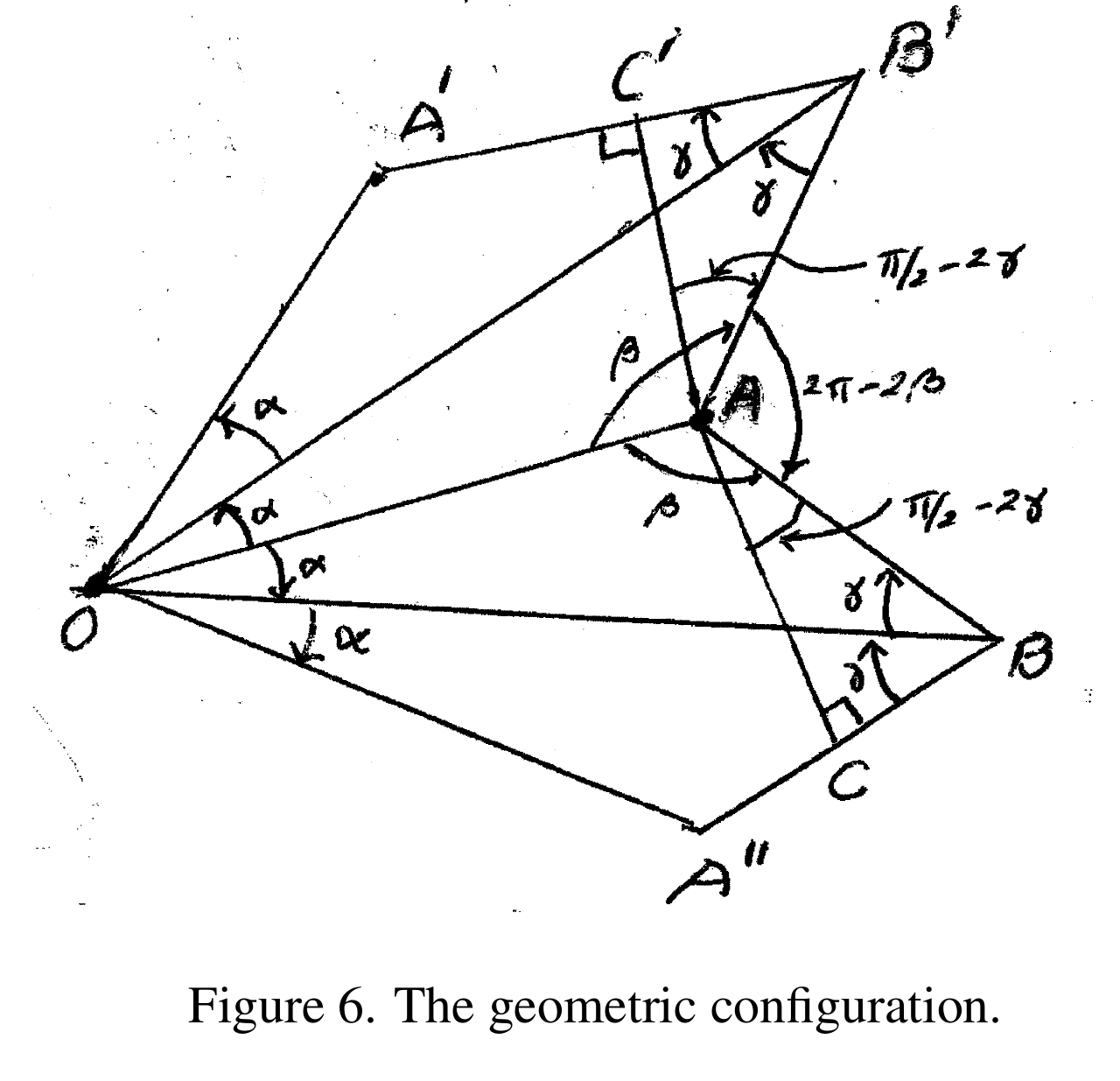}
\vspace{-.1in}
\end{tabular}
\medskip

\noindent \textbf{Lemma 3.4} Let $K^{\ast }$ be formed from the quadrilateral 
$Q = T \cup \overline{T}$, where $T$ is the isosceles rational triangle $T_{n_0n_0n_{\infty}}$ less its vertices. Then nonadjacent edges of $\partial \, \mathrm{cl}(K^{\ast })$ are 
$G^0$-equivalent if and only if they are parallel, see figure 6. \medskip 

\noindent \textbf{Proof.} In figure 6 let $OAB$ be the triangle $T$ with $\angle AOB = \alpha $, $\angle OAB = \beta $, and $\angle ABO = \gamma $. Let $OABA''$ be the 
quadrilateral formed by reflecting the triangle $OAB$ in its edge $OB$. The quadrilateral $OABA''$ reflected it its edge $OA$ is the quadrilateral $OAB'A'$. Let $AC'$ be perpendicular to $A'B'$ and $AC$ be perpendicular to $A''B$, see figure 6. Then $CAC'$ is a straight line if and only if $\angle C'AB' + \angle B'AB + \angle BAC = \pi $. By construction $\angle C'AB' = \angle BAC = \pi /2 - 2 \gamma $ and 
$\angle B'AB = 2\pi - 2 \beta $. So
\begin{align}
\pi & = 2(\ttfrac{\pi }{2} - 2\gamma ) + 2(\pi - \beta ) = 3\pi -2(\beta +\gamma ) - 2\gamma 
\notag \\
& = 3\pi - 2(\alpha +\beta +\gamma ) + 2(\alpha -\gamma ) = \pi +2 (\alpha - \gamma ), \notag 
\end{align}
if and only if $\alpha = \gamma $. Hence the edges $A''B$ and $A'B'$ are parallel if and only if the triangle $OAB$ is isosceles.  \hfill $\square $ \medskip

\noindent \textbf{Theorem 3.5} Let $K^{\ast }$ be the regular stellated $n$-gon 
formed from the rational quadrilateral $Q_{n_0n_1n_{\infty}}$ with $d_j = 
\mathrm{gcd}(n_j,n)$ for $j =0,1, \infty$. The $G$ orbit space  
formed by first identifiying equivalent edges of the regular stellated $n$-gon 
$K^{\ast }$ less $O$ and then acting on the identification 
space by the group $G$ is ${\widetilde{S}}_{\mathrm{reg}}$, which is a smooth $2$-sphere with $g$ handles, where 
$2g = n+2 -(d_0+d_1+d_{\infty})$ less some points corresponding to the 
image of the vertices of $\mathrm{cl}(K^{\ast })$.  \medskip %

Before we begin proving theorem 3.5 we consider the following special case. Let $K^{\ast } = K^{\ast }_{1,1,4}$ 
be a regular stellated hexagon formed by repeatedly rotating the quadrilateral $Q' =OD'C\overline{D'}$ by $R$ through an angle $2\pi /6$, see figure 7.  \medskip 

Let $G^0$ be the group generated by the reflections 
$S^{(0)}_k = R^k S^{(0)}R^{-k} = R^{2k+1}U$ for $k=0,1, \ldots , 5$. Here $S^{(0)} = RU$ is the reflection which leaves the closed ray ${\ell }^0 = 
\{ t {\mathrm{e}}^{i \pi /6} \setrule \, t \in OD'  \}$ fixed. Define an equivalence relation on 
$\mathrm{cl}(K^{\ast })$ \linebreak
\vspace{-.2in} \par \noindent \hspace{1in}\begin{tabular}{l}
\includegraphics[width=200pt]{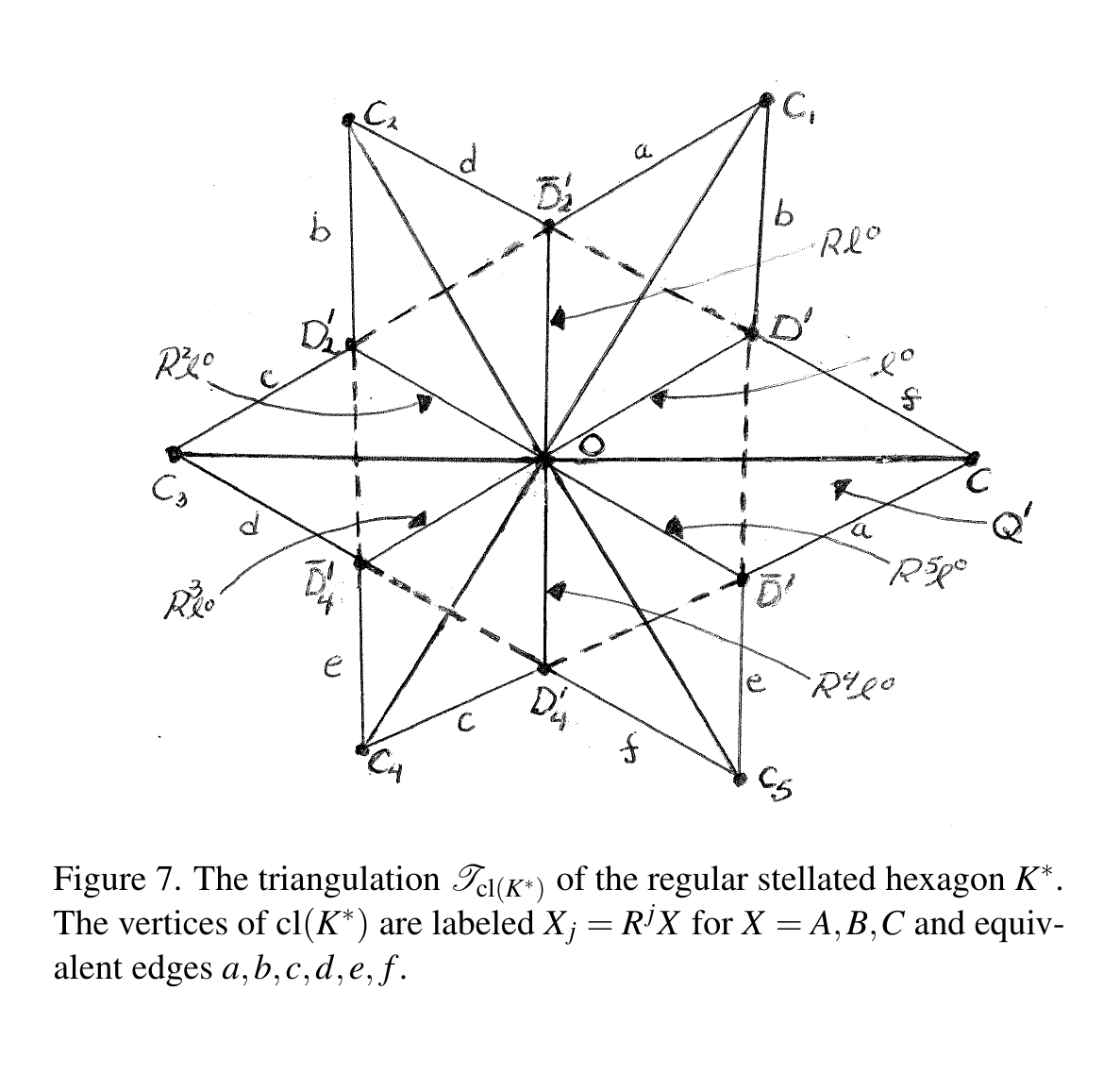}
\vspace{-.1in}
\end{tabular}  
\par \noindent by saying that two points $x$ and $y$ in $\mathrm{cl}(K^{\ast })$ are \emph{equivalent}, $x \sim y$, 
if and only if 1) $x$ and $y$ lie on $\partial \, \mathrm{cl}(K^{\ast})$ with $x$ on 
the closed edge $E$ and $y = S^{(0)}_m(x) \in S^{(0)}_m(E)$ for some reflection 
$S^{(0)}_m \in G^0$ or 2) if $x$ and $y$ lie in the interior of 
$\mathrm{cl}(K^{\ast })$ and $x =y$. Let $\mathrm{cl}(K^{\ast })^{\sim }$ be 
the space of equivalence classes and let 
\begin{equation}
\pi : \mathrm{cl}(K^{\ast }) \rightarrow \mathrm{cl}(K^{\ast })^{\sim } : p \mapsto [p]
\label{eq-s3threestar}
\end{equation}
be the identification map which sends a point $p \in \mathrm{cl}(K^{\ast })$ to the 
equivalence class $[p]$, which contains $p$. Give $\mathrm{cl}(K^{\ast })$ the topology induced from $\C $. Placing the quotient topology on 
$\mathrm{cl}(K^{\ast })^{\sim }$ turns it into a connected topological 
manifold without boundary. Let $K^{\ast }$ be $\mathrm{cl}(K^{\ast })$ less its vertices. The identification space 
$(K^{\ast } \setminus \{ O \})^{\sim} = \pi (K^{\ast } \setminus \{ O \})$ is a connected $2$-dimensional smooth manifold without boundary. \medskip 

Let $G= \langle R, U \setrule \, R^6 = e = U^2 \, \, \& \, \, RU = UR^{-1} \rangle $. The 
usual $G$-action 
\begin{displaymath}
G \times \mathrm{cl}(K^{\ast }) \subseteq G \times \C \rightarrow \mathrm{cl}(K^{\ast }) 
\subseteq \C :(g,z) \mapsto g(z)
\end{displaymath}
preserves equivalent edges of $\mathrm{cl}(K^{\ast })$ and is free on 
$K^{\ast } \setminus \{ O \}$. Hence it induces a $G$ action on 
$(K^{\ast } \setminus \{ O \} )^{\sim }$, which is free and proper. Thus its orbit map 
\begin{displaymath}
\sigma :  (K^{\ast } \setminus \{ O \} )^{\sim } \rightarrow 
(K^{\ast } \setminus \{ O \} )^{\sim }/G = {\widetilde{S}}_{\mathrm{reg}}: 
z \mapsto zG
\end{displaymath}
is surjective, smooth, and open. The orbit space ${\widetilde{S}}_{\mathrm{reg}} = 
\sigma ((K^{\ast } \setminus \{ O \} )^{\sim})$ is a connected $2$-dimensional smooth manifold. The identification space $(K^{\ast} \setminus \{ O \})^{\sim}$ has the orientation induced from an orientation of 
$K^{\ast } \setminus \{ O \} $, which comes from $\C $. So ${\widetilde{S}}_{\mathrm{reg}}$ 
has a complex structure, since each element of $G$ is a conformal mapping of $\C $ 
into itself. \medskip

Our aim is to specify the topology of ${\widetilde{S}}_{\mathrm{reg}}$. The regular stellated hexagon $K^{\ast } \setminus \{ O \}$ less the origin has 
a triangulation ${\mathcal{T}}_{K^{\ast } \setminus \{ O \}}$ made up of 
$12$ open triangles $R^j(\bigtriangleup OCD')$ and $R^j(\bigtriangleup OC{\overline{D}}')$ 
for $j =0,1, \ldots ,5$; $24$ open edges $R^j(OC)$, $R^j(O{\overline{D}}')$, 
$R^j(C{\overline{D}}')$, and $R^j(CD')$ for $j =0,1, \ldots , 5$; and $12$ vertices 
$R^j(D')$ and $R^j(C)$ for $j = 0,1, \ldots , 5$, see figure 7. \medskip 

Consider the set ${\mathcal{E}}^0$ of unordered pairs of equivalent closed edges of $\mathrm{cl}(K^{\ast })$, that is, ${\mathcal{E}}^0$ is the set $[E, S^{(0)}_k(E)]$ for $k=0,1, \ldots , 5$, where $E$ is a closed edge of $\mathrm{cl}(K^{\ast })$. Table 1 lists the elements of ${\mathcal{E}}^0$. \medskip 

\par \noindent \hspace{.1in} \begin{tabular}{lcl}
$a = \big[ \overline{D'}C, S^{(0)}_0( \overline{D'}C) = \overline{D'_2}C_1 \big] $ 
&\quad & 
$b = \big[ D'C_1, S^{(0)}_1(D'C_1) = D'_2C_2 \big] $ \\
\rule{0pt}{12pt}$d = \big[ \overline{D'_2}C_2, S^{(0)}_2(\overline{D'_2}C_2) = 
\overline{D'_4}C_3 \big] $ 
& \quad & 
$c = \big[ D'_2C_3, S^{(0)}_3(D'_2C_3) = D'_4C_4 \big] $ \\ 
\rule{0pt}{12pt}$e = \big[ \overline{D'_4}C_4, S^{(0)}_4(\overline{D'_4}C_4) = 
\overline{D'}C_5 \big] $ 
& \quad & 
$f = \big[ D'_4C_5, S^{(0)}_5(D'_4C_5) = D'C \big] $
\end{tabular} \medskip 

\par \noindent \hspace{.4in} \footnotesize{Table 1.} \parbox[t]{3.3in}{\footnotesize {Elements of the set ${\mathcal{E}}^0$. Here $D'_k =R^k(D')$ and 
$\overline{D'_k}= R^k(\overline{D'} )$ for $k=0,2,4$ and $C_k = R^k(C)$ for 
$k= \{ 0,1, \ldots , 5 \} $, see figure 7.} } %
\bigskip 

\noindent \normalsize $G$ acts on ${\mathcal{E}}^0$, namely, 
$g \cdot [E, S^{(0)}_k(E)] = [g(E), gS^{(0)}_kg^{-1}\big( g(E) \big) ] $, for $g \in G$. 
Since $G^0$ is the group generated by the reflections 
$S^{(0)}_k$, $k=0,1, \ldots , 5$,  it is a normal subgroup of $G$. 
Hence the action of $G$ on ${\mathcal{E}}^0$ restricts to an action of 
$G^0$ on ${\mathcal{E}}^0$ and permutes 
$G^0$-orbits in ${\mathcal{E}}^0$. Thus the set of $G^0$-orbits in ${\mathcal{E}}^0$ is $G$-invariant. \medskip 

We now look at the $G^0$-orbits on 
${\mathcal{E}}^0$. We compute the $G^0$-orbit of 
$d \in {\mathcal{E}}^0$ as follows. We have %
\begin{align}
(UR) \cdot d & = \big[ UR(\overline{D'_2}C_2), UR(S^{(0)}_2(\overline{D'_2}C_2) ) \big] = 
\big[ UR(\overline{D'_2}C_2), UR(\overline{D'_4}C_3)) \big] \notag 
\\
& = \big[ U(D'_2 C_3), U(D'_4C_4) \big] = 
\big[ \overline{D'_4}C_5, \overline{D'_2}C_2 \big] = d. \notag
\end{align}
Since 
\begin{align}
R^2 \cdot d & = R^2 \cdot \big[ \overline{D'_2}C_2, S^{(0)}_2(\overline{D'_2}C_2) \big]  = \big[ R^2(\overline{D'_2}C_2), R^2S^{(0)}_2R^{-2}(R^2(\overline{D'_2}C_2) ) \big] \notag 
\\
& = \big[ \overline{D'_4}C_4, S^{(0)}_4( \overline{D'_4} C_4) \big] = 
\big[ \overline{D'_4}C_4, \overline{D'}C_5 \big] = e \notag
\end{align}
and 
\begin{align}
R^4 \cdot d & =  
\big[ R^4(\overline{D'_4}C_2), R^4S^{(0)}_2R^{-4}(R^4(\overline{D'_2}C_2) )\big] \notag \\
&= \big[ \overline{D'}C , S^{(0)}_6(\overline{D'}C) \big] = 
\big[ \overline{D'}C, S^{(0)}_0(\overline{D'}C) \big] 
= \big[ \overline{D'}C, \overline{D'_2}C_1 \big] = a . \notag
\end{align}
So the $G^0$ orbit $G^0\cdot d$ of $d \in 
{\mathcal{E}}^0$ is $(G^0/\langle UR | \, (UR)^2  = e \rangle )\cdot d = H^0 \cdot d = \{ a, d, e \} $. Here $H^0 = \langle V = R^2 \setrule \, V^3 = e  \rangle $, 
since $G^0 = \langle V = R^2, UR \setrule \, V^3 = e = (UR)^2 \, \, \& \, \, V(UR) = (UR)V^{-1} \rangle $. Similarly, the $G^0$-orbit $G^0 \cdot f$ of $f \in {\mathcal{E}}^0$ is $H^0 \cdot f = \{ b, c, f \} $. Since 
$G^0 \cdot d \, \cup \, G^0 \cdot f = {\mathcal{E}}^0$, we have found all $G^0$-orbits on ${\mathcal{E}}^0$. 
The $G$-orbit of $OC$ is $R^j(OC)$ for $j=0,1, \ldots , 5$, since 
$U(OC) = OC$;  while the $G$-orbit of $OD'$ is $R^j(OD')$, $R^j(O\overline{D'})$ 
for $j=0,1, \ldots , 5$, since $U(OD') = O\overline{D'}$. \medskip 

Suppose that $B$ is an end point of the closed edge $E$ of $\mathrm{cl}(K^{\ast })$. 
Then $E$ lies in a unique $[E, S^{(0)}_m(E)]$ of ${\mathcal{E}}^0$. Let 
$G^0 \cdot [E, S^{(0)}_m(E)]$ be the $G^0$-orbit of 
$[E, S^{(0)}_m(E)]$. Then $g' \cdot B$ is an 
end point of the closed edge $g'(E)$ of $g'\cdot [E, S^{(0)}_m(E)] \in 
{\mathcal{E}}^0$ for every $g' \in G^0$. So $\mathcal{O}(B) = 
\{ g' \cdot B  \setrule \, g' \in G^0 \}$ the 
$G^0$-\emph{orbit} of the vertex $B$. It 
follows from the classification of $G^0$-orbits on ${\mathcal{E}}^0$ 
that $\mathcal{O}(D')= \{ D', D'_2, D'_4 \} $ and 
$\mathcal{O}(\overline{D'}) = \{ \overline{D'}, {\overline{D'}}_2, {\overline{D'}}_4 \} $ 
are $G^0$-orbits of the vertices of $\mathrm{cl}(K^{\ast })$, which are 
permuted by the action of $G$ on ${\mathcal{E}}^0$. 
Also $\mathcal{O}(C) = \{ C,C_1, \ldots , C_5 \}$ and $\mathcal{O}(D' \& \overline{D'}) = 
\{ D', \overline{D'}, D'_2, {\overline{D'}}_2, D'_4, {\overline{D'}}_4 \} $ are $G$-orbits 
of vertices of $\mathrm{cl}(K^{\ast })$, which are end points of the $G$-orbit of 
the rays $OC$ and $OD'$, respectively. \medskip  

\vspace{-.1in}\par \noindent \hspace{.85in}\begin{tabular}{l}
\includegraphics[width=200pt]{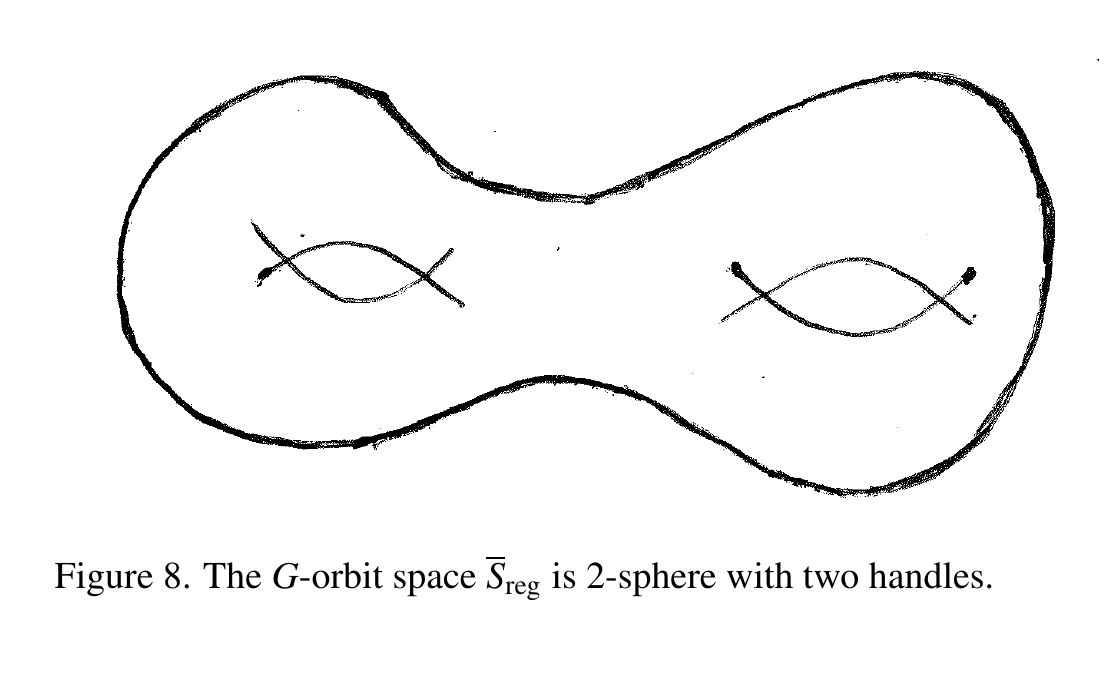}
\vspace{-.1in}
\end{tabular} \medskip 

To determine the topology of the $G$ orbit space ${\widetilde{S}}_{\mathrm{reg}}$ we find a triangulation of ${\widetilde{S}}_{\mathrm{reg}}$. Note that the triangulation 
${\mathcal{T}}_{K^{\ast } \setminus \{ O \}}$ of 
$K^{\ast } \setminus \{ O \}$, illustrated in figure 7, is $G$-invariant. Its image under the 
identification map $\pi $ is a $G$-invariant triangulation 
${\mathcal{T}}_{(K^{\ast } \setminus \{ O \}}$ of 
$(K^{\ast } \setminus \{ O \})^{\sim }$. After identification of equivalent edges, 
each vertex $\pi (v)$, each open edge $\pi (E)$, having $\pi (O)$ as an end point, or each open edge $\pi ([F,F'])$, where $[F,F']$ is a pair of equivalent edges of $\mathrm{cl}(K^{\ast })$, and each open triangle $\pi (T)$ in ${\mathcal{T}}_{(K^{\ast } \setminus \{ O \})^{\sim }}$ lies in a unique $G$ orbit. It follows that 
$\sigma (\pi (v))$, $\sigma (\pi (E))$ or $\sigma (\pi ([F,F']))$, and 
$\sigma (\pi (T))$ is a vertex, an open edge, and an open triangle, respectively, of a triangulation 
${\mathcal{T}}_{{\widetilde{S}}_{\mathrm{reg}}} = 
\sigma ({\mathcal{T}}_{(K^{\ast } \setminus \{ O \})^{\sim }})$ of 
${\widetilde{S}}_{\mathrm{reg}}$. The triangulation 
${\mathcal{T}}_{{\widetilde{S}}_{\mathrm{reg}}}$ has $4$ vertices, 
corresponding to the $G$ orbits $\sigma (\pi (\mathcal{O}(D')))$, 
$\sigma (\pi (\mathcal{O}(\overline{D'})))$, 
$\sigma (\pi (\mathcal{O}(C)))$, and $\sigma (\pi (\mathcal{O}(D' \& \overline{D'})))$; 
$18$ open edges 
corresponding to $\sigma (\pi (R^j(OC)))$, $\sigma (\pi (R^j(OD')))$, 
and $\sigma (\pi (R^j(CD')))$ for $j=0,1, \ldots , 5$; and $12$ open triangles 
$\sigma (\pi (R^j(\bigtriangleup OCD')))$ and 
$\sigma (\pi (R^j(\bigtriangleup OC\overline{D'})))$ for $j=0,1, \ldots , 5$. \linebreak 
Thus the Euler characteristic $\chi ({\widetilde{S}}_{\mathrm{reg}})$ of
${\widetilde{S}}_{\mathrm{reg}}$ is 
$4 - 18 + 12 = -2$. Since ${\widetilde{S}}_{\mathrm{reg}}$ is a $2$-dimensional smooth 
real manifold, $\chi ({\widetilde{S}}_{\mathrm{reg}}) = 2 - 2g$, where $g$ is the genus of 
${\widetilde{S}}_{\mathrm{reg}}$. Hence $g=2$. So ${\widetilde{S}}_{\mathrm{reg}}$ is a smooth $2$-sphere with $2$ handles,  less a finite number of points, which lies in a compact topological space $\widetilde{S} = 
\mathrm{cl}(K^{\ast })^{\sim }/G$, that is its closure.  
\hfill $\square $ \medskip %

\noindent \textbf{Proof of theorem 3.5} We now begin the construction of ${\widetilde{S}}_{\mathrm{reg}}$ by 
identifying equivalent edges of $\mathrm{cl}(K^{\ast })$. Let $[E, S^{(0)}_m(E)]$ be an unordered pair of equivalent closed edges of $\mathrm{cl}(K^{\ast })$. We say that $x$ and $y$ in $\mathrm{cl}(K^{\ast })$ 
are \emph{equivalent}, $x \sim y$, if 1) $x$ and $y$ lie in 
$\partial \, \mathrm{cl}(K^{\ast })$ with $x \in E$ and $y = S^{(0)}_m(x) \in S^{(0)}_m(E)$ for some $m \in \{ 0,1, \ldots , n-1 \}$ or 
2) $x$ and $y$ lie in $\mathrm{int}\, \mathrm{cl}(K^{\ast })$ and $x = y$. The relation 
$\sim $ is an equivalence relation on $\mathrm{cl}(K^{\ast})$. Let 
$\mathrm{cl}(K^{\ast })^{\sim }$ be the set of equivalence classes and let 
\begin{equation}
\pi : \mathrm{cl}(K^{\ast }) \rightarrow  \mathrm{cl}(K^{\ast })^{\sim} : p \mapsto [p]
\label{eq-s3fourstar}
\end{equation} 
be the map which sends $p$ to the equivalence class $[p]$, that 
contains $p$. Compare this argument with that of Richens and Berry \cite{richens-berry}. 
Give $\mathrm{cl}(K^{\ast })$ the topology induced from $\C $ and put the quotient topology on $\mathrm{cl}(K^{\ast })^{\sim }$. \medskip 
 
\noindent \textbf{Claim 3.6}  Let $K^{\ast }$ be $\mathrm{cl}(K^{\ast })$ less its vertices. Then 
$(K^{\ast } \setminus \{ O \} )^{\sim } = \pi (K^{\ast } \setminus \{ O \} )$ is a smooth manifold. Also 
$\mathrm{cl}(K^{\ast })^{\sim }$ is a topological manifold.\medskip 

\noindent \textbf{Proof.} To show that $(K^{\ast } \setminus \{ O \} )^{\sim }$ is a smooth manifold, let $E_{+}$ be an open edge of $K^{\ast }$. For $p_{+} \in E_{+}$ let 
$D_{p_{+}}$ be a disk in $\C $ with center at $p_{+}$, which does not contain a vertex of 
$\mathrm{cl}(K^{\ast })$. Set $D^{+}_{p_{+}} = K^{\ast } \cap D_{p_{+}}$. 
Let $E_{-}$ be an open edge of $K^{\ast }$, which is equivalent to 
$E_{+}$ via the reflection $S^{(0)}_m$, that is, 
$[\mathrm{cl}(E_{+}), \mathrm{cl}(E_{-}) = S^{(0)}_m(\mathrm{cl}(E_{+}))] \in 
{\mathcal{E}}^0$ is an unordered pair of $S^{(0)}_m$ equivalent edges. 
Let $p_{-} = S^{(0)}_m(p_{+})$ and set $D^{-}_{p_{-}} 
= S^{(0)}_m(D^{+}_{p_{+}})$. Then $V_{[p]} = \pi (D^{+}_{p_{+}} \cup D^{-}_{p_{-}}) $ is an open neighborhood of $[p] = [p_{+}] = [p_{-}]$ in $(K^{\ast } \setminus \{ O \})^{\sim }$, which is a smooth $2$-disk, since the identification mapping $\pi $ 
is the identity on $\mathrm{int}\, K^{\ast }$. It follows that 
$(K^{\ast } \setminus \{ O \} )^{\sim }$ is a smooth $2$-dimensional manifold without boundary. \medskip

We now handle the vertices of $\mathrm{cl}(K^{\ast })$. Let $v_{+}$ be a vertex of 
$\mathrm{cl}(K^{\ast })$ and set $D_{v_{+}} =
\widetilde{D} \cap \mathrm{cl}(K^{\ast })$, where $\widetilde{D}$ is a disk in $\C $ 
with center at the vertex $v_{+}= r_0{\mathrm{e}}^{i\pi {\theta }_0}$. The map 
\begin{displaymath}
W_{v_{+}}: D_{+} \subseteq \C \rightarrow D_{v_{+}} \subseteq \C: 
r{\mathrm{e}}^{i\pi \theta } \mapsto |r - r_0| {\mathrm{e}}^{i \pi s (\theta - {\theta }_0)} 
\end{displaymath}
with $r \ge 0$ and $0 \le \theta \le 1$ is a homeomorphism, which sends the wedge with angle $\pi $ to the wedge with angle $\pi s$. The latter wedge is formed by the closed edges $E'_{+}$ and $E_{+}$ of $\mathrm{cl}(K^{\ast })$, which are adjacent at the vertex 
$v_{+}$ such that ${\mathrm{e}}^{i\pi s}E'_{+} = E_{+}$ with the edge $E'_{+}$ being 
swept out through $\mathrm{int}\, \mathrm{cl}(K^{\ast })$ during its rotation to the edge $E_{+}$. Because $\mathrm{cl}(K^{\ast })$ is a rational regular stellated 
$n$-gon, the value of $s$ is a rational number for each vertex of 
$\mathrm{cl}(K^{\ast })$. Let $E_{-} = S^{(0)}_m(E_{+})$ be 
an edge of $\mathrm{cl}(K^{\ast })$, which is equivalent to $E_{+}$ and set 
$v_{-} = S(v_{+})$. Then $v_{-}$ is a vertex of $\mathrm{cl}(K^{\ast })$, which is 
the center of the disk $D_{v_{-}} = S^{(0)}_m(D_{v_{+}})$. Set $D_{-} = 
{\overline{D}}_{+}$. Then $D = D_{+} \cup D_{-}$ is a disk in $\C $. The map 
$W: D \rightarrow \pi (D_{v_{+}} \cup D_{v_{-}} )$, where $W|_{D_{+}} = \pi \comp W_{v_{+}}$ and $W|_{D_{-}} = \pi \comp S^{(0)}_m \comp W_{v_{+}} \comp 
{\mbox{}}^{\overline{\rule{5pt}{0pt}}}$, is a homeomorphism of $D$ into a neighborhood 
$\pi ( D_{v_{+}} \cup D_{v_{-}})$ of $[v] = [v_{+}] = [v_{-}]$ in 
$\mathrm{cl}(K^{\ast })^{\sim}$. Consequently, the identification space 
$\mathrm{cl}(K^{\ast })^{\sim }$ is a topological manifold. \hfill $\square $ \medskip  

We now describe a triangulation of $K^{\ast } \setminus \{ O \} $. Let 
$T' = T_{1,n_1, n-(1+n_1)}$ be the open rational triangle $\bigtriangleup OCD'$ 
with vertex at the origin $O$, longest side $OC$ on the real axis, and 
interior angles $\frac{1}{n} \pi $, $\frac{n_1}{n} \pi $, and $\frac{n-1-n_1}{n} \pi $.
 Let $Q'$ be the quadrilateral $T' \cup \overline{T'}$. Then $Q'$ is a subset of the quadrilateral $Q = ODC\overline{D}$, see figure 5. Moreover $K^{\ast } = 
\bigcup^{n-1}_{\ell = 0}R^{\ell }(Q')$. The $2n$ triangles $\mathrm{cl}(R^j(T')) \setminus 
\{ O \}$ and $\mathrm{cl}(R^j(\overline{T'}))\setminus \{ O \}$ with $j = 0,1, \ldots , n-1$ form a triangulation ${\mathcal{T}}_{K^{\ast } \setminus \{ O \}}$ of $K^{\ast } \setminus \{ O \} $ with $2n$ vertices $R^j(C)$ and $R^j(D')$ for $j = 0,1, \ldots , n-1$; $4n$ open edges $R^j(OC)$, $R^j(OD')$, $R^j(CD')$, and $R^j(C\overline{D'})$ for 
$j =0,1, \ldots , n-1$; and $2n$ open triangles $R^j(T')$, $R^j(\overline{T'})$ with 
$j=0,1,\ldots , n-1$. The image of the triangulation 
${\mathcal{T}}_{K^{\ast } \setminus \{ O \}}$ under the identification map 
$\pi $ (\ref{eq-s3fourstar}) is a triangulation 
${\mathcal{T}}_{K^{\ast } \setminus \{ O \})^{\sim}}$ of the 
identification space $(K^{\ast } \setminus \{ O \} )^{\sim}$. \medskip 

The action of $G$ on $\mathrm{cl}(K^{\ast })$ preserves the set of unordered pairs of 
$S^{(j)}_m$ equivalent edges of $\mathrm{cl}(K^{\ast })$ for $j=0,1, \infty $.  Hence $G$ induces an action on 
$\mathrm{cl}(K^{\ast })^{\sim }$, which is proper, since $G$ is finite. The $G$ action 
is free on $K^{\ast} \setminus \{ O \} $ and thus on $(K^{\ast } \setminus \{ O \})^{\sim }$ 
by lemma A2. We have proved \medskip 

\noindent \textbf{Lemma 3.7} The $G$-orbit space $\widetilde{S} = 
\mathrm{cl}(K^{\ast })^{\sim }/G$ is a compact connected topological manifold with 
${\widetilde{S}}_{\mathrm{reg}} = (K^{\ast }\setminus \{ O \})^{\sim }/G$ being a smooth manifold. Let \linebreak    
\begin{displaymath}
\sigma : \mathrm{cl}(K^{\ast })^{\sim } \rightarrow  \widetilde{S} = 
\mathrm{cl}(K^{\ast })^{\sim }/G : 
z \mapsto zG.
\end{displaymath}
Then $\sigma $ is the $G$ orbit map, which is surjective, continuous, and open. 
The restriction of the map $\sigma $ to 
$K^{\ast } \setminus \{ O \}$ has image ${\widetilde{S}}_{\mathrm{reg}}$ and is a smooth open mapping. \medskip 

We now determine the topology of the orbit space ${\widetilde{S}}_{\mathrm{reg}}$. 
For $j =0,1, \infty$ and $\ell = 0, 1, \ldots , d_j-1$ let $A^j_{\ell }$ 
be an end point of a closed edge $E$ of $\mathrm{cl}(K^{\ast })$, which lies on the unordered pair $[E, S^{(j)}_{\ell }(E)] \in {\mathcal{E}}^j$. Then 
$H^j \cdot A^{(j)}_{\ell }$ is an end point of the edge $H^j\cdot E$ of the unordered pair 
$H^j\cdot [E, S^{(j)}_{\ell }(E)]$ of ${\mathcal{E}}^j$. See appendix A for the definition of the 
group $H_j$. Fix $j$. The 
sets $\mathcal{O}(A^{(j)}_{\ell }) = \{ H^j \cdot A^{(j)}_{\ell } \} $ 
with $\ell = 0,1, \ldots , d_j-1$ are permuted by $G$. The action of $G$ on 
$K^{\ast } \setminus \{ O \} $ preserves the set of open edges of the triangulation 
${\mathcal{T}}_{K^{\ast }\setminus \{ O \} } $. 
There are $3n$-orbits: $R^j(OC)$; $R^j(O\overline{D'})$, since 
$OD' = R(O\overline{D'})$; and $R^j(CD)$, since $C\overline{D'} = U(CD)$ 
for $j=0,1, \ldots , n-1$. So the image of the triangulation 
${\mathcal{T}}_{K^{\ast }\setminus \{ O \} }$ under the continuous open map  
\begin{equation}
\mu= \sigma \comp \pi |_{K^{\ast } \setminus \{ O \}}: K^{\ast }\setminus \{ O \}  \rightarrow 
{\widetilde{S}}_{\mathrm{reg}}
\label{eq-s3fourdagger}
\end{equation}
is a triangulation ${\mathcal{T}}_{ {\widetilde{S}}_{\mathrm{reg}}}$ of the $G$-orbit space 
${\widetilde{S}}_{\mathrm{reg}}$ with $d_0+d_1+d_{\infty}$ vertices 
$\mu (\mathcal{O}(A^{(j)}_{\ell }))$, where $j =0,1, \infty $ and $\ell = 0,1, \ldots , d_j-1$; $3n$ open edges $\mu (R^j(OC))$, $\mu (R^j(O\overline{D'}))$, and 
$\mu (R^j(CD))$ for $j=0,1, \ldots , n-1$; and $2n$ open triangles 
$\mu (R^j(T'))$ and $\mu (R^j(\overline{T'}))$ for $j =0,1, \ldots n-1$. Thus the Euler characteristic $\chi ({\widetilde{S}}_{\mathrm{reg}})$ of ${\widetilde{S}}_{\mathrm{reg}}$ is 
$d_0+d_1+d_{\infty} - 3n +2n = d_0+d_1+d_{\infty} - n$. But 
${\widetilde{S}}_{\mathrm{reg}}$ is a smooth manifold. So 
$\chi ({\widetilde{S}}_{\mathrm{reg}}) = 2 -2g$, 
where $g$ is the genus of ${\widetilde{S}}_{\mathrm{reg}}$. Hence 
$2g = n+2 -(d_0+d_1+d_{\infty})$. Compare this argument with that of Weyl \cite[p.174]{weyl}.  This 
proves theorem 3.5. \hfill $\square $  \medskip 

Since the quadrilateral $Q$ is a fundamental domain for the action of $G$ on 
$K^{\ast }$, the $G$ orbit map $\overline{\mu } = \sigma \comp \pi : K^{\ast } \subseteq \C 
\rightarrow \widetilde{S}$ restricted to $Q$ is a bijective continuous open mapping. 
But ${\delta }_Q: \mathcal{D} \subseteq {\mathcal{S}}_{\mathrm{reg}} 
\rightarrow Q \subseteq \C $ is a bijective continous open mapping of the 
fundamental domain $\mathcal{D}$ of the $\mathcal{G}$ action on 
$\mathcal{S}$. Consequently, the $\mathcal{G}$ orbit space $\mathcal{S}$ is 
homeomorphic to the $G$ orbit space $\widetilde{S}$. 
The mapping $\overline{\mu }$ is holomorphic except possibly at $0$ and 
the vertices of $\mathrm{cl}(K^{\ast })$. So the mapping 
$\overline{\mu } \comp {\delta }_{K^{\ast }} : {\mathcal{S}}_{\mathrm{reg}} \rightarrow 
{\widetilde{S}}_{\mathrm{reg}}$ is a holomorphic diffeomorphism. \medskip

\section{An affine model of ${\widetilde{S}}_{\mathrm{reg}}$}

We construct an affine model of the Riemann surface 
${\widetilde{S}}_{\mathrm{reg}}$. \medskip  

We return to the regular stellated $n$-gon $K^{\ast } = K^{\ast}_{n_0n_1n_{\infty}}$, which is formed from the quadrilateral $Q = Q_{n_0n_1n_{\infty}}$ less its vertices. Repeatedly reflecting in the edges of ${K^{\ast}}$ and then in the edges of the resulting reflections of 
${K^{\ast}}$ et cetera,  we obtain a covering of $\C \setminus {\mathbb{V}}^{+}$ by certain translations of $K^{\ast}$. Here ${\mathbb{V}}^{+}$ is the union of the translates of the vertices of $\mathrm{cl}(K^{\ast})$ and its center $O$. Let $\mathfrak{T}$ be the group generated by these translations. The semidirect product $\mathfrak{G} = G \ltimes \mathfrak{T}$ acts freely, properly and transitively on $\C \setminus {\mathbb{V}}^{+}$. It preserves equivalent edges of $\C \setminus {\mathbb{V}}^{+}$ and it acts freely and properly on 
$(\C \setminus {\mathbb{V}}^{+})^{\sim }$, the space formed by identifying equivalent 
edges in $\C \setminus {\mathbb{V}}^{+}$. The orbit space 
$(\C \setminus {\mathbb{V}}^{+})^{\sim }/\mathfrak{G}$ is holomorphically diffeomorphic to 
${\widetilde{S}}_{\mathrm{reg}}$ and is the desired affine model of ${\widetilde{S}}_{\mathrm{reg}}$. We now justify these assertions. \medskip

First we determine the group $\mathcal{T}$ of translations. \medskip  

\noindent \textbf{Lemma 4.1} Each of the $2n$ sides of the regular stellated $n$-gon $K^{\ast}$ is perpendicular to exactly one of the directions 
\begin{equation}
{\mathrm{e}}^{[\frac{1}{2} - \frac{n_1}{n}  + 2j \frac{1}{n}] \pi  i} \, \, \, \mathrm{or} \, \, \,  {\mathrm{e}}^{[-\frac{1}{2} -\frac{1}{n} + \frac{n_1}{n} + (2j+1)\frac{1}{n}] \pi i},  
\label{eq-s3threedoublestar}
\end{equation}
for $j =0, 1, \ldots , n-1$. \medskip  

\noindent \textbf{Proof.} From figure 9 we have $\angle D'CO = \frac{n_1}{n} \pi $. So 
$\angle COH = \frac{1}{2}\pi - \frac{n_1}{n}\pi $. Hence the line ${\ell }_0$, containing 
the edge $CD'$ of $K^{\ast}$, is perpendicular to the direction 
${\mathrm{e}}^{[\frac{1}{2} - \frac{n_1}{n}]\pi }$. Since $\bigtriangleup CO\overline{D'} $ 
is the reflection of $\bigtriangleup COD'$ in the line segment $OC$, the line 
${\ell }_1$, containing the edge $C\overline{D'}$ of $K^{\ast}$, is perpendicular to 
the direction ${\mathrm{e}}^{[-\frac{1}{2} + \frac{n_1}{n}]\pi }$. 
Because the regular stellated $n$-gon ${K^{\ast}}$ is formed by repeatedly rotating the quadrilateral $Q'= OD'C\overline{D'}$ through an angle $\frac{2\pi }{n}$, we find that equation (\ref{eq-s3threedoublestar}) holds. \hfill $\square $ 
\par\noindent \hspace{1in}\begin{tabular}{l}
\includegraphics[width=200pt]{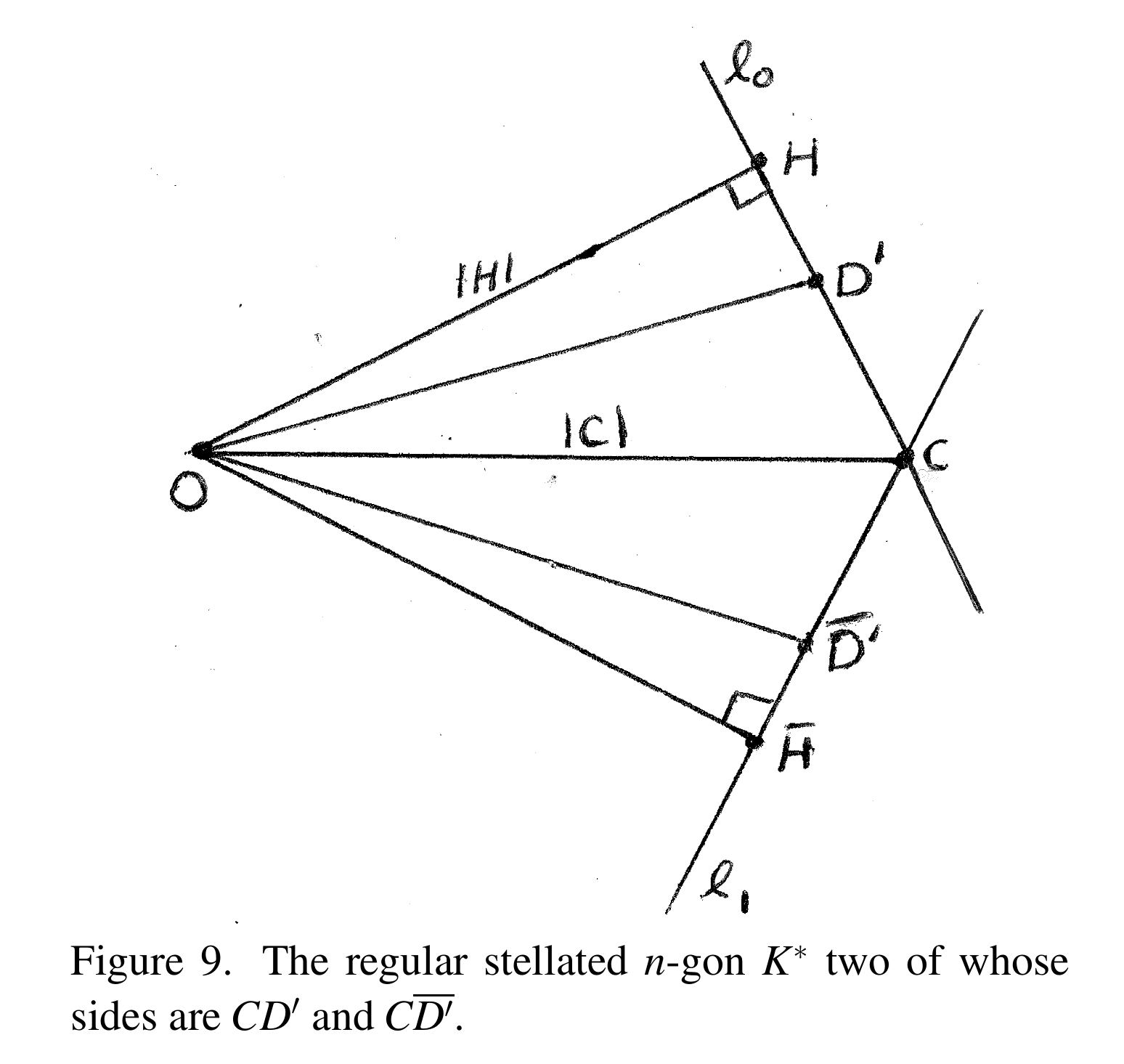} \\ 
\vspace{-.2in}
\end{tabular}  

Since $\angle COH = \frac{1}{2} \pi -\frac{n_1}{n} \pi $, it follows that 
$|H| = |C| \sin \pi \frac{n_1}{n}$ is the distance from the center $O$ of $K^{\ast }$ 
to the line ${\ell }_0$ containing the side $CD'$, or to the line 
${\ell }_1$ containing the side $C \overline{D'}$. So 
$u_0 = (|C| \sin \pi \frac{n_1}{n}){\mathrm{e}}^{ [\frac{1}{2} -\frac{n_1}{n}] \pi i}$ 
is the closest point $H$ on ${\ell }_0$ to $O$ and 
$u_1 = (|C| \sin \pi \frac{n_1}{n} ) {\mathrm{e}}^{ [ -\frac{1}{2} +\frac{n_1}{n} ] \pi i} $ 
is the closest point $\overline{H}$ on ${\ell }_1$ to $O$. Since the regular stellated $n$-gon ${K^{\ast}}$ is formed by repeatedly rotating the quadrilateral $Q' = OD'C\overline{D'}$ through an angle $\frac{2\pi }{n}$, the point 
\begin{equation}
u_{2j } = R^ju_0 = 
(|C| \sin \pi \frac{n_1}{n} ){\mathrm{e}}^{[ \frac{1}{2} -\frac{n_1}{n} +2j \frac{1}{n} ] \pi i}  
\label{eq-s4zeroa}
\end{equation}
lies on the line ${\ell }_{2j}= R^j{\ell }_0$, which contains the edge 
$R^j(CD')$ of ${K^{\ast}}$; while 
\begin{equation}
u_{2j+1} = R^ju_1 = 
(|C| \sin \pi \frac{n_1}{n} ){\mathrm{e}}^{ [ -\frac{1}{2} + \frac{n_1}{n} -\frac{1}{n} 
+ (2j +1)\frac{1}{n} ] \pi i}
\label{eq-s4zerob}
\end{equation}
lies on the line ${\ell }_{2j+1} = R^j{\ell }_1$, which contains the edge 
$R^j(C \overline{D'} )$ of ${K^{\ast}}$ for every $j \in \{ 0,1, \ldots , n-1 \}$. 
Also the line segments $Ou_{2j}$ and $Ou_{2j+1}$ are perpendicular to the line 
${\ell }_{2j}$ and ${\ell }_{2j+1}$, respectively, for $j \in \{ 0,1, \ldots , n-1 \} $.  \medskip   

\noindent \textbf{Corollary 4.1a} For $j =0, 1, \ldots , n-1$ we have 
\begin{equation}
\overline{u_{2j}} = u_{2(n-j) +1} \, \, \, \mathrm{and} \, \, \, 
\overline{u_{2j+1}} = u_{2(n-j)} . 
\label{eq-s4zeroc}
\end{equation}

\noindent \textbf{Proof.} We compute. From (\ref{eq-s4zeroa}) it follows that 
\begin{align}
\overline{u_{2j}} & = U(u_{2j}) = UR^j(u_0) = R^{-j}(U(u_0)) \notag \\
& = R^{-j}(u_1) = R^{n-j}(u_1) = u_{2(n-j) +1}, \quad \mbox{using (\ref{eq-s4zerob});} \notag 
\end{align}
while from (\ref{eq-s4zerob}) we get 
\begin{align}
\overline{u_{2j+1}} & = U(u_{2j+1}) = UR^j(u_1) = R^{-j}(U(u_1)) 
 = R^{n-j}(u_0) = u_{2(n-j)}. \tag*{$\square $}
\end{align}

\noindent \textbf{Corollary 4.1b} For $\ell $, $j \in \{ 0, 1, \ldots , 2n-1 \} $ we have 
\begin{equation}
u_{(\ell +2j) \bmod 2n} = R^j u_{\ell }.
\label{eq-s4two}
\end{equation}

\noindent \textbf{Proof.} If $\ell =2k$, then $u_{\ell } = R^k u_0$, by definition. So 
\begin{displaymath}
R^ju_{\ell } = R^{j+k}u_0 = u_{(2k+2j) \bmod 2n} = u_{(\ell +2j) \bmod 2n} . 
\end{displaymath}
If $\ell  = 2k+1$, then $u_{\ell } = R^ku_1$, by definition. So 
\begin{align} 
R^ju_{\ell }= R^{j+k}u_1 = u_{(2(j+k)+1) \bmod 2n } = u_{(\ell +2j) \bmod 2n} . 
\tag*{$\square $}
\end{align}

For $j=0,1, \ldots , 2n-1$ let ${\tau }_j$ be the translation 
\begin{equation}
{\tau }_j: \C \rightarrow \C : z \mapsto z + 2u_j. 
\label{eq-s4three}
\end{equation}

\noindent \textbf{Corollary 4.1c} For $k$, $j \in \{ 0, 1, \ldots , 2n-1 \} $ we have 
\begin{equation}
{\tau }_{(k+2j) \bmod 2n} \comp R^j = R^j \comp {\tau }_k. 
\label{eq-s4four}
\end{equation}

\noindent \textbf{Proof.} For every $z \in \C $, we have 
\begin{align}
{\tau }_{(k +2j) \bmod 2n} (z) & = z + 2  u_{(k+2j) \bmod 2n }, \quad 
\mbox{using (\ref{eq-s4three})} \notag \\
& = z + 2 R^j u_k \quad \mbox{by (\ref{eq-s4two})} \notag \\
& = R^j(R^{-j}z + 2 u_k) = R^j \comp {\tau }_k (R^{-j}z). \tag*{$\square $}
\end{align}

Reflecting the regular stellated $n$-gon $K^{\ast} $ in its edge $CD'$ contained in 
${\ell }_0$ gives a congruent regular stellated $n$-gon $K^{\ast }_0$ with the center $O$ of $K^{\ast }$ becoming the center $2u_0$ of $K^{\ast }_0$. \medskip

\noindent \textbf{Lemma 4.2} The collection of all the centers of the regular stellated $n$-gons formed by reflecting ${K^{\ast}}$ in its edges and then reflecting in the edges of the reflected regular stellated $n$-gons et cetera is 
\begin{align}
\{ {\tau }^{{\ell }_0}_0 \comp \cdots \comp {\tau }^{{\ell }_{2n-1}}_{2n-1}(0) \in \C \setrule 
({\ell }_0, \ldots , {\ell }_{2n-1}) \in ({\Z }_{\ge 0})^{2n} \} 
& = \notag \\
&\hspace{-2in} 
= \big\{ 2\, \hspace{-20pt} \sum^{\infty}_{\hspace{10pt}{\ell }_0, \ldots , {\ell }_{2n-1} = 0 }
\hspace{-20pt}\big( {\ell }_0 u_0 + \cdots {\ell }_{2n-1}u_{2n-1} \big) \big\}, \notag 
\end{align}
where for $j=0, 1, \ldots , 2n-1$ we have 
\begin{displaymath}
{\tau }^{{\ell }_j}_j = 
\overbrace{{\tau }_j \comp \cdots \comp {\tau }_j}^{{\ell}_j}:\C \rightarrow \C: z \mapsto 
z + 2{\ell }_ju_j. 
\end{displaymath} 

\noindent \textbf{Proof.} For each $k_0 =0, 1, \ldots , 2n-1$ the center of the $2n$ regular stellated congruent 
$n$-gon $K^{\ast}_{k_0}$ formed by reflecting in an edge of 
$K^{\ast}$ contained in the line ${\ell }_{k_0}$ is ${\tau }_{k_0}(0) = 2u_{k_0}$. Repeating the reflecting process in each edge of $K^{\ast}_{k_0}$ gives $2n$ congruent regular stellated $n$-gons $K^{\ast}_{k_0k_1}$ with center at ${\tau}_{k_1}\big( {\tau }_{k_0}(0) \big) = 2(u_{k_1} + u_{k_0})$, where $k_1=0,1, \ldots 2n-1$. Repeating this construction proves the lemma. 
\hfill $\square $ \medskip 

The set $\mathbb{V}$ of vertices of the regular stellated $n$-gon $K^{\ast}$ is 
\begin{displaymath}
\{ V_{2j} = C{\mathrm{e}}^{2j(\frac{1}{n} \pi \, i )}, \, V_{2j+1} = D'{\mathrm{e}}^{(2j+1)(\frac{1}{n} \pi \, i)} \, \, \mbox{for $0 \le j \le n-1$}  \} ,
\end{displaymath}
see figure 5. Clearly the set $\mathbb{V}$ is $G$ invariant.  \medskip 

\noindent \textbf{Corollary 4.2a} The set 
\begin{align}
{\mathbb{V}}^{+} & = 
\{ v_{{\ell }_0 \cdots {\ell }_{2n-1}} = 
{\tau }^{{\ell }_0}_0 \comp \cdots \comp {\tau }^{{\ell }_{2n-1}}_{2n-1}(V) 
\setrule \notag \\
& \hspace{.5in}  V \in \mathbb{V} \cup \{ 0 \} \, \, \& \, \, 
({\ell }_0, \ldots , {\ell }_{2n-1}) \in ({\Z }_{\ge 0})^{2n} \} 
\label{eq-s4threestardot}
\end{align}
is the collection of vertices and centers of the congruent regular stellated 
$n$-gons $K^{\ast }$, $K^{\ast }_{k_1}$, $K^{\ast}_{k_0k_1}, \ldots $. \medskip 

\noindent \textbf{Proof.} This follows immediately from lemma 4.2.  \hfill $\square $ \medskip 

\noindent \textbf{Corollary 4.2b} The union of $K^{\ast}, K^{\ast}_{k_0}, K^{\ast}_{k_0k_1}, \ldots K^{\ast }_{k_0k_1 \cdots k_{\ell }}, \ldots $, where $\ell \ge 0$, $0 \le j \le \ell $, and 
$0 \le k_j \le 2n-1$,  covers 
$\C \setminus {\mathbb{V}}^{+}$, that is, 
\begin{displaymath}
K^{\ast } \cup \bigcup_{\ell \ge 0} \, \bigcup_{0 \le j \le \ell } \, \, \, 
\bigcup_{0 \le k_j\le 2n -1} K^{\ast }_{k_0 k_1 \cdots k_{\ell }} = 
\C \setminus {\mathbb{V}}^{+}. 
\end{displaymath}  

\noindent \textbf{Proof.} This follows immediately from 
$K^{\ast }_{k_0k_1 \cdots k_{\ell }} = 
{\tau }_{k_{\ell }}\comp \cdots \comp {\tau }_{k_0}(K^{\ast })$.  \hfill $\square $ \medskip 

Let $\mathcal{T}$ be the abelian subgroup of the $2$-dimensional Eulcidean group 
$\mathrm{E}(2)$ generated by the translations ${\tau }_j$ (\ref{eq-s4three}) for 
$j =0, 1, \ldots 2n-1$. It follows from corollary 4.2b that the regular stellated 
$n$-gon $K^{\ast }$ with its vertices and center removed is the fundamental domain for the action of the abelian group $\mathcal{T}$ on $\C \setminus {\mathbb{V}}^{+}$. 
The group $\mathcal{T}$ is isomorphic to the abelian subgroup 
$\mathfrak{T}$ of $(\C , +) $ generated 
by ${\{ 2u_j \} }^{2n-1}_{j=0}$. \medskip

Next we define the group $\mathfrak{G}$ and show that it acts freely, properly, and transitively on 
$\C \setminus {\mathbb{V}}^{+}$. \medskip  

Consider the group $\mathfrak{G} =
G \ltimes \mathfrak{T} \subseteq G \times \mathfrak{T}$, which is 
the semidirect product of the dihedral group $G$, generated by the rotation 
$R$ through $2\pi /n$ and the reflection $U$ subject to the relations 
$R^n = e = U^2$ and $RU = UR^{-1}$, and the abelian group $\mathfrak{T}$. 
An element $(R^jU^{\ell }, 2u_k)$ of $\mathfrak{G}$ is the affine linear map 
\begin{displaymath}
(R^jU^{\ell }, 2u_k) : \C \rightarrow \C : z \mapsto R^j U^{\ell } z + 2u_k.
\end{displaymath}
Multiplication in $\mathfrak{G}$ is defined by 
\begin{equation}
(R^jU^{\ell }, 2u_k) \cdot (R^{j'}U^{{\ell}'}, 2u_{k'}) = \big( R^{j-j'}U^{\ell +{\ell }'} , 
(R^jU^{\ell }) (2u_{k'}) + 2u_k \big) ,
\label{eq-s4five}
\end{equation}
which is the composition of the affine linear map $(R^{j'}U^{{\ell }'}, 2u_{k'})$ followed by 
$(R^jU^{\ell }, 2u_k)$. The mappings $G \rightarrow \mathfrak{G}: R^j \mapsto (R^jU^{\ell },0)$ and $\mathfrak{T} \rightarrow \mathfrak{G}: 2u_k \mapsto (e, 2u_k)$ are injective, which allows us to identify the groups $G$ and $\mathfrak{T}$ with their image in $\mathfrak{G}$. Using (\ref{eq-s4five}) 
we may write an element $(R^jU^{\ell }, 2u_k)$ of $\mathfrak{G}$ as 
$(e, 2u_k) \cdot (R^jU^{\ell },0)$. So 
\begin{displaymath}
(e, 2u_{(j +2k) \bmod 2n}) \cdot (R^kU^{\ell },0) = (R^kU^{\ell }, 2u_{(j +2k)\bmod 2n}) , 
\end{displaymath} 
For every $z \in \C $ we have 
\begin{align}
R^kU^{\ell }z + 2u_{(j+2k)\bmod 2n} & = R^k U^{\ell }z + R^kU^{\ell }(2u_j), \quad 
\mbox{using (\ref{eq-s4two}),} \notag 
\end{align}
that is, 
\begin{displaymath}
(R^kU^{\ell } , 2u_{(j+2k)\bmod 2n}) = (R^kU^{\ell },R^kU^{\ell }(2u_j)) = 
(R^kU^{\ell },0) \cdot (e, 2u_j). 
\end{displaymath}
Hence 
\begin{equation}
(e, 2u_{(j+2k)\bmod 2n}) \cdot (R^kU^{\ell },0) = (R^kU^{\ell },0) \cdot (e, 2u_j), 
\label{eq-s4six}
\end{equation} 
which is just equation (\ref{eq-s4four}). The group $\mathfrak{G}$ acts on $\C $ as 
$\mathrm{E}(2)$ does, namely, by affine linear orthogonal mappings. 
Denote this action by  
\begin{displaymath}
\psi : \mathfrak{G} \times \C \rightarrow \C : ((g, \tau ), z ) \mapsto 
\tau (g(z)). 
\end{displaymath} 

\noindent \textbf{Lemma 4.3} The set of vertices 
${\mathbb{V}}^{+}$ (\ref{eq-s4threestardot}) is invariant under the $\mathfrak{G}$ action. \medskip 

\noindent \textbf{Proof.} Let $v \in {\mathbb{V}}^{+}$. Then for some 
$({\ell }'_0, \ldots , {\ell }'_{2n-1}) \in {\Z }^{2n}_{\ge 0}$ and some 
$w \in \mathbb{V} \cup \{ 0 \} $
\begin{displaymath}
v = {\tau }^{{\ell }'_0}_0 \comp \cdots \comp {\tau }^{{\ell }'_{2n-1}}_{2n-1}(w) = 
{\psi }_{(e, 2u' )}(w) ,
\end{displaymath}
where $u' = \sum^{2n-1}_{j =0} {\ell }'_j u_j$. For $(R^jU^{\ell }, 2u) \in \mathfrak{G}$ 
with $j =0,1, \ldots , n-1$ and $\ell =0,1$ we 
have  
\begin{align}
{\psi }_{(R^jU^{\ell }, 2u)}v & = {\psi }_{(R^jU^{\ell }, 2u) } \comp {\psi }_{(e, 2u')} (w) 
= {\psi }_{(R^jU^{\ell },2u) \cdot (e, 2u')}(w) \notag \\
& = {\psi }_{(R^jU^{\ell }, R^j U^{\ell }(2u') +2u)}(w) = 
{\psi }_{(e, 2(R^jU^{\ell }u' +u)) \cdot (R^jU^{\ell },0)}(w) \notag \\ 
& = {\psi }_{(e, 2(R^jU^{\ell }u'+u))} \big( {\psi }_{(R^jU^{\ell },0)} (w) \big)  = 
{\psi }_{(e, 2(R^jU^{\ell }u'+u))}(w'),  
\label{eq-s4seven}
\end{align}
where $w' = {\psi }_{(R^jU^{\ell }, 0)}(w) = R^jU^{\ell }(w) \in 
\mathbb{V} \cup \{ 0 \} $. If $\ell =0$, then  

\begin{align}
R^ju' & = R^j (\sum^{2n-1}_{k=0} {\ell }'_k u_k ) = 
\sum^{2n-1}_{k=0} {\ell }'_kR^j (u_k) = 
\sum^{2n-1}_{k=0} {\ell }'_k  u_{(k+ 2j) \bmod 2n };  \notag
\end{align}
while if $\ell =1$, then 
\begin{align}
R^jU(u') & = \sum^{2n-1}_{k=0} {\ell }'_k R^j(U(u_k)) = 
\sum^{2n-1}_{k=0} {\ell }'_kR^j (u_{k'(k)})  
= \sum^{2n-1}_{k=0} {\ell }'_k  u_{(k'(k)+2j)\bmod 2n } . \notag
\end{align}
Here $k'(k) =${\tiny $\left\{ \begin{array}{cl} \hspace{-5pt}2n-k +1, & 
\hspace{-8pt} \mbox{if $k$ is even} \\
\hspace{-5pt} 2n-k-1, & \hspace{-8pt} \mbox{if $k$ is odd}, \end{array} \right. $} \hspace{-5pt}see corollary 4.1a. So $(e, 2(R^jU^{\ell }u' +u)) \in \mathfrak{T}$, which implies ${\psi }_{(e, 2(R^jU^{\ell }u'+u))}(w') \in 
{\mathbb{V}}^{+}$, as desired. \hfill $\square $ \medskip 

\noindent \textbf{Lemma 4.4} The action of $\mathfrak{G}$ on 
$\C \setminus {\mathbb{V}}^{+}$ is free.  \medskip 

\noindent \textbf{Proof.} Suppose that for some $v \in \C \setminus {\mathbb{V}}^{+}$ and 
some $(R^jU^{\ell }, 2u) \in \mathfrak{G}$ we have $v = {\psi }_{(R^jU^{\ell },2u)}(v)$. Then $v$ lies in some $K^{\ast }_{k_0k_1\cdots k_{\ell }} $. So for some 
$v' \in K^{\ast}$ we have 
\begin{align}
v & = {\tau }^{{\ell }'_0}_0 \comp \cdots {\tau }^{{\ell }'_{2n-1}}_{2n-1}(v') 
= {\psi }_{(e, 2u')}(v'),  \notag 
\end{align} 
where $u' = \sum^{2n-1}_{j=0} {\ell }'_j u_j$ 
for some $({\ell }'_0, \ldots , {\ell }'_{2n-1}) \in ({\Z }_{\ge 0})^{2n}$. Thus 
\begin{displaymath}
{\psi }_{(e, 2u)}(v') = {\psi }_{(R^jU^{\ell }, 2u) \cdot (e,2u')}(v') = 
{\psi }_{(R^jU^{\ell }, 2R^jU^{\ell }u +2u)}(v'). 
\end{displaymath}
This implies $R^j U^{\ell }=e$, that is, $j=\ell =0$. So $2u =2R^ju' +2u = 2u' +2u$, that is, $u =0$. Hence $(R^jU^{\ell },u) = (e,0)$, which is the identity element of 
$\mathfrak{G}$. \hfill $\square $ \medskip 

\noindent \textbf{Lemma 4.5} The action of $\mathcal{T}$ (and hence 
$\mathfrak{G}$) on $\C \setminus {\mathbb{V}}^{+}$ is transitive. \medskip 

\noindent \textbf{Proof.} Let $K^{\ast }_{k_0 \cdots k_{\ell }}$ and 
$K^{\ast }_{k'_0 \cdots k'_{{\ell }'}}$ lie in 
\begin{displaymath}
\C \setminus {\mathbb{V}}^{+} = 
K^{\ast } \cup \bigcup_{\ell \ge 0} \, \bigcup_{0 \le j \le \ell } \, \, 
\bigcup_{0 \le k_j\le 2n -1} K^{\ast }_{k_0 k_1 \cdots k_{\ell }}. 
\end{displaymath}
Since $K^{\ast }_{k_0 \cdots k_{\ell }} = 
{\tau }_{k_{\ell }} \comp \cdots \comp {\tau }_{k_0}(K^{\ast })$ and 
$K^{\ast }_{k'_0 \cdots k'_{{\ell }'}} = 
{\tau }_{k'_{{\ell }'}} \comp \cdots \comp {\tau }_{k'_0}(K^{\ast })$, it follows that  
$({\tau }_{k'_{{\ell }'}} \comp \cdots \comp {\tau }_{k'_0}) \comp 
({\tau }_{k_{\ell }} \comp \cdots \comp {\tau }_{k_0})^{-1}(K^{\ast }_{k_0 \cdots k_{\ell }}) =
K^{\ast }_{k'_0 \cdots k'_{{\ell }'}}$. \hfill $\square $ \medskip 

The action of $\mathfrak{G}$ on $\C \setminus {\mathbb{V}}^{+}$ is proper because 
$\mathfrak{G}$ is a discrete subgroup of $\mathrm{E}(2)$ with no accumulation points.  \medskip 

We now define an edge of $\C \setminus {\mathbb{V}}^{+}$ and what it means for an unordered 
pair of edges to be equivalent. We show that the group $\mathfrak{G}$ acts freely and properly 
on the identification space of equivalent edges. \medskip 

Let $E$ be an open edge of $K^{\ast }$. Since  
$E_{k_0 \cdots k_{\ell }} = {\tau }_{k_0} \cdots {\tau }_{k_{\ell }}(E) \in 
K^{\ast }_{k_0 \cdots k_{\ell }}$, it follows that $E_{k_0 \cdots k_{\ell }}$ is 
an open edge of $K^{\ast }_{k_0 \cdots k_{\ell }}$. Let 
\begin{displaymath}
\mathfrak{E} = \{ E_{k_0\cdots k_{\ell }} \setrule 
\ell \ge 0, \, \, 0 \le j \le \ell  \, \, \& \, \, 0 \le k_j \le 2n-1 \} . 
\end{displaymath}
Then $\mathfrak{E}$  is the set of open edges of 
$\C \setminus {\mathbb{V}}^{+}$ by lemma 4.2b. Since ${\tau }_{k_{\ell }} \comp \cdots 
\comp {\tau }_{k_0}(0)$ is the center of $K^{\ast }_{k_0 \cdots k_{\ell}}$, the element  
$(e, {\tau }_{k_{\ell }} \comp \cdots \comp {\tau }_{k_0}) \cdot 
(g, ({\tau }_{k_{\ell }} \comp \cdots \comp {\tau }_{k_0})^{-1})$ of $\mathfrak{G}$ is a 
rotation-reflection of $K^{\ast }_{k_0 \cdots k_{\ell}}$, which sends an edge of 
$K^{\ast }_{k_0 \cdots k_{\ell }}$ to another edge of $g \ast K^{\ast }_{k_0 \cdots k_{\ell}}$. Thus $\mathfrak{G}$ sends $\mathfrak{E}$ into itself. 
For $j =0,1, \infty $ let ${\mathfrak{E}}^j_{k_0 \cdots k_{\ell }}$ 
be the set of unordered pairs $[E_{k_0\cdots k_{\ell }}, E'_{k_0\cdots k_{\ell }}]$ of 
equivalent open edges of $K^{\ast }_{k_0 \cdots k_{\ell }}$, that is, 
$E_{k_0\cdots k_{\ell }} \cap E'_{k_0\cdots k_{\ell }} = \varnothing $, so the open edges 
$E_{k_0 \cdots k_{\ell }} = {\tau }_{k_0} \cdots {\tau }_{k_{\ell }}(E)$ and 
$E'_{k_0 \cdots k_{\ell }} = {\tau }_{k_0} \cdots {\tau }_{k_{\ell }}(E')$ of 
$\mathrm{cl}(K^{\ast }_{k_0 \cdots k_{\ell }})$ are not adjacent, which implies that the 
open edges $E$ and $E'$ of $K^{\ast }$ are not adjacent, and for some generator 
$S^{(j)}_m$ of the group $G^j$ of reflections we have 
\begin{displaymath}
E'_{k_0 \cdots k_{\ell }} = 
({\tau }_{k_0} \comp \cdots \comp {\tau }_{k_0})\big( S^{(j)}_m(({\tau }_{k_{\ell}} \comp \cdots \comp {\tau }_{k_0})^{-1}(E_{k_0 \cdots k_{\ell }} )) \big) . 
\end{displaymath}
Let ${\mathfrak{E}}^j = \cup_{\ell \ge 0} \cup_{0 \le j \le \ell } \, \cup_{0 \le k_j \le 2n-1} 
{\mathfrak{E}}^j_{k_0 \cdots k_{\ell }}$. Then ${\mathfrak{E}}^j$ is the 
set of unordered pairs of equivalent edges of $\C \setminus {\mathbb{V}}^{+}$. Define an action 
$\ast $ of $\mathfrak{G}$ on ${\mathcal{E}}^j$ by 
\begin{align}
(g, \tau ) \ast [ E_{k_0 \cdots k_{\ell }}, E'_{k_0 \cdots k_{\ell }} ] & =  
\big( [({\tau }' \comp \tau )( g({\tau }')^{-1}(E_{k_0 \cdots k_{\ell }})), 
({\tau }' \comp \tau )(g (({\tau }')^{-1}(E'_{k_0 \cdots k_{\ell }})) ] \big)  \notag \\
&= [(g,\tau) \ast E_{k_0\cdots k_{\ell }}, (g,\tau ) \ast E'_{k_0\cdots k_{\ell }}], \notag
\end{align}
where ${\tau }' = {\tau }_{k_{\ell }} \comp \cdots {\tau }_{k_0}$.  \medskip 

Define a relation $\sim $ on $\C \setminus {\mathbb{V}}^{+}$ as follows. We say that 
$x$ and $y \in \C \setminus {\mathbb{V}}^{+}$ are related, $x \sim y$, if 
1) $x \in F = \tau (E) \in {\mathfrak{E}}^0$ and $y \in F' = \tau (E') \in {\mathfrak{E}}^0$ 
such that $[F, F'] = [\tau (E), \tau (E')] \in {\mathfrak{E}}^0$, where $[E, E'] \in 
{\mathcal{E}}^0$ with $E' = S^{(0)}_m (E)$ for some $S^{(0)}_m \in G^0$ 
and $y = \tau \big( S^{(0)}_m({\tau }^{-1}(x)) \big)$ or 2) $x$, $y \in 
\big( \C \setminus {\mathbb{V}}^{+} \big) \setminus \mathfrak{E}$ and $x =y$. Then 
$\sim $ is an equivalence relation on $\C \setminus {\mathbb{V}}^{+}$. Let 
$(\C \setminus {\mathbb{V}}^{+})^{\sim }$ be the set of equivalence classes and 
let $\Pi $ be the map 
\begin{equation}
\Pi : \C \setminus {\mathbb{V}}^{+} \rightarrow (\C \setminus {\mathbb{V}}^{+})^{\sim }: 
p \mapsto [p], 
\label{eq-s4sevenstar}
\end{equation}
which assigns to every $p \in \C \setminus {\mathbb{V}}^{+}$ the 
equivalence class $[p]$ containing $p$. \medskip 

\noindent \textbf{Lemma 4.6} $\Pi |_{K^{\ast }}$ is the map $\pi $ (\ref{eq-s3threestar}). \medskip 

\noindent \textbf{Proof.} This follows immediately from the definition of the maps 
$\Pi $ and $\pi $. \hfill $\square $ \medskip 

\noindent \textbf{Lemma 4.7} The usual action of $\mathfrak{G}$ on $\C $, restricted 
to $\C \setminus {\mathbb{V}}^{+}$, is compatible with the equivalence relation 
$\sim $, that is, if $x$, $y \in \C \setminus \mathbb{V}$ and $x \sim y$, then 
$(g, \tau )(x) \sim (g, \tau )(y)$ for every $(g, \tau ) \in \mathfrak{G}$. \medskip 

\noindent \textbf{Proof.} Suppose that $x \in F = {\tau }'(E)$, where ${\tau }' \in 
\mathcal{T}$. Then $y \in F' = {\tau }'(E')$, since $x \sim y$. So for some 
$S^{(0)}_m \in G^{\, 0}$ we have $({\tau }')^{-1}(y) = 
S^{(0)}_m({\tau }^{-1}(x))$. Let $(g, \tau ) \in \mathfrak{G}$. Then 
\begin{displaymath}
(g, \tau )\big( ({\tau }')^{-1}(y) \big) = 
g( ({\tau }')^{-1}(y) ) + u_{\tau } = g\big( S^{(0)}_m({\tau }^{-1}(x)) \big) + u_{\tau }. 
\end{displaymath}
So $(g, \tau )(y) \in (g, \tau )\ast F'$. But  $(g, \tau )(x) \in (g, \tau ) \ast F$ and 
$[ (g, \tau ) \ast F , (g, \tau ) \ast F' ]  = 
(g , \tau )\ast [F,F']$. Hence $(g, \tau )(x) \sim (g, \tau )(y)$. \hfill $\square $ \medskip 

Because of lemma 4.7, the usual $\mathfrak{G}$-action on 
$\C \setminus {\mathbb{V}}^{+}$ induces an action of $\mathfrak{G}$ on 
$(\C \setminus {\mathbb{V}}^{+})^{\sim }$. \medskip 

\noindent \textbf{Lemma 4.8} The action of $\mathfrak{G}$ on 
$(\C \setminus {\mathbb{V}}^{+})^{\sim }$ is free and proper. \medskip 

\noindent \textbf{Proof.} The following argument shows that it is free. 
Using lemma A2 we see that an 
element of $\mathfrak{G}$, which lies in the isotropy group ${\mathfrak{G}}_{[F,F']}$ 
for $[F,F'] \in {\mathfrak{E}}^0$, interchanges the edge $F$ with the equivalent 
edge $F'$ and thus fixes the equivalence class $[p]$ for every $p \in F$. Hence 
the $\mathfrak{G}$ action on $(\C \setminus {\mathbb{V}}^{+})^{\sim }$ is free. It is 
proper because $\mathfrak{G}$ is a discrete subgroup of the Euclidean group 
$\mathrm{E}(2)$ with no accumulation points. \hfill $\square $ \medskip  

\noindent \textbf{Theorem 4.9} The $\mathfrak{G}$-orbit space 
$(\C \setminus {\mathbb{V}}^{+})^{\sim }/ \mathfrak{G}$ is holomorphically 
diffeomorphic to the $G$-orbit space $(K^{\ast } \setminus \{ O \})^{\sim }/G = 
{\widetilde{S}}_{\mathrm{reg}}$. \medskip 

\noindent \textbf{Proof.} This claim follows from the fact that the fundamental 
domain of the $\mathfrak{G}$-action on $ \C \setminus {\mathbb{V}}^{+}$ is 
$K^{\ast } \setminus \{ O \} $, which is the fundamental domain of the 
$G$-action on $K^{\ast } \setminus \{ O \}$. Thus $\Pi ( \C \setminus {\mathbb{V}}^{+})$ 
is a fundamental domain of the $\mathfrak{G}$-action on $(\C \setminus 
{\mathbb{V}}^{+})^{\sim }$, which is equal to $\pi (K^{\ast } \setminus \{ O \} ) 
= (K^{\ast } \setminus \{ O \})^{\sim }$ by lemma 4.6. Hence the $\mathfrak{G}$-orbit space $(\C \setminus {\mathbb{V}}^{+})^{\sim }/ 
\mathfrak{G}$ is equal to the $G$-orbit space ${\widetilde{S}}_{\mathrm{reg}}$. 
So the identity map from $\Pi (\C \setminus {\mathbb{V}}^{+})$ to 
$(K^{\ast }\setminus \{ O \})^{\sim }$ induces a holomorphic diffeomorphism 
of orbit spaces. \hfill $\square $ \medskip 

Because the group $\mathfrak{G}$ is a discrete subgroup of the $2$-dimensional 
Euclidean group $\mathrm{E}(2)$, the Riemann surface 
$(\C \setminus {\mathbb{V}}^{+})^{\sim} /\mathfrak{G}$ is an \emph{affine} model of the 
affine Riemann surface ${\mathcal{S}}_{\mathrm{reg}}$. 

\section{The developing map and geodesics}

In this section we show that the mapping 
\begin{equation}
\delta : \mathcal{D} \subseteq {\mathcal{S}}_{\mathrm{reg}} \rightarrow Q \subseteq \C :
(\xi , \eta ) \rightarrow F_Q \big( \widehat{\pi }(\xi ,\eta ) \big) 
\label{eq-s4onenw}
\end{equation}
straightens the holomorphic vector field $X$ (\ref{eq-s2eight}) on the fundamental domain $\mathcal{D} \subseteq {\mathcal{S}}_{\mathrm{reg}}$, see Bates and Cushman \cite{bates-cushman} and Flaschka \cite{flaschka}. We verify that $X$ is the geodesic vector field for a flat Riemannian metric $\Gamma $ on $\mathcal{D}$. \medskip 

First we rewrite equation (\ref{eq-s2eightdot}) as 
\begin{equation}
T_{(\xi , \eta )} \widehat{\pi } \big( X(\xi ,\eta ) \big) = \eta \frac{\partial }{\partial \xi }, 
\quad \mbox{for $(\xi ,\eta ) \in \mathcal{D}$.}
\label{eq-s4twonw}
\end{equation}
From the definition of the mapping $F_Q$ (\ref{eq-s1two}) we get  
\begin{displaymath}
\dee z = \dee F_Q = 
\frac{1}{ \big( {\xi }^{n-n_0}(1-\xi )^{n-n_1}\big)^{\raisebox{-2pt}{$\scriptstyle 1/n$}}} \dee  \xi = \frac{1}{\eta } \dee \xi  ,
\end{displaymath}
where we use the same complex $n^{\mathrm{th}}$ root as in the definition of $F_Q$. This implies 
\begin{equation}
\frac{\partial }{\partial z} = T_{\xi }F_Q \big( \eta \frac{\partial }{\partial \xi } \big) ,  
\quad \mbox{for $(\xi , \eta ) \in \mathcal{D}$}
\label{eq-s4threenw}
\end{equation}
For each $(\xi , \eta ) \in \mathcal{D}$ using (\ref{eq-s4twonw}) and (\ref{eq-s4threenw}) we get
\begin{align}
&\hspace{-25pt}T_{(\xi ,\eta )}\delta \big( X(\xi ,\eta ) \big)  = 
\big( T_{\xi }F_Q \comp T_{(\xi ,\eta )}\widehat{\pi }\big) \big( X(\xi ,\eta ) \big) = 
T_{\xi }F_Q (\eta \frac{\partial }{\partial \xi }) = 
\frac{\partial }{\partial z} \rule[-10pt]{.5pt}{24pt}\, \raisebox{-9pt}{$\scriptscriptstyle z = \delta (\xi ,\eta )$} \hspace{-25pt} . \notag
\end{align}
So the holomorphic vector field $X$ (\ref{eq-s2eight}) on $\mathcal{D}$ and the holomorphic vector field $\frac{\partial }{\partial z}$ on $Q$ are $\delta $-related. Hence $\delta $ sends an integral curve of the vector field $X$ starting at $(\xi ,\eta ) \in 
\mathcal{D}$ onto an integral curve of the vector field $\frac{\partial }{\partial z}$ starting at 
$z =\delta (\xi ,\eta ) \in Q$. Since an integral curve of $\frac{\partial }{\partial z}$ is a horizontal line segment in $Q$, we have proved \medskip 

\noindent \textbf{Claim 5.1} The holomorphic mapping $\delta $ (\ref{eq-s4onenw}) 
straightens the holomorphic vector field $X$ (\ref{eq-s2eight}) on the fundamental domain 
$\mathcal{D} \subseteq {\mathcal{S}}_{\mathrm{reg}}$. \medskip

We can say more. Let $u = \mathrm{Re}\, z$ and $v = \mathrm{Im}\, z$. Then 
\begin{equation}
\gamma = \dee u \circdot \dee u + \dee v \circdot \dee v = 
\dee z \circdot \, \overline{\dee z}
\label{eq-s4fournw}
\end{equation}
is the flat Euclidean metric on $\C $. Its restriction 
$\gamma |_{\C \setminus {\mathbb{V}}^{+}}$ 
to $\C \setminus {\mathbb{V}}^{+}$ is invariant under the group $\mathfrak{G}$, 
which is a subgroup of the Euclidean group $\mathrm{E}(2)$.  \medskip 

Consider the flat Riemannian metric $\gamma |_Q$ on $Q$, where 
$\gamma $ is the metric (\ref{eq-s4fournw}) on $\C $. Pulling back $\gamma |_Q$ 
by the mapping $F_Q$ (\ref{eq-s1two}) gives a metric 
\begin{displaymath}
\widetilde{\gamma } = F^{\ast }_Q \gamma |_Q =  
{|{\xi }^{n-n_0}(1-\xi )^{n-n_1}|}^{-2/n} \dee \xi \circdot \, \overline{\dee \xi }
\end{displaymath} 
on $\C \setminus \{ 0, 1 \} $. Pulling the metric $\widetilde{\gamma }$ back by the 
projection mapping $\widetilde{\pi }:{\C }^2 \rightarrow \C :(\xi ,\eta ) \mapsto \xi $ gives 
\begin{displaymath}
\widetilde{\Gamma } = {\widetilde{\pi }}^{\ast }\widetilde{\gamma } = 
{|{\xi }^{n-n_0}(1-\xi )^{n-n_1}|}^{-2/n} \dee \xi \circdot \, \overline{\dee \xi } 
\end{displaymath}
on ${\C }^2 $. Restricting $\widetilde{\Gamma }$ to the affine Riemann surface ${\mathcal{S}}_{\mathrm{reg}}$ gives 
$\Gamma = \frac{1}{\eta } \dee \xi \, \circdot \, \frac{1}{\overline{\eta }} \overline{\dee \xi }$. \medskip 

\noindent \textbf{Lemma 5.2} $\Gamma $ is a flat Riemannian metric on 
${\mathcal{S}}_{\mathrm{reg}}$. \medskip

\noindent \textbf{Proof.} We compute. For every $(\xi , \eta ) \in {\mathcal{S}}_{\mathrm{reg}}$ we have 
\begin{align}
\mbox{\footnotesize $\Gamma (\xi ,\eta )\big( X(\xi ,\eta ), X(\xi ,\eta ) \big)$} & = \notag \\
&\hspace{-1.25in} = 
\mbox{\footnotesize $\frac{1}{\eta } \dee \xi \big( \eta \frac{\partial }{\partial \xi } + 
\ttfrac{n-n_0}{n} \frac{{\xi }{(1-\xi )} (1-\frac{2n-n_0-n_1}{n} \xi )}{{\eta }^{n-2}} \frac{\partial }{\partial \eta } \big) \cdot 
\frac{1}{\overline{\eta }} \overline{\dee \xi }\big( \overline{\eta} \overline{\frac{\partial }{\partial \xi }} + 
\ttfrac{n-n_0}{n} \frac{\overline{\xi (1-{\xi })(1-\frac{2n-n_0-n_1}{n}{\xi })}}{{\overline{\eta }}^{n-2}}
\, \overline{\frac{\partial }{\partial \eta }} \big) $} \notag \\
&\hspace{-1.25in} = \mbox{\footnotesize $\frac{1}{\eta } \dee \xi \big( \eta \frac{\partial }{\partial \xi } \big) \cdot  \frac{1}{\overline{\eta }} \overline{\dee \xi } \big( \overline{\eta } \overline{\frac{\partial }{\partial \xi }} \big) =1 $.} \notag
\end{align}
Thus $\Gamma $ is a Riemannian metric on ${\mathcal{S}}_{\mathrm{reg}}$. It is flat by 
construction. \hfill $\square $ \medskip 

Because $\mathcal{D}$ has nonempty interior and the map $\delta $ (\ref{eq-s4onenw}) is holomorphic, 
it can be analytically continued to the map  
\begin{align}
& {\delta }_Q: {\mathcal{S}}_{\mathrm{reg}} \subseteq {\C }^2 \rightarrow Q \subseteq \C : 
(\xi ,\eta ) \mapsto F_Q \big( \widehat{\pi }(\xi ,\eta ) \big), 
\label{eq-s4fivenw}
\end{align}
since $\delta = {\delta }_Q|_{\mathcal{D}}$. By construction 
${\delta }^{\ast }_Q (\gamma |_Q) = \Gamma $. 
So the mapping ${\delta }_Q$ 
is an isometry of $({\mathcal{S}}_{\mathrm{reg}}, \Gamma )$ onto $(Q, \gamma |_Q )$. 
In particular, the map $\delta $ is an isometry of 
$(\mathcal{D}, \Gamma |_{\mathcal{D}})$ onto $(Q, \gamma |_Q)$. Moreover, $\delta $ is 
a local holomorphic diffeomorphism, because for every $(\xi ,\eta ) \in \mathcal{D}$, the complex linear mapping $T_{(\xi ,\eta )}\delta $ is an isomorphism, since it sends 
$X(\xi ,\eta )$ to 
$\frac{\partial }{\partial z}\rule[-6pt]{.5pt}{15pt}\, \raisebox{-6pt}{$\scriptscriptstyle z 
= \delta (\xi ,\eta )$}$  \hspace{-25pt}. Thus $\delta $ is a \emph{developing map} in the sense of differential geometry, see Spivak \cite[p.97]{spivak} note on \S 12 of Gauss \cite{gauss}. The map $\delta $ is \emph{local} because the integral curves of $\frac{\partial }{\partial z}$ on $Q$ are only defined for a finite time, since they are horizontal line segments in $Q$. Thus the integral curves of $X$ (\ref{eq-s2eight}) on $\mathcal{D}$ are defined for a finite time. Since the integral curves of $\frac{\partial }{\partial z}$ are geodesics on $(Q, \gamma |_Q)$, the image of a local integral curve of 
$\frac{\partial }{\partial z}$ under the local inverse of the mapping $\delta $ is a local integral curve of $X$. This latter local integral curve is a geodesic on 
$(\mathcal{D}, \Gamma |_{\mathcal{D}})$, since $\delta$ is an isometry. Thus we have proved \medskip 

\noindent \textbf{Claim 5.3} The holomorphic vector field $X$ (\ref{eq-s2eight}) on 
the fundamental domain $\mathcal{D}$ is the geodesic vector field for the flat Riemannian metric $\Gamma |_{\mathcal{D}}$ on $\mathcal{D}$. \medskip 

\noindent \textbf{Corollary 5.3a} The holomorphic vector field $X$ on the affine Riemann 
surface ${\mathcal{S}}_{\mathrm{reg}}$ is the geodesic vector field for the flat Riemannian 
metric $\Gamma $ on ${\mathcal{S}}_{\mathrm{reg}}$. \medskip 

\noindent \textbf{Proof.} The corollary follows by analytic continuation from the conclusion of claim 5.3, since $\mathrm{int}\, \mathcal{D}$ is a nonempty open subset of 
${\mathcal{S}}_{\mathrm{reg}}$ and both the vector field $X$ and the Riemannian metric 
$\Gamma $ are holomorphic on ${\mathcal{S}}_{\mathrm{reg}}$. 

\section{Discrete symmetries and billiard motions}

Let $\mathcal{G}$ be the group of homeomorphisms of the affine Riemann surface 
$\mathcal{S}$ (\ref{eq-s2one}) generated by the mappings 
\begin{displaymath}
\mathcal{R}: \mathcal{S} \rightarrow \mathcal{S}: (\xi ,\eta ) \mapsto (\xi , {\mathrm{e}}^{2\pi i/n}\eta )\, \, \,  \mathrm{and}\, \, \,  \mathcal{U}:\mathcal{S} \rightarrow \mathcal{S}: (\xi , \eta ) \mapsto (\overline{\xi }, \overline{\eta }). 
\end{displaymath}
Clearly, the relations ${\mathcal{R}}^n = {\mathcal{U}}^2 = e$ hold. For every $(\xi, \eta ) \in \mathcal{S}$ we have  
\begin{align}
\mathcal{U}{\mathcal{R}}^{-1}(\xi , \eta )  & = \mathcal{U}(\xi ,{\mathrm{e}}^{-2\pi i/n}\eta ) 
= (\overline{\xi }, {\mathrm{e}}^{2\pi i /n}\overline{\eta }) 
= \mathcal{R}(\overline{\xi }, \overline{\eta }) = \mathcal{R} \, \mathcal{U}(\xi , \eta ) . \notag 
\end{align} 
So the additional relation $\mathcal{U}{\mathcal{R}}^{-1} = \mathcal{R}\, \mathcal{U}$ holds. Thus $\mathcal{G}$  is isomorphic to the dihedral group $D_{2n}$. \medskip 

\noindent \textbf{Lemma 6.1} $\mathcal{G}$ is a group of isometries of 
$({\mathcal{S}}_{\mathrm{reg}}, \Gamma )$. \medskip 

\noindent \textbf{Proof.} For every $(\xi , \eta ) \in {\mathcal{S}}_{\mathrm{reg}}$ we get 
\begin{align}
{\mathcal{R}}^{\ast }\Gamma (\xi ,\eta )\big( X(\xi ,\eta ), X(\xi ,\eta ) \big) & = 
\Gamma \big( \mathcal{R}(\xi , \eta ) \big) \big( T_{(\xi ,\eta )}\mathcal{R}\big( X(\xi ,\eta ) \big) , 
T_{(\xi ,\eta )}\mathcal{R}\big( X(\xi ,\eta ) \big) \big) \notag \\
&\hspace{-1.25in} = \Gamma (\xi , {\mathrm{e}}^{2\pi i/n}\eta ) \big( 
{\mathrm{e}}^{2\pi i/n} \eta \frac{\partial }{\partial \xi } + \ttfrac{n-n_0}{n} 
\frac{\xi (1-\xi )(1-\ttfrac{2n-n_0-n_1}{n-n_0}\xi )}{{\eta }^{n-2}} 
{\mathrm{e}}^{2\pi i/n}\frac{\partial }{\partial \eta }, 
\notag \\
&\hspace{-.25in} {\mathrm{e}}^{2\pi i/n} \eta \frac{\partial }{\partial \xi } + 
\frac{n-n_0}{n} \ttfrac{\xi (1-\xi )(1-\tfrac{2n-n_0-n_1}{n-n_0} \xi )}{{\eta }^{n-2}} 
{\mathrm{e}}^{2\pi i/n}\frac{\partial }{\partial \eta } \big) \notag \\
& \hspace{-1.25in} = \frac{1}{|{\mathrm{e}}^{2\pi i/n}\eta |^2} \dee \xi \big( {\mathrm{e}}^{2\pi i/n} \eta \frac{\partial }{\partial \xi } \big)  \cdot \, 
\overline{\dee \xi } \big( \overline{{\mathrm{e}}^{2\pi i /n} \eta \frac{\partial }{\partial \xi } } \big) = 1 \notag \\
&\hspace{-1.25in} = \frac{1}{|\eta |^2} \dee \xi (\eta \frac{\partial }{\partial \xi } \big)  \cdot 
\overline{\dee \xi }(\overline{\eta  \frac{\partial }{\partial \xi }}) = 
\Gamma (\xi , \eta )\big( X(\xi ,\eta ), X(\xi ,\eta ) \big)  \notag 
\end{align}
and 
\begin{align}
{\mathcal{U}}^{\ast } \Gamma (\xi ,\eta ) \big( X(\xi , \eta ), X(\xi , \eta ) \big) 
& = \Gamma \big(\mathcal{U}(\xi ,\eta ) \big) \big( T_{(\xi ,\eta )}\mathcal{U} \big(X(\xi ,\eta ) \big) , 
T_{(\xi ,\eta )}\mathcal{U} \big( X(\xi ,\eta ) \big) \big) \notag \\
& \hspace{-1.25in} = 
\frac{1}{|\eta |^2} \overline{\dee \xi }(\overline{\eta \frac{\partial }{\partial \xi }} ) 
\cdot \overline{\overline{\dee \xi }}(\overline{\overline{\eta \frac{\partial }{\partial \xi }}} ) = 
\Gamma (\xi ,\eta ) \big( X(\xi , \eta ), X(\xi ,\eta ) \big) . \tag*{$\square $} 
\end{align}

Recall that the group $G$, generated by the linear mappings 
\begin{displaymath}
R: \C \rightarrow \C : z \mapsto {\mathrm{e}}^{2\pi i /n}z \, \, \, \mathrm{and} \, \, \,  
U: \C \rightarrow \C: z \mapsto \overline{z}, 
\end{displaymath}
is isomorphic to the dihedral group $D_{2n}$. \medskip 

\noindent \textbf{Lemma 6.2} $G$ is a group of isometries of $(\C , \gamma )$. \medskip 

\noindent \textbf{Proof.} This follows because $R$ and $U$ are Euclidean motions. \hfill $\square $ \medskip  

We would like the developing map ${\delta }_Q$ (\ref{eq-s4fivenw}) to intertwine the actions of $\mathcal{G}$ and $G$ and the geodesic flows on 
$({\mathcal{S}}_{\mathrm{reg}}, \Gamma )$ and $(Q, \gamma |_Q)$. There are several difficulties. The first is: the group $G$ does \emph{not} preserve the quadrilateral $Q$. To overcome this difficulty we extend the mapping ${\delta }_Q$ (\ref{eq-s4fivenw}) to the mapping ${\delta }_{K^{\ast }}$ (\ref{eq-s3three}) of the affine Riemann surface ${\mathcal{S}}_{\mathrm{reg}}$ onto the regular stellated $n$-gon $K^{\ast }$. \medskip 

\noindent \textbf{Lemma 6.3} The mapping ${\delta }_{K^{\ast }}$ (\ref{eq-s3three}) intertwines the action $\Phi $ (\ref{eq-s2ten}) of $\mathcal{G}$ on 
${\mathcal{S}}_{\mathrm{reg}}$ with the action
\begin{equation}
\Psi : G \times K^{\ast } \rightarrow K^{\ast }: (g,z) \mapsto g(z) 
\label{eq-s5one}
\end{equation}
of $G$ on the regular stellated $n$-gon $K^{\ast }$. \medskip 

\noindent \textbf{Proof.} From the definition of the mapping ${\delta }_{K^{\ast }}$ we see that for each $(\xi ,\eta ) \in \mathcal{D}$ we have 
${\delta }_{K^{\ast }} \big( {\mathcal{R}}^j(\xi ,\eta ) \big) = 
R^j {\delta }_{K^{\ast }}(\xi ,\eta )$ for every $j \in \Z $. By analytic 
continuation we see that the preceding 
equation holds for every $(\xi ,\eta ) \in {\mathcal{S}}_{\mathrm{reg}}$. 
Since $F_Q(\overline{\xi }) = \overline{F_Q(\xi )}$ by construction and 
$\widehat{\pi }( \overline{\xi }, \overline{\eta } ) = \overline{\xi }$ (\ref{eq-s2seven}), from the definition of the mapping $\delta $ (\ref{eq-s4onenw}) we get 
$\delta (\overline{\xi }, \overline{\eta } ) = \overline{\delta (\xi ,\eta )}$ for every 
$(\xi ,\eta ) \in \mathcal{D}$. In other words, 
${\delta }_{K^{\ast }}\big( \mathcal{U}(\xi ,\eta ) \big) = 
U{\delta }_{K^{\ast }}(\xi, \eta )$ for every $(\xi ,\eta ) \in \mathcal{D}$. 
By analytic continuation we see 
that the preceding equation holds for all $(\xi ,\eta ) \in {\mathcal{S}}_{\mathrm{reg}}$. Hence on ${\mathcal{S}}_{\mathrm{reg}}$ we have 
\begin{equation}
{\delta }_{K^{\ast }} \comp {\Phi }_g = 
{\Psi }_{\varphi (g)} \comp {\delta }_{K^{\ast }} \quad \mbox{for every $g \in \mathcal{G}$.}
\label{eq-s5twostar}
\end{equation} 
The mapping $\varphi : \mathcal{G} \rightarrow G$ sends the generators 
$\mathcal{R}$ and $\mathcal{U}$ of the group $\mathcal{G}$ to the generators $R$ and $U$ of the group $G$, respectively. So it is an isomorphism. \hfill $\square $ \medskip 

There is a second more serious difficulty: the integral curves of 
$\frac{\partial }{\partial z}$ run off the quadrilateral $Q$ in finite time. 
We fix this by requiring that when an integral curve reaches a point $P$ on the boundary 
$\partial Q $ of $Q$, which is not a vertex, it undergoes a specular reflection at $P$. (If the integral curve reaches a vertex of $Q$ in forward or backward time, then the motion ends). This motion can be continued as a straight line motion, which extends the motion on the original segment in $Q$ or $S(Q)$. To make this precise, we give 
$Q$ the orientation induced from $\C $ and suppose that the incoming (and hence 
outgoing) straight line motion has the same orientation as $\partial Q$. If the 
incoming motion makes an angle $\alpha $ with respect to the inward pointing 
normal $N$ to $\partial Q$ at $P$, then the outgoing motion makes an angle 
$\alpha $ with the normal $N$, see Richens and Berry \cite{richens-berry}. Specifically, if the incoming motion to $P$ is an integral curve of 
$\frac{\partial }{\partial z}$, then the outgoing motion, after reflection at $P$, is an integral curve of $R^{-1} \frac{\partial }{\partial z} = 
{\mathrm{e}}^{-2\pi i/n}\frac{\partial }{\partial z}$. Thus the outward motion makes a turn of $-2\pi /n$ at $P$ towards the interior of $Q$, see figure 10 (left). In figure 10 (right) 
the incoming motion has the opposite orientation from $\partial Q$.  
\vspace{-.1in}\par \noindent \hspace{.5in} \begin{tabular}{l}
\vspace{-.15in} \\
\includegraphics[width=250pt]{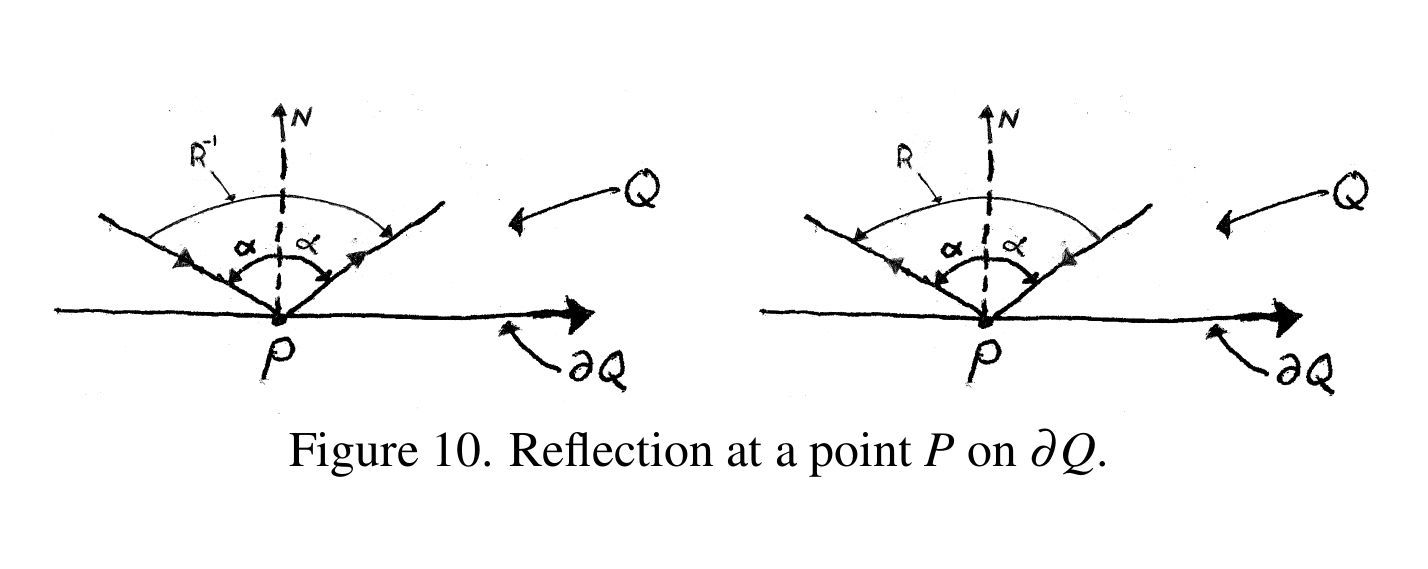}
\end{tabular}
\vspace{-.1in} \par \noindent This extended motion on $Q$ is called a billiard motion. A billiard motion starting in the interior of $\mathrm{cl}(Q)$ is defined for \emph{all} time and remains in 
$\mathrm{cl}(Q)$ less its vertices, since each of the segments of the billiard motion is 
a straight line parallel to an edge of $\mathrm{cl}(Q)$ and does not hit a vertex of 
$\mathrm{cl}(Q)$, see figure 12.  \medskip 

We can do more. If we apply a reflection $S$ in the edge of $Q$ in its boundary 
$\partial Q$, which contains the reflection point $P$, to the initial reflected motion at 
$P$, and then again to the extended straight line motion in $S(Q)$ when it reaches %
\linebreak  

\vspace{.1in} \par \noindent \hspace{1in}\begin{tabular}{l}
\vspace{-.3in} \\
\includegraphics[width=190pt]{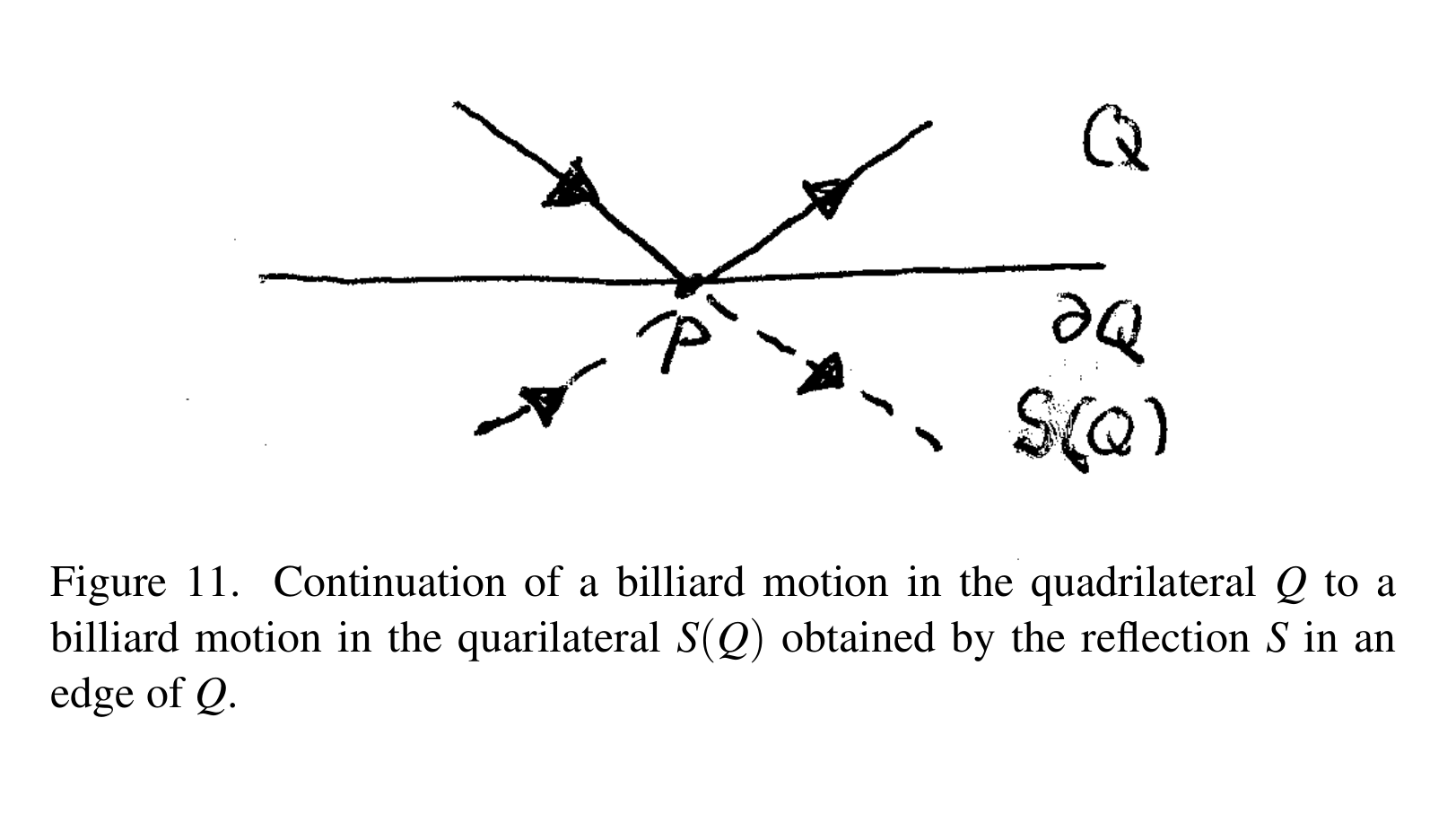}
\end{tabular}
\par \noindent $\partial S(Q)$, et cetera, we see that the extended motion becomes a billiard motion in the 
regular stellated $n$-gon $K^{\ast }= Q \cup {\amalg}_{0\le k \le n-1}SR^k(Q) \big)$, see figure 12. 

\par\noindent \hspace{.6in}\begin{tabular}{l}
\vspace{-.25in}\\
\includegraphics[width=250pt]{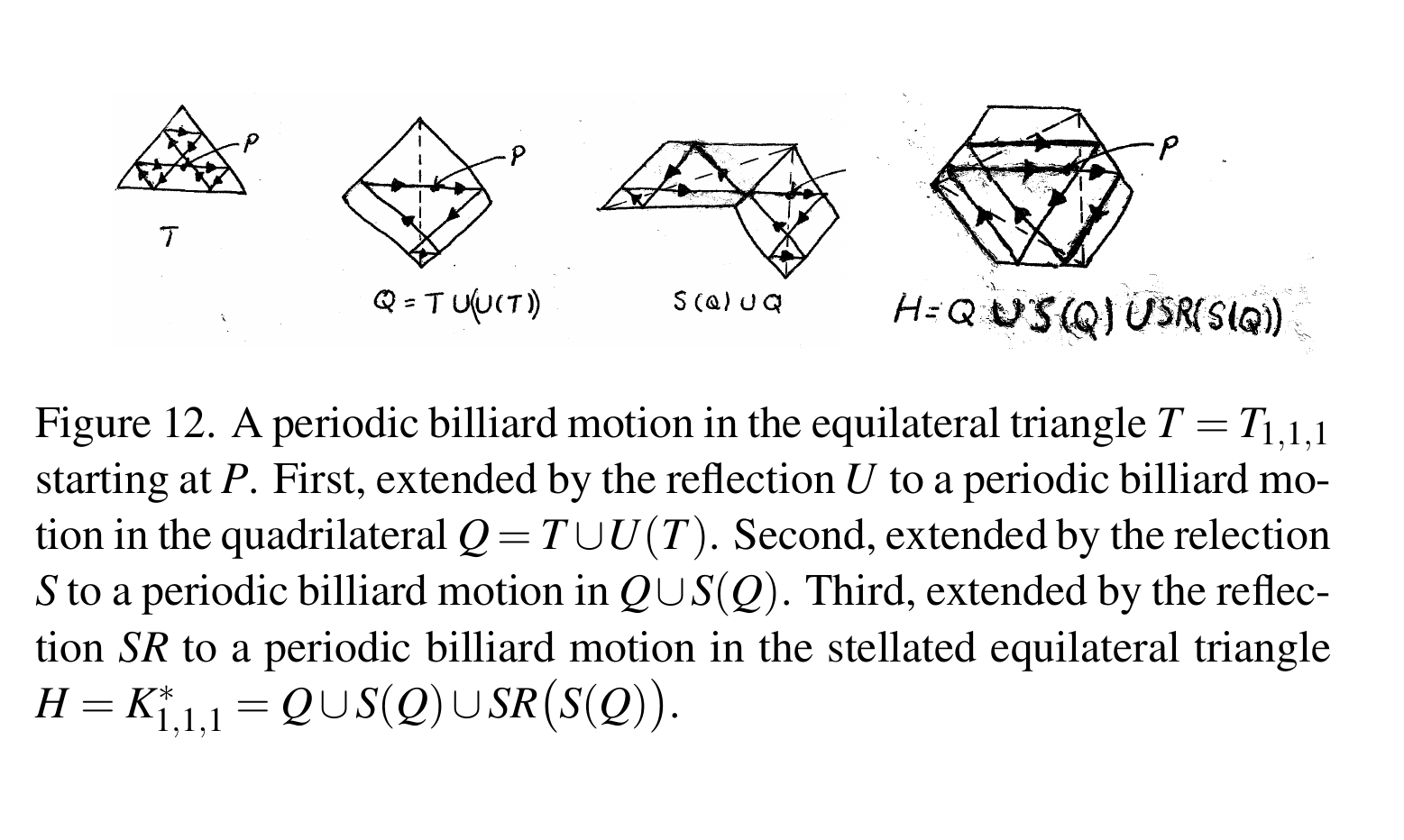}
\end{tabular} 

\vspace{-.15in} \par \noindent So we have verified \medskip 

\noindent \textbf{Claim 6.4} A billiard motion in the regular stellated $n$-gon $K^{\ast }$, which starts at a point in the interior of $K^{\ast } \setminus \{ O \} $ does not hit a vertex of $\mathrm{cl}(K^{\ast })$ and is invariant under the action of the isometry subgroup $\widehat{G}$ of the isometry group $G$ of 
$(K^{\ast }, \gamma |_{K^{\ast }})$ generated by the rotation $R$. \medskip 

Let $\widehat{\mathcal{G}}$ be the subgroup of $\mathcal{G}$ generated by the rotation $\mathcal{R}$. We now show \medskip 

\noindent \textbf{Lemma 6.5} The holomorphic vector field $X$ (\ref{eq-s2eight}) 
on ${\mathcal{S}}_{\mathrm{reg}}$ is $\widehat{\mathcal{G}}$-invariant. \medskip 

\noindent \textbf{Proof.} We compute. For every $(\xi ,\eta ) \in {\mathcal{S}}_{\mathrm{reg}}$ 
and for $\mathcal{R} \in \widehat{\mathcal{G}}$ we have 
\begin{align}
T_{(\xi ,\eta )}{\Phi }_{\mathcal{R}} \big( X(\xi ,\eta ) \big) & = 
{\mathrm{e}}^{2\pi i/n} \big[ \eta \frac{\partial }{\partial \xi } + 
\ttfrac{n-n_0}{n} \frac{\xi (1-\xi)(1-\ttfrac{2n-n_0-n_1}{n} \xi)}{ {\eta }^{n-2} }
\frac{\partial }{\partial \eta }  \big] \notag  \\
& \hspace{-.5in} = ({\mathrm{e}}^{2\pi i/n} \eta) \frac{\partial }{\partial \xi } + 
\ttfrac{n-n_0}{n} \frac{\xi (1-\xi )(1-\ttfrac{2n-n_0-n_1}{n} \xi )}
{ ({\mathrm{e}}^{2\pi i/n }\eta )^{n-2} } \frac{\partial }{\partial ({\mathrm{e}}^{2\pi i/n}\eta )}  
\notag \\
& \hspace{-.5in} = X(\xi , {\mathrm{e}}^{2\pi i/n} \eta ) = X \comp {\Phi }_{\mathcal{R}}(\xi ,\eta ). 
\notag 
\end{align}
Hence for every $j \in \Z $ we get 
\begin{equation}
T_{(\xi ,\eta )} {\Phi }_{{\mathcal{R}}^j } \big( X(\xi ,\eta ) \big) = 
X \comp {\Phi }_{{\mathcal{R}}^j }(\xi ,\eta )  
\label{eq-s5twoanw}
\end{equation}
for every $(\xi ,\eta ) \in {\mathcal{S}}_{\mathrm{reg}}$. In other words, the vector field 
$X$ is invariant under the action of $\widehat{\mathcal{G}}$ on ${\mathcal{S}}_{\mathrm{reg}}$. 
\hfill $\square $ \medskip 

\noindent \textbf{Corollary 6.5a} For every $(\xi ,\eta ) \in \mathcal{D}$ we have 
\begin{equation}
X|_{{\Phi }_{{\mathcal{R}}^j}(\mathcal{D}) } = T{\Phi }_{{\mathcal{R}}^j } \comp X|_{\mathcal{D}}. 
\label{eq-s5twobnw}
\end{equation}

\noindent \textbf{Proof.} Equation (\ref{eq-s5twobnw}) is a rewrite of equation 
(\ref{eq-s5twoanw}). \hfill $\square $ \medskip 

\noindent \textbf{Corollary 6.5b} Every geodesic on $({\mathcal{S}}_{\mathrm{reg}}, \Gamma )$ 
is $\widehat{\mathcal{G}}$-invariant. \medskip 

\noindent \textbf{Proof.} This follows immediately from the lemma. \hfill $\square $ \medskip 

\noindent \textbf{Lemma 6.6} For every $(\xi ,\eta ) \in {\mathcal{S}}_{\mathrm{reg}}$ 
and every $j \in \Z$ we have 
\begin{equation}
T_{{\Phi }_{{\mathcal{R}}^j}(\xi , \eta ) }{\delta }_{K^{\ast }} \big( X(\xi ,\eta ) \big) = 
\frac{\partial }{\partial z}\rule[-10pt]{.5pt}{24pt}\, 
\raisebox{-9pt}{$\scriptscriptstyle {\delta }_{K^{\ast }}({\Phi }_{{\mathcal{R}}^j}(\xi ,\eta )) 
= R^j z$.} 
\label{eq-s5seven**}
\end{equation}

\noindent \textbf{Proof.} From equation (\ref{eq-s5twostar}) we get
${\delta }_{K^{\ast }} \comp {\Phi }_{\mathcal{R}} = {\Psi }_R \comp {\delta }_{K^{\ast }}$ 
on ${\mathcal{S}}_{\mathrm{reg}}$. Differentiating the preceding equation and then evaluating the result at $X(\xi ,\eta ) \in 
T_{(\xi ,\eta )}{\mathcal{S}}_{\mathrm{reg}}$ gives  
\begin{displaymath}
\big( T_{{\Phi }_{\mathcal{R}}(\xi , \eta ) }{\delta }_{K^{\ast }} \comp T_{(\xi ,\eta )}{\Phi }_{\mathcal{R}} \big) 
X(\xi ,\eta ) = \big( T_{{\delta }_{K^{\ast }} (\xi ,\eta )}{\Psi }_R \comp T_{(\xi ,\eta )} 
{\delta }_{K^{\ast }} \big) X(\xi ,\eta ) 
\end{displaymath}
for all $(\xi ,\eta ) \in {\mathcal{S}}_{\mathrm{reg}}$. When $(\xi ,\eta ) \in \mathcal{D}$, by definition ${\delta }_{K^{\ast }}(\xi ,\eta ) = \delta (\xi ,\eta )$. So for every 
$(\xi ,\eta ) \in {\mathcal{S}}_{\mathrm{reg}}$ 
\begin{displaymath}
T_{(\xi ,\eta )}{\delta }_{K^{\ast }} \big( X(\xi ,\eta ) \big) = 
T_{(\xi ,\eta )} \delta \big( X(\xi ,\eta ) \big) = 
\frac{\partial }{\partial z}\rule[-10pt]{.5pt}{24pt}\, \raisebox{-9pt}{$\scriptscriptstyle z = 
{\delta }(\xi ,\eta ) $} = \frac{\partial }{\partial z}\rule[-10pt]{.5pt}{24pt}\, \raisebox{-9pt}{$\scriptscriptstyle z = {\delta }_{K^{\ast }}(\xi ,\eta ) $}\hspace{-25pt} . 
\end{displaymath}
Thus 
\begin{equation}
T_{{\Phi }_{\mathcal{R}}(\xi , \eta ) } {\delta }_{K^{\ast }} \big( T_{(\xi ,\eta )}{\Phi }_{\mathcal{R}}  
X(\xi ,\eta ) \big) = T_{{\delta }_{K^{\ast }}(\xi ,\eta )}{\Psi }_R \big( \frac{\partial }{\partial z}\rule[-10pt]{.5pt}{24pt}\, \raisebox{-9pt}{$\scriptscriptstyle z = 
{\delta }_{K^{\ast }}(\xi ,\eta ) $} \hspace{-30pt} \big) , 
\label{eq-s5six}
\end{equation}
for every $(\xi ,\eta ) \in \mathcal{D}$. By analytic continuation (\ref{eq-s5six}) holds for every $(\xi ,\eta ) \in {\mathcal{S}}_{\mathrm{reg}}$. Now 
$T_{(\xi , \eta )}{\Phi }_{\mathcal{R}}$ 
sends $T_{(\xi ,\eta )}{\mathcal{S}}_{\mathrm{reg}}$ to $T_{{\Phi }_{\mathcal{R}}(\xi ,\eta )}
{\mathcal{S}}_{\mathrm{reg}}$. Since $T_{(\xi ,\eta )}{\Phi }_{\mathcal{R}}X(\xi ,\eta ) 
= {\mathrm{e}}^{2\pi i/n} X(\xi ,\eta )$ for every 
$(\xi ,\eta )\in {\mathcal{S}}_{\mathrm{reg}}$, 
it follows that ${\mathrm{e}}^{2\pi i/n} X(\xi ,\eta )$ lies in 
$T_{ {\Phi }_{\mathcal{R}}(\xi , \eta ) }{\mathcal{S}}_{\mathrm{reg}}$. Also since 
$T_{{\delta }_{K^{\ast }}(\xi ,\eta )} {\Psi }_R$ sends $T_{{\delta }_{K^{\ast }}(\xi ,\eta ) } K^{\ast }$ to $T_{ {\Psi }_R({\delta }_{K^{\ast }}(\xi , \eta )}K^{\ast }$, we get 
\begin{displaymath}
T_{{\delta }_{K^{\ast } }(\xi ,\eta )} {\Psi }_R 
\big( \frac{\partial }{\partial z}\rule[-10pt]{.5pt}{24pt}\, \raisebox{-9pt}{$\scriptscriptstyle 
z = {\delta }_{K^{\ast }(\xi ,\eta ) } $} \big) = 
R \frac{\partial }{\partial z}\rule[-10pt]{.5pt}{24pt}\, \raisebox{-9pt}{$\scriptscriptstyle 
Rz = {\Psi }_R({\delta }_{K^{\ast }}(\xi ,\eta ) ) $} .
\end{displaymath}
For every $(\xi ,\eta ) \in {\mathcal{S}}_{\mathrm{reg}}$ we obtain 
\begin{equation}
T_{{\Phi }_{\mathcal{R}}(\xi ,\eta ) } 
{\delta }_{K^{\ast }}\big( X(\xi ,\eta ) \big) =  
\frac{\partial }{\partial z}\rule[-10pt]{.5pt}{24pt}\, \raisebox{-9pt}{$\scriptscriptstyle 
Rz = {\Psi }_{\mathcal{R}}({\delta }_{K^{\ast }}(\xi ,\eta ) ) $} ,
\label{eq-s5sevennw}
\end{equation}
that is, equation (\ref{eq-s5seven**}) holds with $j=0$. A similar calculation shows that 
equation (\ref{eq-s5sevennw}) holds with $\mathcal{R}$ replaces by 
${\mathcal{R}}^j$. This verifies equation (\ref{eq-s5seven**}). \hfill $\square $ \medskip 
 
\vspace{-.15in}We now show \medskip 

\noindent \textbf{Theorem 6.7} The image of a $\widehat{\mathcal{G}}$ invariant geodesic on 
$({\mathcal{S}}_{\mathrm{reg}}, \Gamma )$ under the \linebreak 
developing map ${\delta }_{K^{\ast }}$ (\ref{eq-s3three}) is a billiard motion in $K^{\ast }$. \medskip 

\noindent \textbf{Proof.} Because ${\Phi }_{{\mathcal{R}}^j}$ and 
${\Psi }_{R^j}$ are isometries of $({\mathcal{S}}_{\mathrm{reg}}, \Gamma )$ and 
$(K^{\ast }, \gamma |_{K^{\ast }})$, respectively, it follows from equation 
(\ref{eq-s5twostar}) that the surjective map 
${\delta }_{K^{\ast }}: ({\mathcal{S}}_{\mathrm{reg}}, \Gamma ) 
\rightarrow (K^{\ast }, \gamma |_{K^{\ast }} )$ (\ref{eq-s3three}) is an isometry. Hence 
${\delta }_{K^{\ast }}$ is a developing map. Using the local inverse of 
${\delta }_{K^{\ast }}$ 
and equation (\ref{eq-s5seven**}), it follows that a billiard motion in 
$\mathrm{int}(K^{\ast } \setminus \{ 0 \} )$ is mapped onto a geodesic in 
$({\mathcal{S}}_{\mathrm{reg}}, \Gamma )$, which is possibly broken at the 
points $({\xi }_i, {\eta }_i) = {\delta }^{-1}_{K^{\ast}}(p_i)$. Here $p_i \in \partial K^{\ast}$ 
are the points where the billiard motion undergoes a reflection. But the geodesic on 
${\mathcal{S}}_{\mathrm{reg}}$ is smooth at $({\xi }_i, {\eta }_i)$ since the geodesic 
vector field $X$ is holomorphic on ${\mathcal{S}}_{\mathrm{reg}}$. Thus the image of 
the geodesic under the developing map ${\delta }_{K^{\ast }}$ is a billiard motion. 
\hfill $\square $ \medskip 

\vspace{-.1in}\begin{tabular}{l}
\includegraphics[width=300pt]{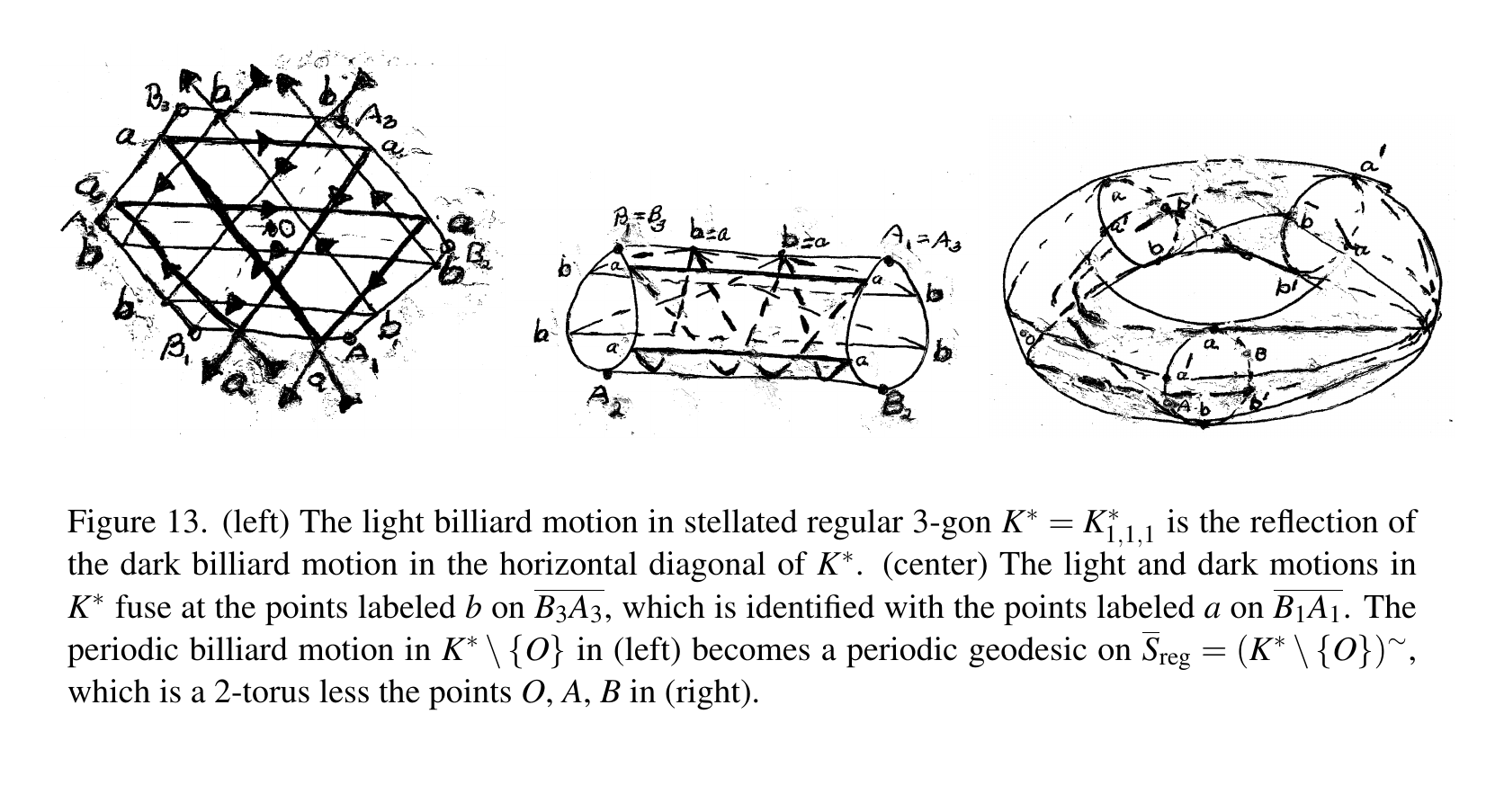}
\end{tabular}    

Next we follow a $G$-invariant set of billiard motions in 
$(K^{\ast },\gamma |_{K^{\ast }})$, which is the union of an $R$-invariant billiard motion and its $U$ reflection. 
After identification of equivalent edges of $\mathrm{cl}(K^{\ast})$, see figure 13 (left) and (center) and then 
dividing out the induced $G$ action, we get a motion on the Riemann surface 
${\widetilde{S}}_{\mathrm{reg}}$, which is a geodesic for the induced Riemannian metric 
$\widehat{\gamma}$ on the $\mathfrak{G}$-orbit space 
$(\C \setminus {\mathbb{V}}^{+})^{\sim}/\mathfrak{G}$, see figure 13 (right). We now justify these assertions. \medskip 

A \emph{billiard} \emph{motion} ${\gamma }_z$ in the regular stellated $n$-gon 
$K^{\ast }$, which starts at a point $z$ in the interior of 
$\mathrm{cl}(K^{\ast })\setminus \{ O \} $ and does not hit a vertex of 
$\mathrm{cl}(K^{\ast })$, is made up of line segments, 
each of which is parallel to an edge of $\mathrm{cl}(K^{\ast }$). It is invariant under 
the subgroup $\widehat{G}$ of $G$ generated by the rotation $R$. 
Let ${\mathrm{Rfl}}^{{\gamma }_z}  = 
\{ p \in \partial \, \mathrm{cl}(K^{\ast}) \setrule \, p = {\gamma }_z(T_p) \, \, 
\mbox{for some $T_p \in \R $} \} $ be the set of reflection points in the boundary of 
$\mathrm{cl}(K^{\ast })$ of the billiard motion ${\gamma }_z$. Since 
${\gamma }_z$ is invariant under the group $\widehat{G}$, the set 
${\mathrm{Rfl}}^{{\gamma }_z}$ of reflection points is invariant under group 
$\widehat{G}$. Because ${\gamma }_z$ does not hit a vertex of 
$\mathrm{cl}(K^{\ast })$, $z$ is not fixed by the reflection $U$. 
The billiard motion ${\gamma }_{\, \overline{z}}$ starting at $\overline{z} = U(z)$ is invariant under the group $\widehat{G}$ and, by uniqueness of billiard motions with a given starting point, is equal to the billiard motion $U({\gamma }_z) = 
{\gamma }_{\, \overline{z}}$. So $U({\mathrm{Rfl}}^{{\gamma }_z}) = 
{\mathrm{Rfl}}^{{\gamma }_{\, \overline{z}}}$. From $U(z) \ne z$, it follows that 
${\mathrm{Rfl}}^{{\gamma }_z} \cap {\mathrm{Rfl}}^{{\gamma }_{\, \overline{z}}} = 
\varnothing $. Let $E^{{\gamma }_z}$ be the set of closed edges of 
$\mathrm{cl}(K^{\ast })$, which the billiard motion 
${\gamma }_z$ reflects off of. In other words, $E^{{\gamma }_z} = 
\{ \mbox{$E$ an edge of $\mathrm{cl}(K^{\ast })$} 
\setrule \, p \in E \, \, \mbox{for some $p \in {\mathrm{Rfl}}^{{\gamma }_z}$} \} $. 
\medskip 

\noindent \textbf{Lemma 6.8} $E^{{\gamma }_{\overline{z}}} = 
U(E^{{\gamma }_z})$. \medskip

\noindent \textbf{Proof.} Suppose that $E \in E^{{\gamma }_{\overline{z}}}$. 
Then for some $p \in {\mathrm{Rfl}}^{{\gamma }_{\overline{z}}}$ we have 
$p \in E$. Since ${\mathrm{Rfl}}^{{\gamma }_{\overline{z}}} = 
U({\mathrm{Rfl}}^{{\gamma }_z})$, $U(p) \in U({\mathrm{Rfl}}^{{\gamma }_{\overline{z}}}) = 
{\mathrm{Rfl}}^{{\gamma }_z}$ and $U(p) \in U(E)$. Thus $U(E) \in E^{{\gamma }_z}$. So $U(E^{{\gamma }_{\overline{z}}}) \subseteq E^{{\gamma }_z}$. A similar argument 
shows that $U(E^{{\gamma }_z}) \subseteq E^{{\gamma }_{\overline{z}}}$. 
Hence $E^{{\gamma }_{\overline{z}}} = U ( U(E^{{\gamma }_{\overline{z}}}))
\subseteq U(E^{{\gamma }_z}) \subseteq E^{{\gamma }_{\overline{z}}}$, which 
implies $E^{{\gamma }_{\overline{z}}} = U(E^{{\gamma }_z})$. \hfill $\square $ \medskip  

\noindent \textbf{Lemma 6.9} The sets $E^{{\gamma }_z}$ and 
$E^{{\gamma }_{\overline{z}}}$ are $\widehat{G}$-invariant. 
\medskip 

\noindent \textbf{Proof.} Let $E \in E^{{\gamma }_z}$ and $p \in E \cap 
{\mathrm{Rfl}}^{{\gamma }_z}$. Since ${\mathrm{Rfl}}^{{\gamma }_z}$ is 
$\widehat{G}$-invariant, it follows that $R(p) \in {\mathrm{Rfl}}^{{\gamma }_z}$ and 
$R(p) \in R(E)$. Hence $R(E) \in E^{{\gamma }_z}$. So $E^{{\gamma }_z}$ is 
$\widehat{G}$-invariant. Similarly, $E^{{\gamma }_{\overline{z}}}$ is 
$\widehat{G}$-invariant. \hfill $\square $ \medskip 

\noindent \textbf{Lemma 6.10} Let $S_0$ be the reflection $R^{n_0}U$ and set 
$S_m = R^mS_0R^{-m}$ for $m \in \{ 0,1, \ldots , n-1 \} $. Then 
$S_m({\mathrm{Rfl}}^{{\gamma }_z}) = U({\mathrm{Rfl}}^{{\gamma }_z})$. \medskip 

\noindent \textbf{Proof.} If $p \in {\mathrm{Rfl}}^{{\gamma }_z}$, then 
$S_m(p) \in U({\mathrm{Rfl}}^{{\gamma }_z})$, for $U(p) \in 
U({\mathrm{Rfl}}^{{\gamma }_z})$, which implies $R^{n_0}((U(p))) \in 
U({\mathrm{Rfl}}^{{\gamma }_z})$, since $U({\mathrm{Rfl}}^{{\gamma }_z})$ is 
$\widehat{G}$-invariant. Hence $S_0(p) \in U({\mathrm{Rfl}}^{{\gamma }_z})$. 
If $p \in {\mathrm{Rfl}}^{{\gamma }_z}$, then 
$R^{-m}(p) \in {\mathrm{Rfl}}^{{\gamma }_z}$, 
since ${\mathrm{Rfl}}^{{\gamma }_z}$ is $\widehat{G}$-invariant. So 
$S_0(R^{-m}(p)) \in  U( {\mathrm{Rfl}}^{{\gamma }_z} )$, which implies 
$R^mS_0(R^{-m}(p)) \in U( {\mathrm{Rfl}}^{{\gamma }_z} )$, because 
$U( {\mathrm{Rfl}}^{{\gamma }_z} )$ is $\widehat{G}$-invariant. So 
$S_m({\mathrm{Rfl}}^{{\gamma }_z}) \subseteq U( {\mathrm{Rfl}}^{{\gamma }_z} )$. 
A similar argument shows that $S_m(U( {\mathrm{Rfl}}^{{\gamma }_z} )) \subseteq 
{\mathrm{Rfl}}^{{\gamma }_z}$. Thus 
\begin{displaymath}
{\mathrm{Rfl}}^{{\gamma }_z} = S_m(S_m({\mathrm{Rfl}}^{{\gamma }_z})) \subseteq 
S_m(U( {\mathrm{Rfl}}^{{\gamma }_z} )) \subseteq {\mathrm{Rfl}}^{{\gamma }_z} .
\end{displaymath}
So $S_m(U( {\mathrm{Rfl}}^{{\gamma }_z} )) = {\mathrm{Rfl}}^{{\gamma }_z}$, that is, 
$U( {\mathrm{Rfl}}^{{\gamma }_z} ) = S_m( {\mathrm{Rfl}}^{{\gamma }_z}) $.  
\hfill $\square $ \medskip 

\noindent \textbf{Lemma 6.11} Every reflection $S_m$ interchanges an edge in $E^{{\gamma }_z}$ with 
an edge in $E^{{\gamma }_{\, \overline{z}}}$, specifically, $S_m(E^{{\gamma }_z}) = E^{{\gamma }_{\overline{z}}}$. \medskip 

\noindent \textbf{Proof.} Let $E \in E^{{\gamma }_z}$. Then there is a 
$p \in {\mathrm{Rfl}}^{{\gamma }_z}$ such that $p \in E$. So $S_m(p) \in S_m(E)$. 
But $S_m(p) \in U({\mathrm{Rfl}}^{{\gamma }_z})$, which shows that 
$S_m(E) \in U(E^{{\gamma }_z})$. Hence $S_m(E^{{\gamma }_z}) 
\subseteq U(E^{{\gamma }_z})$. A similar argument shows that 
$S_m(U(E^{{\gamma }_z})) \subseteq E^{{\gamma }_z}$. Thus 
$E^{{\gamma }_z} = S_m(S_m(E^{{\gamma }_z})) \subseteq S_m(U(E^{{\gamma }_z})) 
\subseteq E^{{\gamma }_z}$. So $S_m(U(E^{{\gamma }_z})) = E^{{\gamma }_z}$, which implies $S_m(E^{{\gamma }_z}) = U(E^{{\gamma }_z}) = 
E^{{\gamma }_{\overline{z}}}$. \hfill $\square $ \medskip %

An \emph{extended billiard motion} ${\lambda }_z$ in $K^{\ast }$ starting at a point $z \in 
\mathrm{int}(K^{\ast } \setminus \{ 0 \} )$ is the union of a billiard motion ${\gamma }_z$ 
in $(\mathrm{int}\, K^{\ast }) \setminus \{ O \}$ starting at $z$ and a billiard motion 
${\gamma }_{\, \overline{z}}$ in $(\mathrm{int}\, K^{\ast }) \setminus \{ O \} $ starting at 
$\overline{z} = Uz$. The motion ${\lambda }_z$ is invariant under the group generated by the rotation $R$ and the reflection $U$. So ${\lambda }_z$ is $G$-invariant. The set of points of an extended billiard motion in $K^{\ast } \setminus \{ O \}$, which lie on 
$\partial K^{\ast}$ is $G$-invariant and is the disjoint union of reflection points 
${\mathrm{Rfl}}^{{\gamma }_z}$ for the billiard motion ${\gamma }_z$ 
and ${\mathrm{Rfl}}^{{\gamma }_{\, \overline{z}}} = U({\mathrm{Rfl}}^{{\gamma }_z})$ 
for its $U$ reflection ${\gamma }_{\, \overline{z}}$. From lemma 6.10 it follows that the equivalence relation $\sim $ among the closed edges of $\mathrm{cl}(K^{\ast })$ interchanges these subsets. Identifying equivalent points in 
${\mathrm{Rfl}}^{{\gamma }_z}$ and ${\mathrm{Rfl}}^{{\gamma }_{\, \overline{z}}}$ with 
the equivalent edges, in which they are contained, gives a \emph{continuous} motion 
${\lambda }^{\sim}_z = \Pi ({\lambda }_z)$ in the smooth space 
$(K^{\ast } \setminus \{ O \})^{\sim}$, which is $G$-invariant. Here $\Pi $ is the map 
(\ref{eq-s4sevenstar}). \medskip 

\noindent \textbf{Theorem 6.12} Under the restriction of the mapping 
\begin{equation}
\nu = \sigma \comp \Pi : \C \setminus {\mathbb{V}}^{+}   \rightarrow  (\C \setminus {\mathbb{V}}^{+} )^{\sim }/\mathfrak{G} = {\widetilde{S}}_{\mathrm{reg}}
\label{eq-s5sevenverynw}
\end{equation} 
to $K^{\ast }\setminus \{ O \} $ the image of an extended billiard motion 
${\lambda }_z$ is a smooth geodesic ${\widehat{\lambda }}_{\nu (z)}$ on 
$({\widetilde{S}}_{\mathrm{reg}} , \widehat{\gamma } ) $, where 
${\nu }^{\ast }(\widehat{\gamma }) = \gamma |_{\C \setminus {\mathbb{V}}^{+}}$. \medskip 

\noindent \textbf{Proof.} Since the Riemannian metric 
$\gamma $ on $\C $ is 
invariant under the group of Euclidean motions, the Riemannian metric 
${\gamma }|_{K^{\ast } \setminus \{ O \}}$ on $K^{\ast } \setminus \{ O \}$ is 
$G$-invariant. Hence ${\gamma }_{K^{\ast } \setminus \{ O \} }$ is invariant 
under the reflection $S_m$ for $m \in \{ 0,1, \ldots , n-1 \} $. So 
${\gamma }|_{K^{\ast } \setminus \{ O \} }$ pieces together to give a Riemannian 
metric ${\gamma }^{\sim }$ on the identification space 
$(K^{\ast } \setminus \{ O \} )^{\sim}$. In other words, the pull back of 
${\gamma }^{\sim }$ under the map ${\Pi }|_{K^{\ast } \setminus \{ O \} }: 
K^{\ast } \setminus \{ O \} \rightarrow (K^{\ast } \setminus \{ O \})^{\sim }$, 
which identifies equivalent edges of $K^{\ast }$, is the metric 
${\gamma }|_{K^{\ast } \setminus \{ O \} }$. Since ${\Pi }|_{K^{\ast } \setminus \{ O \} }$ 
intertwines the $G$-action on $K^{\ast } \setminus \{ O \}$ with the $G$-action on 
$(K^{\ast } \setminus \{ O \})^{\sim }$, the metric ${\gamma }^{\sim }$ is $G$-invariant. 
It is flat because the metric $\gamma $ is flat. So ${\gamma }^{\sim }$  
induces a flat Riemannian metric $\widehat{\gamma }$ 
on the orbit space $(K^{\ast } \setminus \{ O \} )^{\sim } /G = 
{\widetilde{S}}_{\mathrm{reg}}$. Since the extended billiard motion ${\lambda}_z$ is 
a $G$-invariant broken geodesic on $(K^{\ast } \setminus \{ O \}, 
{\gamma }_{K^{\ast } \setminus \{ O \} })$, which is made up of two continuous pieces, 
it gives rise to a \emph{continuous} broken geodesic 
${\lambda }^{\sim }_{\Pi (z)}$ on $((K^{\ast } \setminus \{ O \})^{\sim }, {\gamma }^{\sim })$, which is $G$-invariant. Thus ${\widehat{\lambda }}_{\nu (z) } = \nu ({\lambda}_z)$ is a piecewise smooth geodesic on the smooth $G$-orbit space 
$((K^{\ast } \setminus \{ O \})^{\sim }/G = 
{\widetilde{S}}_{\mathrm{reg}}, \widehat{\gamma })$.  
\medskip %

We need only show that ${\widehat{\lambda }}_{\nu (z)}$ is smooth. To see this we 
argue as follows. Let $s \subseteq K^{\ast }$ be a closed 
segment of a billiard motion ${\gamma }_z$, which is contained in the extended 
billiard motion ${\lambda }_z$ that does not meet a vertex of $\mathrm{cl}(K^{\ast })$. 
Then ${\gamma }_z$ is a horizontal straight line motion in $\mathrm{cl}(K^{\ast })$. 
Suppose that $E_{k_0}$ is the edge of $K^{\ast }$, perpendicular to the 
direction $u_{k_0}$, which is first met by ${\gamma }_z$ and let 
$P_{k_0}$ be the meeting point. Let $S_{k_0}$ be the reflection in $E_{k_0}$. 
The continuation of the motion ${\gamma }_z$ at $P_{k_0}$ is the 
horizontal line $RS_{k_0}({\gamma }_z)$ in $K^{\ast }_{k_0}$. 
Recall that $K^{\ast }$ is the translation 
of $K^{\ast }$ by ${\tau }_{k_0}$. Since $O_{k_0} = {\tau }_{k_0}(0)$ is the center of 
$K^{\ast }_{k_0}$, the extended motion is the same as the motion $U({\gamma }_z)$
translated by ${\tau }_k$. Using a suitable sequence of \linebreak 

\mbox{}\vspace{-.3in}\begin{tabular}{l}
\includegraphics[width=300pt]{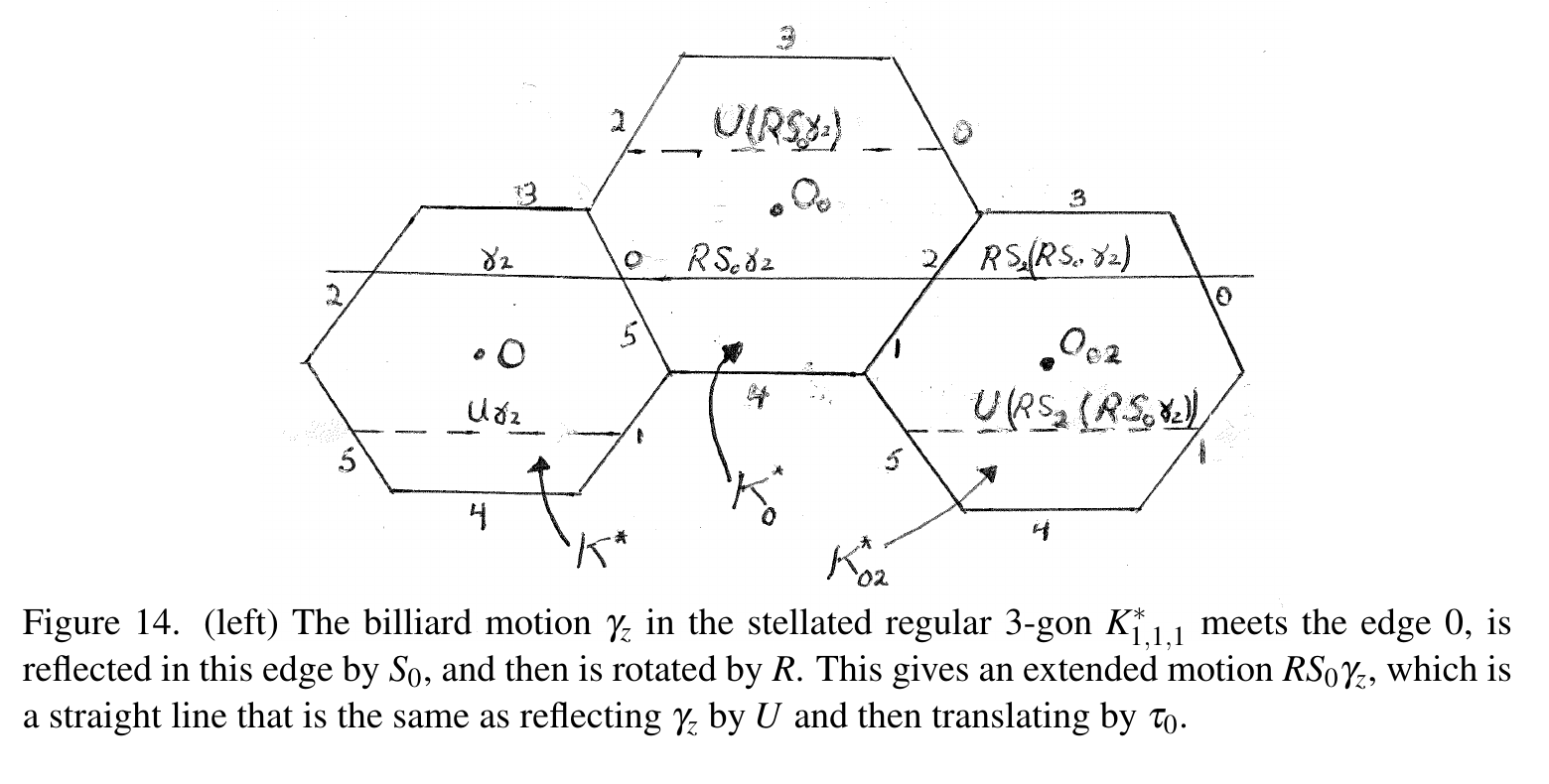}
\end{tabular}
\vspace{.25in}\par \noindent reflections in the edges of a suitable 
$K^{\ast }_{k_0 \cdots k_{\ell }}$ followed by a rotation $R$, which gives rise to a reflection $U$ and a translation in $\mathcal{T}$ corresponding to their origins, we can extend 
$s$ to a smooth straight line $\lambda $ in $\C \setminus {\mathbb{V}}^{+}$, see figure 14. The line $\lambda $ is a geodesic in $(\C \setminus {\mathbb{V}}^{+}, \gamma |_{\C \setminus {\mathbb{V}}^{+}})$, 
which in $K^{\ast }$ has image ${\widehat{\lambda }}_{\nu (z)}$ under the $\mathfrak{G}$-orbit map that is a smooth geodesic on 
$({\widehat{S}}_{\mathrm{reg}}, \widehat{\gamma } ) $. 
The geodesic $\nu (\lambda )$ starts at $\nu (z)$. Thus the smooth geodesic 
${\widehat{\lambda}}_{\nu (z)}$ and the possibly broken geodesic $\nu (\lambda )$ are equal. In other words, $\nu (\lambda )$ is a smooth geodesic. \hfill $\square $ \medskip 

Thus the affine orbit space ${\widetilde{S}}_{\mathrm{reg}} = 
(\C \setminus {\mathbb{V}}^{+}) /\mathfrak{G}$ with flat Riemannian metric 
$\widehat{\gamma }$ is the \emph{affine} analogue of the Poincar\'{e} model of 
the affine Riemann surface ${\mathcal{S}}_{\mathrm{reg}}$ as an orbit space of a discrete subgroup of $\mathrm{PGl}(2, \C)$ acting on the unit disk in $\C $ with the Poincar\'{e} metric. 

\section{Appendix. Group theoretic properties}

In this appendix we discuss some group theoretic properties of the set of equivalent edges of 
$\mathrm{cl}(K^{\ast })$, which we use to determine the topology of 
${\widetilde{S}}_{\mathrm{reg}}$. \medskip 

Let $\mathcal{E}$ be  the set of unordered pairs $[E,E']$ of nonadjacent edges of $\mathrm{cl}(K^{\ast })$. Define an action 
\hspace{-2pt}\raisebox{1pt}{ \tiny $\bullet $} of $G$ on $\mathcal{E}$ by 
\begin{displaymath}
g \mbox{\hspace{-2pt}\raisebox{1pt}{ \tiny $\bullet $}} [E, E'] = [g(E), g(E')] 
\end{displaymath}
for every unordered pair $[E, E']$ of nonadjacent edges of $\mathrm{cl}(K^{\ast})$. For every $g \in G$ the edges $g(E)$ and $g(E')$ are nonadjacent. This follows 
because the edges $E$ and $E'$ are nonadjacent and the elements of $G$ are invertible mappings of $\C $ into itself. So $\varnothing = g (E \cap E') = 
g(E) \cap g(E')$. Thus the mapping \hspace{-2pt}\raisebox{1pt}{\tiny $\bullet $} is well defined. It is an action because for every $g$ and $h \in G$ we have 
\begin{align}
g \mbox{\hspace{-2pt}\raisebox{1pt}{ \tiny $\bullet $}} 
(h \mbox{\hspace{-2pt}\raisebox{1pt}{ \tiny $\bullet $}} [E, E']) & = 
g \mbox{\hspace{-2pt}\raisebox{1pt}{ \tiny $\bullet $}} [h(E), h(E')] = [g(h(E), g(h(E')] \notag \\
& = [(gh)(E), (gh)(E') ] = 
(gh) \mbox{\hspace{-2pt}\raisebox{1pt}{ \tiny $\bullet $}} [E, E']. \notag 
\end{align}
The action \hspace{-2pt}\raisebox{1pt}{ \tiny $\bullet $} of $G$ on $\mathcal{E}$ induces an action $\cdot $ of the group $G^j$ of reflections on the set ${\mathcal{E}}^j$ of equivalent edges of $\mathrm{cl}(K^{\ast })$, which is defined by 
\begin{displaymath}
g_j \cdot [E, S^{(j)}_k(E)] = [g_j(E), g_j(S^{(j)}_k(E))] = 
[g_j(E), (g_jS^{{(j)}}_kg^{-1}_j)(g_j(E))] ,
\end{displaymath}
for every $g_j \in G^j$, every edge $E$ of $\mathrm{cl}(K^{\ast })$, and every generator $S^{(j)}_k$ of $G^j$, where $k = 0, 1, \ldots ,$  
$n-1$. Since $g_jS^{(j)}_kg^{-1}_j = S^{(j)}_r$ by corollary 3.3b, the mapping $\cdot $ 
is well defined. \medskip 

\noindent \textbf{Lemma A1} The group $G$ action 
\raisebox{1pt}{\tiny $\bullet $} sends a 
$G^j$-orbit on ${\mathcal{E}}^j$ to another $G^j$-orbit on 
${\mathcal{E}}^j$. \medskip %

\noindent \textbf{Proof.} Consider the $G^j$-orbit of 
$[E, S^{(j)}_m(E)] \in {\mathcal{E}}^j$. For every $g \in G$ we have 
\begin{displaymath}
g \, \, \mbox{\hspace{-2pt}\raisebox{1pt}{\tiny $\bullet $}} 
\big( G^j \cdot [E, S^{(j)}_m(E)] \big) = 
(g G^j g^{-1}) \, \, \mbox{\hspace{-2pt}\raisebox{1pt}{\tiny $\bullet $}} 
\big( g \, \mbox{\raisebox{1pt}{\tiny $\bullet $}}\,  [E, S^{(j)}_m(E)] \big) = 
G^j \cdot \big( g \, \, \mbox{\hspace{-2pt}\raisebox{1pt}{\tiny $\bullet $}}  
\, [E, S^{(j)}_m(E)] \big) , 
\end{displaymath}
because $G^j$ is a normal subgroup of $G$ by corollary 3.3c. Since 
\begin{displaymath}
g \, \, \mbox{\hspace{-2pt}\raisebox{1pt}{\tiny $\bullet $}} \, [E, S^{(j)}_m(E)] = 
[g(E), g(S^{(j)}_m(E))] = [g(E), gS^{(j)}_m g^{-1}(g(E))] 
\end{displaymath}
and $gS^{(j)}_m g^{-1} =S^{(j)}_r$ by corollary 3.3b, it follows that 
$g \, \, \mbox{\hspace{-2pt}\raisebox{1pt}{\tiny $\bullet $}}\, [E, S^{(j)}_m(E)] \in 
{\mathcal{E}}^j$. \hfill $\square $ \medskip 

\noindent \textbf{Lemma A2} For every $j=0,1,\infty $ and every 
$k=0,1, \ldots , n-1$ the isotropy group $G^j_{e^j_k}$ of the
$G^j$ action on ${\mathcal{E}}^j$ at $e^j_k = [E, S^{(j)}_k(E)]$ is 
$\langle S^{(j)}_k \setrule \, (S^{(j)}_k)^2 =e \rangle $. 
\medskip 

\noindent \textbf{Proof.} Every $g \in G^j_{e^j_k}$ satisfies 
\begin{displaymath}
e^j_k = [E, S^{(j)}_k(E)] = g \cdot e^j_k = g \cdot [E, S^{(j)}_k(E)]
\end{displaymath} 
if and only if 
\begin{displaymath}
[E, S^{(j)}_k(E)] = [g(E), gS^{(j)}_kg^{-1}(g(E))] = [g(E), S^{(j)}_r(g(E))]
\end{displaymath}
if and only if one of the statements 1) $g(E) =E$ \& $S^{(j)}_k(E) = S^{(j)}_r(g(E))$ or 
2) $E = g(S^{(j)}_r(E))$ \& $g(E) = S^{(j)}_k(E)$ holds. From $g(E) = E$ in 1) we get 
$g = e$ using lemma 3.2. To see this we argue as follows. If $g \ne e $, then 
$g = R^p(S^{(j)})^{\ell }$ for some $\ell = 0 ,1$ and some 
$p \in \{ 0,1, \ldots , n -1 \} $, see equation (\ref{eq-s3threedotdagger}). 
Suppose that $g = R^p$ with $p \ne 0$. Then $g(E) \ne E$, which contradicts our 
hypothesis. Now suppose that $g = R^pS^{(j)}$. Then $E = g(E) = R^pS^{(j)}(E)$, which gives $R^{-p}(E) = S^{(j)}(E)$. Let $A$ and $B$ be end points of the edge $E$. Then 
the reflection $S^{(j)}$ sends $A$ to $B$ and $B$ to $A$, while the rotation $R^{-p}$ 
sends $A$ to $A$ and $B$ to $B$. Thus $R^{-p}(E) \ne S^{(j)}(E)$, which is a contradiction. Hence $g =e$. If $g(E) = S^{(j)}_k(E)$ in 2), then $(S^{(j)}_kg)(E) = E$. So $S^{(j)}_kg = e$ by lemma 3.2, that is, $g = S^{(j)}_k$. \hfill $\square $ \medskip 

For every $j=0,1, \infty $ and every $m = 0,1, \ldots , \frac{n}{d_j}-1$ let 
$G^j_{e^j_{md_j}} = 
\{ g_j \in G^j \setrule \, g_j \cdot e^j_{md_j} = e^j_{md_j} \} $ be the isotropy 
group of the $G^j$ action on ${\mathcal{E}}^j$ at 
$e^j_{md_j} = [E, S^{(j)}_{md_j}(E)]$. Since $G^j_{e^j_{md_j}} = 
\langle S^{(j)}_{md_j} \setrule \, (S^{(j)}_{md_j})^2 = e \rangle $ is an abelian subgroup of $G^j$, it is a normal subgroup. Thus $H^j = 
G^j/G^j_{e^j_{md_j}}$ is a subgroup of $G^j$ of 
order $(2n/d_j)/2 = n/d_j$. This proves \medskip 

\noindent \textbf{Lemma A3} For every $j=0,1, \infty$ and each 
$m = 0,1, \ldots , \frac{n}{d_j}-1$ the 
$G^j$-orbit of $e^j_{md_j}$ in ${\mathcal{E}}^j$ is equal to the $H^j$-orbit of $e^j_{md_j}$ 
in ${\mathcal{E}}^j$. \medskip 

\noindent \textbf{Lemma A4} For $j=0,1, \infty $ we have $H^j = \langle V = R^{d_j} 
\setrule \, V^{n/d_j} = e \rangle $. \medskip 

\noindent \textbf{Proof.} Since 
\begin{equation}
S^{(j)}_k = R^kS^{(j)}R^{-k} = R^k(R^{n_j}U)R^{-k} = R^{2k+n_j}U = R^{2k}S^{(j)}, 
\label{eq-s3threedotdagger}
\end{equation}
we get $S^{(j)}_{md_j} = R^{(2m +\frac{n_j}{d_j})d_j} U = (R^{d_j})^mS^{(j)}$. 
Because the group $G^j$ is generated by the reflections $S^{(j)}_k$ for 
$k = 0,1, \ldots , n-1$, it follows that 
\begin{displaymath}
G^j \subseteq \langle V=R^{d_j}, S^{(j)}_{md_j} \setrule \, V^{n/d_j} = e = 
(S^{(j)}_{md_j})^2 \, \, \& \, \, VS^{(j)}_{md_j} = S^{(j)}_{md_j}V^{-1} \rangle = K_j. 
\end{displaymath}
$K_j$ is a subgroup of $G$ of order $2n/d_j$. Clearly the isotropy group 
$G^j_{e^j_{md_j}} = \langle S^{(j)}_{md_j} \setrule \, (S^{(j)}_{md_j})^2 =e \rangle $ 
is an abelian subgroup of $K^j$. Hence $H^j = G^j/G^j_{e^j_{md_j}} 
\subseteq K^j/G^j_{e^j_{md_j}} = L^j$, where $L^j$ is a subgroup of $K^j$ of order 
$(2n/d_j)/2 = n/d_j$. Thus the group $L^j$ has the same order as its subgroup 
$H^j$. So $H^j= L^j$. But $L^j = \langle V=R^{d_j} \setrule V^{n/d_j} = e \rangle $. 
\hfill $\square $ \medskip  

Let $f^j_{\ell } = R^{\ell } \cdot e^j_0$. Then 
\begin{align}
f^j_{\ell } & = R^{\ell }\cdot e^j_0 = R^{\ell } \cdot [E, S^{(j)}(E)] \notag \\
& = [R^{\ell }(E), R^{\ell }S^{(j)}R^{-\ell } (R^{\ell }(E))] = 
[R^{\ell }(E), S^{(j)}_{\ell }(R^{\ell }(E))] . \notag 
\end{align}
So 
\begin{align}
V^m \cdot f^j_{\ell } & = V^m \cdot [R^{\ell }(E), R^{\ell }S^{(j)}R^{-\ell } (R^{\ell }(E))] \notag \\
& =[V^m(R^{\ell }(E)), V^mS^{(j)}_{\ell }V^{-m}(V^m(R^{\ell }(E))] \notag \\
& = [R^{md_j + \ell }(E), S^{(j)}_{md_j +\ell }(E) ] = 
e^j_{md_j +\ell }.  \notag 
\end{align}
This proves  
\begin{equation}
\bigcup^{d_j-1}_{\ell = 0 } H^j \cdot f^j_{\ell } = 
\bigcup^{d_j-1}_{\ell =0} \bigcup^{\frac{n}{d_j}-1}_{m = 0} V^m \cdot f^j_{\ell } 
= \bigcup^{n-1}_{k=0} e^j_k, 
\label{eq-s3threedoublestar} 
\end{equation}
since every $k \in \{ 0,1, \ldots , n-1 \} $ may be written uniquely as 
$md_j +\ell $ for some $m \in \{ 0,1, \ldots , \frac{n}{d_j}-1 \} $ and some 
$\ell \in \{ 0,1, \ldots , d_j-1 \}$.

\end{document}